\documentclass[11pt]{article}
\usepackage{verbatim}
\usepackage{hyperref}
\hypersetup{
	colorlinks=true,
	linkcolor=red,
	filecolor=magenta,      
	urlcolor=cyan,
}
\usepackage[abs]{overpic}
\usepackage{subcaption}
\RequirePackage{amsmath,amsfonts,amssymb}
\RequirePackage{rotating}
\RequirePackage{graphicx,xcolor}

\newtheorem{rem}{Remark}

\newcommand{\bA}{\mathbf A}

\newcommand{\bk}{\mathbf k}

\newcommand{\bu}{\mathbf u}

\newcommand{\bv}{\mathbf v}

\def\cl {\nonumber \\}
\def\el {\nonumber }

\begin{document}

\title{Filter stabilization for the mildly compressible Euler equations with application to atmosphere dynamics simulations}
\author{Nicola Clinco$^\star$, Michele Girfoglio$^\star$, Annalisa Quaini$^{\star \star}$, Gianluigi Rozza$^\star$}

\maketitle

\begin{center}
\footnotesize
$\star$ mathLab, Mathematics Area, SISSA, via Bonomea, 265, Trieste, I-34136, Italy \\ $\star \star$ Department of Mathematics, University of Houston, Houston TX 77204, USA
\end{center}

\begin{abstract}
We present a filter stabilization technique for the mildly compressible Euler equations that relies on a linear or nonlinear indicator function to identify the regions of the domain where artificial viscosity is needed and determine its amount. For the realization of this technique, we adopt a three step algorithm called Evolve-Filter-Relax (EFR), which at every time step evolves the solution (i.e., solves the Euler equations on a coarse mesh), then filters the computed solution, and finally performs a relaxation step to combine the filtered and non-filtered solutions. We show that the EFR algorithm is equivalent to an eddy-viscosity model in Large Eddy Simulation. Three indicator functions are considered: a constant function (leading to a linear filter), a function proportional to the norm of the velocity gradient (recovering a Smagorinsky-like model), and a function based on approximate deconvolution operators. Through well-known benchmarks for atmospheric flow, we show that the  deconvolution-based filter yields stable solutions that are much less dissipative than the linear filter and the Samgorinsky-like model and we highlight the efficiency of the EFR algorithm. 
\end{abstract}


\section{Introduction}

The Direct Numerical Simulation (DNS) is a simulation in Computational Fluid Dynamics (CFD) that solves the equations governing the fluid motion by resolving the entire range of relevant spatial and temporal scales. In many practical CFD applications, the smallest spatial scales can be several orders of magnitude smaller than the largest scales in the flow. 
An example is atmospheric flow, whose smallest spatial scales are typically of the order of $10^{-4}$ m while the typical domain size is of the order of $10^4-10^5$ m.
For these applications, a DNS is beyond reach for nowadays computing machines and it will be for the foreseeable future. 

One way to keep the computational cost affordable without sacrificing accuracy is to solve for the flow using a coarser mesh and model the effects of the small scales that are not directly solved through a so-called subgrid-scale (SGS) model.
This is the principal idea behind Large Eddy Simulation (LES). Traditionally, SGS models introduce the effects of the unresolved scales with momentum fluxes that are linearly dependent upon the rate of strain of the large scales. This is known eddy-viscosity closure. The most famous eddy-viscosity model is the Smagorinsky model \cite{smagorinsky1963}. Its success is due to several factors: i) it is relatively simple and easy to implement, ii) it is computationally inexpensive compared to other SGS models, and iii) it features parameters that can be tuned for the particular application at hand so that the results are realistic. The main limitation of the Smagorinsky model is the assumption of local balance between the subgrid scale energy production and dissipation. Since such equilibrium conditions do not hold in many practical applications, the Smagorinsky model often results into over-diffusive simulations. A large body of research has been motivated by improving upon the Smagorinsky model.

Some alternative methods introduce artificial diffusion that can be solution-dependent (see, e.g., \cite{ABGRALL2001277,hughes1995,klockner_warburton_hesthaven_2011,rispoliSaavedra2006,perssonPeraire2006}) or residual-based (see, e.g., \cite{Guermond_pasquetti_popov_JCP_2011,guermondPasquetti2008,Guermond_Popov_2014,KURGANOV20128114,marrasNazarovGiraldo2015}). These methods are driven by the intent to have an artificial viscosity that vanishes where the solution is smooth and/or decreases as the grid is refined. 
Other methods add a set of equations to the discrete governing equations formulated on a coarse mesh (coarse is meant with respect to the resolution required by DNS).
This extra-problem can be devised in different ways, for example by a functional
splitting of the solved and unresolved scales
as in variational multiscale methods (see, e.g., \cite{bazilevsCaloCottrellHughesRealiScovazzi2007, codina2002,Codinaetal2017,hughesFeijoo1998}). In this paper, we propose an extra
problem that acts as a differential (linear or nonlinear) low-pass filter added sequentially to the mildly compressible Euler equations for stratified flows. This sequential algorithm is called Evolve-Filter-Relax (EFR) since, at every time step, one first evolves the solution, i.e., solves the Euler equations on a coarse mesh, then filters the computed solution, and finally performs a relaxation step to combine the filtered and non-filtered solutions. 
This techniques is also known with the name of filter stabilization because it reduces or eliminates unphysical fluctuations in the computed solution. We will show that the EFR algorithm is an eddy-viscosity model.

Introduced in \cite{layton_JMFM}, EFR algorithms
have been widely applied to the incompressible Navier-Stokes equations \cite{BQV,Bowers2012,abigail_CMAME,Ervin2012,layton_JMFM,GIRFOGLIO201927,LAYTON20113183,Olshanskii2013}. It was shown in \cite{Bowers2012,abigail_CMAME,Ervin2012,layton_JMFM,LAYTON20113183} that numerical results obtained with nonlinear differential filters are more precise in localizing where eddy viscosity is needed and are overall more accurate than results obtained with plain Smagorinsky-type models or variational multiscale methods. Despite these promising results for the incompressible Navier-Stokes equations, the application of filter stabilization to the Euler equations has received much less attention \cite{Chehab2021,Hesthaven_Warburton, Euler_comp_orig, Euler_comp}. In this paper, we consider the EFR algorithm for the Euler equations with both linear and nonlinear filters. 
Developed in \cite{CNM:CNM219,Boyd1998283,Fischer2001265}, stabilization based on linear filters has been widely studied (see, e.g., \cite{B-garnier,Mathew2003,Visbal2002}). However, it was noted in \cite{doi:10.1137/120867482} that a linear stabilization can, at most, give a solution converging to a weak solution that is not the entropy solution, hence the need to investigate nonlinear filter stabilization techniques. 
We will show that the EFR algorithm with a deconvolution-based filter yields stable solutions that are much less dissipative than the Smagorinsky model. We recall that the use of deconvolution operators in SGS models to increase accuracy is well established and
mathematically grounded  \cite{Dunca2005,Stolz1999,Stolz2001}.

The main advantages of the EFR algorithm are: i) modularity, i.e., its implementation does not require any major modification of a legacy solver, and ii) flexibility in the choice of the filter. In addition, if one chooses the 
deconvolution-based filter, the
viscosity introduced by the EFR algorithm vanishes where the solution is smooth and decreases as the mesh is refined.  We will show that we obtain numerical results that agree very well with data published in the literature for  
well-known 2D benchmark problems involving stratified and gravity driven atmospheres. We will also show that the computational cost to solve the additional filter problem is a fraction of the computational cost required by the Euler solver. 

All the simulations in this paper have been carried out with GEA (Geophysical and Environmental Applications) \cite{GEA}, a new open-source atmosphere and ocean modeling framework within the finite volume C++ library OpenFOAM\textsuperscript{\textregistered} \cite{Weller1998}. For more details on GEA, see \cite{GQR_GEA,GQR_OF_clima}. Although we demonstrate numerically the accuracy and efficiency of EFR algorithm using a finite volume method for space discretization, the algorithm itself can be used with any space discretization method.

The outline of the paper is as follows. Sec.~\ref{sec:pd} describes the compressible Euler equations for low Mach stratified flows and introduces the filter stabilization for this model. In Sec.~\ref{sec:sd}, we discuss space discretization and the perturbation terms introduced by the filter stabilization to the Euler equations. Numerical results are presented in Sec.~\ref{sec:num_res} and conclusions are drawn in Sec.~\ref{sec:concl}.

\section{Problem definition}\label{sec:pd}
\subsection{The compressible Euler equations}

We consider mildly compressible Euler equations to describe the motion of the dry atmosphere, i.e., a compressible inviscid fluid, assumed to behave like an ideal gas. Let {$\Omega$} be a spatial domain of interest and $(0,t_f]$ a time interval of interest. Let $\rho$, $\bu = (u, v, w)$, and $p$ be the fluid density, velocity, and pressure. 
Moreover, let $e = c_v T + |\bu|^2/2 + g z$ be the total energy density, where $c_{v}$ is the specific heat capacity at constant volume, $T$ is the absolute temperature, $g$ is the gravitational constant, and $z$ is the vertical coordinate.
The conservation of mass, momentum, and total energy can be written as: 
\begin{align}
&\frac{\partial \rho}{\partial t} + \nabla \cdot (\rho \bu) = 0 &&\text{in}\,\,\Omega\times (0,t_f],   \label{eq:mass}  \\
&\frac{\partial (\rho \bu)}{\partial t} +  \nabla \cdot (\rho \bu \otimes \bu) + \nabla p   + \rho g \widehat{\bk} =\boldsymbol{0} &&\text{in}\,\, \Omega\times (0,t_f],   \label{eq:mom} \\
&\frac{\partial (\rho e)}{\partial t} +  \nabla \cdot (\rho e \bu) + \nabla \cdot (p \bu) = 0 &&\text{in}\,\, \Omega \times (0,t_f], 
\label{eq:ent}
\end{align}
where $\bk$ is the unit vector aligned with the vertical axis $z$.
We close system \eqref{eq:mass}-\eqref{eq:ent} using the following thermodynamics equation of state for $p$:
\begin{equation}
    p = \rho R T, \label{eq:GasLaw}
\end{equation}
where $R$ is the specific gas constant of dry air.

Let us to write the pressure as the sum of a fluctuation $p'$ with respect to a hydrostatic term:
\begin{equation}
    p=p'+\rho g z. \label{eq:pH}
\end{equation}
By plugging \eqref{eq:pH} into \eqref{eq:mom}, we obtain:
\begin{align}
    \frac{\partial (\rho \bu)}{\partial t} +  \nabla \cdot (\rho \bu \otimes \bu) + \nabla p' + gz \nabla \rho =0\quad \text{in}\,\, \Omega\times (0,t_f].  \label{eq:mom_split}
\end{align}

Let $c_{p}$ be the specific heat capacity at constant pressure for dry air and let
\begin{equation}\label{eq:K_h}
K = |\bu|^2/2, \quad h = c_v T + p/\rho = c_p T,    
\end{equation}
be the kinetic energy density and the specific enthalpy, respectively. The total energy density can be written as $e = h - p/\rho + K + gz$. Then, eq.~\eqref{eq:ent} can be rewritten as:
\begin{align}
\frac{\partial (\rho h)}{\partial t} +  \nabla \cdot (\rho \bu h) + 
\frac{\partial (\rho K)}{\partial t} +  \nabla \cdot (\rho \bu K) - \dfrac{\partial p}{\partial t}  +  
\rho g \bu \cdot \widehat{\bk} = 0,
\label{eq:over_ent}
\end{align}
where we have used eq.~\eqref{eq:mass} for further simplification. 

This paper focuses on formulation  \eqref{eq:mass},\eqref{eq:GasLaw}-\eqref{eq:over_ent} of the Euler equations.

A quantity of interest for atmospheric problems is the potential temperature
\begin{align}
\theta = \frac{T}{\pi}, \quad \pi = \left( \frac{p}{p_0} \right)^{\frac{R}{c_{p}}}, \label{eq:theta}
\end{align}
where $p_0 = 10^5$ Pa, which is the atmospheric pressure at the ground. Additionally, 
we define the potential temperature fluctuation $\theta'$ as the difference between $\theta$ and its mean hydrostatic value $\theta_0$:
\begin{align}
\theta'(x,y,z,t) = 
\theta(x,y,z,t) - \theta_0(z) . \label{eq:theta_split} 
\end{align}
See, e.g., \cite{kellyGiraldo2012} for more details. 

\subsection{Filter stabilization as an eddy viscosity model}\label{LeRayModel}

A numerical solution of system \eqref{eq:mass},\eqref{eq:GasLaw}-\eqref{eq:over_ent} computed with a mesh coarser than necessary for a DNS will be affected by non-physical oscillations that will eventually lead to a simulation breakdown. In order to avoid incurring into non-physical solutions, we adapt to the Euler equations an algorithm that has been shown to be accurate, efficient, and robust for the incompressible Navier-Stokes equations \cite{BQV,Bowers2012,abigail_CMAME,Ervin2012,layton_JMFM,GIRFOGLIO201927,LAYTON20113183,Olshanskii2013}. This algorithm consists of three steps: in the first step (called Evolve) one approximates the solution to the Euler equations with a coarse mesh, in the second step (called Filter) the numerical oscillations are smoothed out with a differential filter to obtain a filtered solution, and in the third step (called Relax) one combines the filtered and non-filtered solutions. This Evolve-Filter-Relax (EFR) algorithm, which is a computationally efficient realization of filter stabilization, is described next. 

Let $\Delta t \in \mathbb{R}$, $t^n = n \Delta t$, with $n = 0, ..., N_f$ and $t_f = N_f \Delta t$. Moreover, we denote by $y^n$ the approximation of a generic quantity $y$ at the time $t^n$. We adopt a Backward Differentiation Formula of order 1 (BDF1) for the discretization of the Eulerian time derivatives in \eqref{eq:mass},\eqref{eq:mom_split},\eqref{eq:over_ent}. 
Other time discretization schemes are possible (see, e.g., \cite{BQV,GIRFOGLIO201927,layton_JMFM,LAYTON20113183}).
The EFR algorithm reads as follows: given $\rho^0$, $\bu^0$, $h^0$,  $p^0$, and $T^0$, set $K^0 = |\bu^0|^2/2$ and for $n \geq 0$ perform the following steps:


\begin{itemize}
  \item[-] \emph{Step 1 - Evolve}: find density $\rho^{n+1}$ and intermediate variables  $\bv^{n+1},l^{n+1},K^{n+1}_{\bv},q^{n+1}, q'^{,n+1}$, $T^{n+1}_l$ such that: 
\begin{align}
& \frac{\rho^{n+1} - \rho^n}{\Delta t} + \nabla \cdot (\rho^{n+1} \bv^{n+1}) = 0, \label{eq:mass_td}  \\
&\frac{\rho^{n+1} \bv^{n+1} - \rho^{n}\bu^n}{\Delta t} +  \nabla \cdot (\rho^{n+1} \bv^{n+1} \otimes \bv^{n+1}) + \nabla q'^{,n+1} + gz \nabla \rho^{n+1} = \boldsymbol{0},  \label{eq:mom_td} \\
& \frac{\rho^{n+1} l^{n+1} - \rho^nh^n}{\Delta t} +  \nabla \cdot (\rho^{n+1} \bv^{n+1} l^{n+1}) + \frac{\rho^{n+1} K^{n+1}_{\bv} - \rho^n K^n}{\Delta t} +  \nabla \cdot (\rho^{n+1} \bv^{n+1} K^{n+1}_{\bv})  \cl
&\quad- \frac{q^{n+1} - p^n}{\Delta t} + \rho^{n+1} g \bv^{n+1}  \cdot \widehat{\bk} = 0, 
\label{eq:ent_td} \\
&q^{n+1} = q'^{,n+1} + \rho^{n+1} g z, \label{eq:p_td} \\
& q^{n+1} = \rho^{n+1} R T^{n+1}_l,\label{eq:p_td2}  \\
&l^{n+1} - l^{n} = c_p (T^{n+1}_l - T^n_l), \label{eq:T_td} \\
& K^{n+1}_{\bv} = \frac{|\bv^{n+1}|^2}{2}. \label{eq:K_td}
\end{align}
Notice that in \eqref{eq:T_td} we have chosen to update the value of the intermediate specific enthalpy in an incremental fashion.
  

  \item[-] \emph{Step 2 - Filter}: find filtered variables $\overline{\bv}^{n+1}, \overline{l}^{n+1}$  such that 
     \begin{align}
     &\overline{\bv}^{n+1}=F{\bv}^{n+1},    \label{eq:Fv}\\ 
     &\overline{l}^{n+1}=F l^{n+1}, \label{eq:Fh}
     \end{align}
     where F is a generic filter that could be linear or nonlinear. We will present possible choices for F in Sec.~\ref{sec:filter}.
 \item[-] \emph{Step 3 - Relax}: find end of step ${\bu}^{n+1}, h^{n+1}, K^{n+1}, p^{n+1}, p'^{,n+1}, T^{n+1}$ such that
\begin{align}
 &\bu^{n+1} = (1-\chi){\bv}^{n+1} + \chi \overline{\bv}^{n+1}, \label{eq:Rv} \\
&h^{n+1} = (1 -\xi){l}^{n+1} + \xi \overline{l}^{n+1}, \label{eq:Rh}\\ 
     &p^{n+1} = p'^{,n+1} + \rho^{n+1} g z, \label{eq:Rp} \\
& p^{n+1} = \rho^{n+1} R T^{n+1},\label{eq:Rp_2}  \\
&h^{n+1} - h^{n} = c_p (T^{n+1} - T^n), \label{eq:RT} \\
& K^{n+1} = \frac{|\bu^{n+1}|^2}{2}, \label{eq:RK} 
\end{align}
 where  $ \chi,\xi \in [0,1]$ are relaxation parameters.
\end{itemize}

The connection between the EFR algorithm and LES modeling is easily seen by shifting the index $n+1$ to $n$ in \eqref{eq:Fv}-\eqref{eq:Rh} and plugging them into \eqref{eq:mom_td}-\eqref{eq:ent_td} to obtain:
\begin{align}
&\frac{\rho^{n+1} \bv^{n+1} - \rho^{n}\bv^n}{\Delta t} +  \nabla \cdot (\rho^{n+1} \bv^{n+1} \otimes \bv^{n+1}) + \nabla q'^{,n+1} + gz \nabla \rho^{n+1} + \frac{\chi}{\Delta t} G \bv^n = \boldsymbol{0},  \label{eq:efr2} \\
& \frac{\rho^{n+1} l^{n+1} - \rho^n l^n}{\Delta t} +  \nabla \cdot (\rho^{n+1} \bv^{n+1} l^{n+1}) + \frac{\rho^{n+1} K^{n+1}_{\bv} - \rho^n K^n}{\Delta t} +  \nabla \cdot (\rho^{n+1} \bv^{n+1} K^{n+1}_{\bv})  \cl
&\quad- \frac{q^{n+1} - p^n}{\Delta t} + \rho^{n+1} g \bv^{n+1}  \cdot \widehat{\bk} + \frac{\xi}{\Delta t} G l^n= 0, 
\label{eq:efr3}  
\end{align}
with $G = I - F$, $I$ being the identity operator. System 
\eqref{eq:mass_td}, \eqref{eq:efr2}, \eqref{eq:efr3}, \eqref{eq:p_td}-\eqref{eq:K_td} gives us an implicit discretization of problem \eqref{eq:mass},\eqref{eq:GasLaw}-\eqref{eq:over_ent} with BDF1 and an additional
explicitly treated (linear or nonlinear) dissipation term. 


Let us assume that $\chi = \chi_0 \Delta t$ and $\xi= \xi_0 \Delta t$,
where $\chi_0$ and $\xi_0$ are time-independent constants.
Then, system 
\eqref{eq:mass_td}, \eqref{eq:efr2}, \eqref{eq:efr3}, \eqref{eq:p_td}-\eqref{eq:K_td}
can be seen as a time-stepping scheme for problem:
\begin{align}
&\frac{\partial \rho}{\partial t} + \nabla \cdot (\rho \bv) = 0,   \label{eq:evm1}  \\
&\frac{\partial (\rho \bv)}{\partial t} +  \nabla \cdot (\rho \bv \otimes \bv) + \nabla p' + gz \nabla \rho + \chi_0 G \bv= \boldsymbol{0},  \label{eq:evm2} \\
& \frac{\partial (\rho l)}{\partial t} +  \nabla \cdot (\rho \bu l) + 
\frac{\partial (\rho K)}{\partial t} +  \nabla \cdot (\rho \bu K) - \dfrac{\partial p}{\partial t}  +  
\rho g \bu \cdot \widehat{\bk} + \xi_0 G l= 0, \label{eq:evm3} \\
& p=p'+\rho g z, \label{eq:evm4} \\
& p = \rho R T, \label{eq:evm5} \\
& h = c_p T, \label{eq:evm6} \\
& K = |\bv|^2/2. \label{eq:evm7}
\end{align}
Thus, filter stabilization algorithm \eqref{eq:mass_td}-\eqref{eq:RK} can be interpreted as a splitting scheme for problem (\ref{eq:evm1})-(\ref{eq:evm7}).

 Notice that model (\ref{eq:evm1})-(\ref{eq:evm7}) can be considered as a LES model of the eddy-viscosity type with closure:
\begin{align}
    \nabla \cdot (\rho\overline{\bv  \otimes \bv } -  \rho\overline{\bv} \otimes \overline{\bv}) \approx \chi_0 G \bv, \\
    \nabla \cdot(\rho \overline{\bv l } - \rho \overline{ \bv} \overline{l} )\approx \xi_0 G l. 
\end{align}
This shows the connection between algorithm \eqref{eq:mass_td}-\eqref{eq:RK} and LES
modeling.

\subsection{A possible choice for the filter}\label{sec:filter}

We will consider the following filter for step 2
\eqref{eq:Fv}-\eqref{eq:Fh}:
\begin{equation}\label{eq:FH}
    F = (I + L)^{-1},\quad L=-\nabla \cdot (\delta \nabla) 
\end{equation}
where $\delta>0$ is a linear or nonlinear artificial ``viscosity''. Such a filter applied to $\bv^{n+1}$ as in (\ref{eq:Fv}) amounts to solving the following problem: find $\overline{\bv}^{n+1}$ such that
\begin{equation}
    -\nabla \cdot (\delta \nabla(\overline{\bv}^{n+1})) + \overline{\bv}^{n+1} = \bv^{n+1} ,\quad \delta=\alpha^2 a(\bv^{n+1}), \label{eq:filter}
\end{equation}
where $\alpha$ can be interpreted as the filtering radius
and $a(\cdot) \in (0,1]$ is the so-called indicator function. Note that $\delta$ is not properly a viscosity since it has the dimension of a length a square. However, if we multiply \eqref{eq:filter} by $\rho^{n+1}/\Delta t$, we obtain Stokes problem:
\begin{equation}
    \frac{\rho^{n+1}}{\Delta t} (\overline{\bv}^{n+1}-\bv^{n+1}) -\nabla \cdot (\overline{\mu}\nabla \overline{\bv}^{n+1}) = \boldsymbol{0},\quad \overline{\mu} = \rho^{n+1}\frac{\alpha^2}{\Delta t}a(\bv^{n+1}), \label{eq:ufil}
\end{equation}
where $\overline{\mu}$ is dimensionally a dynamic viscosity. 

The same filter applied to $l^{n+1}$ as in \eqref{eq:Fh} yields:
\begin{equation}
\frac{\rho^{n+1}}{\Delta t}
    (\overline{l}^{n+1}-l^{n+1})  -\nabla \cdot (\overline{\mu}\nabla \overline{l}^{n+1}) = 0. \label{eq:hfil}
\end{equation}

In summary, the EFR algorithm we will use in this paper entails performing the following steps:
\begin{itemize}
\item[-] \emph{Step 1 - Evolve}: find density $\rho^{n+1},\bv^{n+1},l^{n+1},K^{n+1}_{\bv},q^{n+1}, q'^{,n+1}$, $T^{n+1}_l$ such that \eqref{eq:mass_td}-\eqref{eq:K_td} hold. 
\item[-] \emph{Step 2 - Filter}: find filtered variables $\overline{\bv}^{n+1}, \overline{l}^{n+1}$  such that  \eqref{eq:ufil}-\eqref{eq:hfil} hold.
\item[-] \emph{Step 3 - Relax}: set \eqref{eq:Rv}-\eqref{eq:RK}.
\end{itemize}

\subsection{Possible choices for the indicator function}

The success of the EFR algorithm in the simulation of atmospheric flows ultimately depends on the reliability of the indicator function. 
The indicator function has to be such that it takes values close to zero where its argument (i.e., the Euler velocity or specific enthalpy) does not need regularization, while it takes values close to 1 where the argument does need to be regularized.
Different choices for indicator function $a(\cdot)$ have been proposed in the literature for the incompressible Navier-Stokes equations \cite{Borggaard2009,Bowers2012,O-hunt1988,layton_JMFM}. Some indicator functions \cite{Bowers2012,layton_JMFM} are based on physical quantities that
are known to vanish for coherent flow structures. The drawback for these indicator functions is that they do not allow for a rigorous 
convergence theory to verify the robustness of the associated filtering method. Hence,  mathematics-based (instead of physics-based) indicator functions were proposed \cite{abigail_CMAME,layton_JMFM}. 
In this paper, we
will consider and compare three mathematics-based choices. 

The first and easiest choice corresponds to a linear filter, i.e., we take
\begin{equation}
a(\bv) = a_L(\bv) = 1, \label{eq:a_lin}
\end{equation}
in \eqref{eq:ufil}. Besides linearity, another advantage of this choice is that it makes the operator in the filter equations constant in time. However, its efficacy is rather limited, since it introduces the same amount of regularization everywhere in the domain. This is likely to introduce overdiffusion as we will show in Sec.~\ref{sec:num_res}

A second mathematically convenient indicator function is 
\begin{equation}
a(\bv) = a_S(\bv) = \frac{| \nabla \bv |}{\|\nabla \bv\|_{\infty}}, \label{eq:a_smago}
\end{equation}
which has strong monotonicity properties. With $a_S(\cdot)$ as indicator function for the EFR algorithm, we recover a Smagorinsky-like model, which is an improvement over the linear filter obtained with $a_L(\cdot)$. 

Finally, we consider a class of deconvolution-based indicator functions, which were shown to be particularly accurate for realistic incompressible flow problems
\cite{BQV,GIRFOGLIO201927}.
Such functions are defined as:
\begin{equation}
a(\bv) = a_{D}(\bv) = \left|  \bv - D (F(\bv)) \right|, \label{eq:a_deconv}
\end{equation}
where $F$ is the linear Helmholtz filter (i.e., \eqref{eq:FH} with $\delta$ constant in space and time)
and $D$ is the Van Cittert deconvolution:
\begin{equation}\label{eq:D}
D = \sum_{n = 0}^N (I - F)^n. 
\end{equation}
We remark that $D$ is a bounded regularized approximation of $F^{-1}$. Typically, $N$ in \eqref{eq:D} is set to 0, 1 \cite{BQV,GIRFOGLIO201927}. In this paper,  we consider $N = 0$, which means $D=I$. For this choice of $N$, indicator function (\ref{eq:a_deconv}) becomes
\begin{align}
a_D(\bv) = \left|  \bv - F(\bv) \right|. \label{eq:a_D0_a_D1}
\end{align}


\section{Space discretization of the steps in the EFR algorithm}\label{sec:sd}

For space discretization,
we adopt a finite volume method. The Evolve step is the most computationally intensive step in the EFR algorithm and to contain its computational cost we use a splitting scheme thoroughly described in \cite{GQR_OF_clima}. This section focuses on the space discretization of the Filter and Relax steps.

Let us consider a partition of the computational domain $\Omega$ into cells or control volumes $\Omega_i$, with $i = 1, \dots, N_{c}$, where $N_{c}$ is the total number of cells in the mesh. 
Let  \textbf{A}$_j$ be the surface vector of each face of the control volume, 
with $j = 1, \dots, M$. We will start with the space discretization of the Filter problem (\ref{eq:ufil})-(\ref{eq:hfil}) and then write the space discrete version of the Relax Step \eqref{eq:Rv}-\eqref{eq:RK}

The integral form of the eq.~(\ref{eq:ufil}) for each volume $\Omega_i$ is given by: 
\begin{equation*}
\frac{1}{\Delta t}    \int_{\Omega_i} \rho^{n+1} \overline{\bv}^{n+1}\,d\Omega - \int_{\Omega_i} \nabla \cdot (\overline{\mu} \nabla \overline{\bv}^{n+1} )\,d\Omega = \frac{1}{\Delta t} \int_{\Omega_i} \rho^{n+1} \bv^{n+1} \,d\Omega.
\end{equation*}
By using the Gauss-divergence theorem, the above equation becomes:
\begin{equation}\label{eq:fu_gauss}
  \frac{1}{\Delta t}  \int_{\Omega_i} \rho^{n+1} \overline{\bv}^{n+1}\,d\Omega - \int_{\Omega_i} ( \overline{\mu} \nabla \overline{\bv}^{n+1} )\cdot d\bA = \frac{1}{\Delta t}  \int_{\Omega_i} \rho^{n+1} \bv^{n+1} d\Omega.
\end{equation}
Let us denote with $(\overline{\mu} \nabla \overline{\bv}^{n+1})_i$ and $\overline{\bv}^{n+1}_i$ the average stress tensor and filtered density in control volume $\Omega_i$, respectively. 
Similarly, we denote with $\rho^{n+1}_i$ and ${\bv}^{n+1}_i$ the average density and intermediate velocity in $\Omega_i$.
Then, eq.~\eqref{eq:fu_gauss} is approximated as follows:
\begin{equation}
   \frac{1}{\Delta t} \rho^{n+1}_i \overline{\bv}_i^{n+1} - \sum_j (\overline{\mu} \nabla \overline{\bv}^{n+1})_{i,j} \cdot \textbf{A}_j  = \frac{1}{\Delta t} \rho^{n+1}_i {\bv}_i^{n+1}. \label{eq:discreteV}
\end{equation}
We choose to approximate 
the gradient of $\overline{\bv}^{n+1}_i$ at face $j$ with second order accuracy. 
See \cite{jasakphd} for more details. 

Following a similar procedure for (\ref{eq:hfil}), we obtain:
\begin{align}
   \frac{1}{\Delta t} \rho^{n+1}_i \overline{l}_i^{n+1} - \sum_j (\overline{\mu} \nabla \overline{l}^{n+1})_{i,j} \cdot \textbf{A}_j  = \frac{1}{\Delta t} \rho^{n+1}_i {l}_i^{n+1}, \label{eq:discreteh}
\end{align}
where $\overline{l}^{n+1}_i$ and ${l}^{n+1}_i$ are the average filtered and intermediate specific enthalpy in 
$\Omega_i$. For the approximation of the gradient of $\overline{l}^{n+1}_i$ at face $j$, we use the same formula used for the components of  $\overline{\bv}^{n+1}_i$. 

Now, let us turn to the Relax step. The discretized form of each equation in the Relax step is simply given by taking the average of each variable in $\Omega_i$: 
\begin{align}
 &\bu^{n+1}_i = (1-\chi){\bv}^{n+1}_i + \chi \overline{\bv}^{n+1}_i, \label{eq:Rv_sd} \\
&h^{n+1}_i = (1 -\xi){l}^{n+1}_i + \xi \overline{l}^{n+1}_i, \label{eq:Rh_sd}\\ 
     &p^{n+1}_i = p'^{,n+1}_i + \rho^{n+1}_i g z_i, \label{eq:Rp_sd} \\
& p^{n+1}_i = \rho^{n+1} R T^{n+1}_i,\label{eq:Rp_2_sd}  \\
&h^{n+1}_i - h^{n}_i = c_p (T^{n+1}_i - T^n_i), \label{eq:RT_sd} \\
& K^{n+1}_i = \frac{|\bu^{n+1}_i|^2}{2}, \label{eq:RK_sd} 
\end{align}
where $z_i$ is the vertical coordinate of the centriod of cell $\Omega_i$.

\subsection{The EFR algorithm as a solver for the Euler equations with perturbations}

In this section, we will show that by combining the equations at the Evolve, Filter, and Relax steps we obtain the Euler equations perturbed by some extra terms and we discuss such terms. 

Let us use a subindex $h$ to denote the space-discrete solution, where $h$ refers to the mesh size. We rewrite \eqref{eq:mom_td}
\begin{align}
    &\frac{\rho^{n+1}_h \bv^{n+1}_h - \rho^{n}_h\bu^n_h}{\Delta t} +  \nabla \cdot (\rho^{n+1}_h \bv^{*}_h \otimes \bv^{n+1}_h) + H(q'^{,n+1}_h,\rho^{n+1}_h) = \boldsymbol{0},  \label{eq:mom_td_sd}
\end{align}
where $\bv^{*}_h$ is a suitable extrapolation of $\bv^{n+1}_h$ to linearize the convective term and $H(q'^{,n+1}_h,\rho^{n+1}_h) = \nabla q'^{,n+1}_h + gz \nabla \rho^{n+1}_h$. Let us also write the space-discrete version of eq.~\eqref{eq:Rv}
\begin{align}
\bu^{n+1}_h = (1-\chi){\bv}^{n+1}_h + \chi \overline{\bv}^{n+1}_h, \label{eq:Rv_sd}    
\end{align}
and eq.~\eqref{eq:ufil}
\begin{equation}
    \frac{\rho^{n+1}_h}{\Delta t} (\overline{\bv}^{n+1}_h-\bv^{n+1}_h) -\nabla \cdot (\overline{\mu}_h \nabla \overline{\bv}^{n+1}_h) = \boldsymbol{0}, \label{eq:ufil_sd}
\end{equation}
where 
\begin{align}
    \overline{\mu}_h = \rho^{n+1}_h\frac{\alpha^2}{\Delta t}a(\bv^{n+1}_h) \label{eq:mubar_h}.
\end{align}

We multiply \eqref{eq:ufil_sd} by $\chi$, add it to \eqref{eq:mom_td_sd}, and make use of \eqref{eq:Rv_sd} to obtain:
\begin{align}
    &\frac{\rho^{n+1}_h \bu^{n+1}_h - \rho^{n}_h\bu^n_h}{\Delta t} +  \nabla \cdot (\rho^{n+1}_h \bv^{*}_h \otimes \bv^{n+1}_h) + H(q'^{,n+1}_h,\rho^{n+1}_h) - \chi \nabla \cdot (\overline{\mu}_h \nabla \overline{\bv}^{n+1}_h) = \boldsymbol{0}. \el
\end{align}
Using \eqref{eq:Rv_sd} once more, we get:
\begin{align}
    &\frac{\rho^{n+1}_h \bu^{n+1}_h - \rho^{n}_h\bu^n_h}{\Delta t} +  \nabla \cdot (\rho^{n+1}_h \bv^{*}_h \otimes \bu^{n+1}_h) + H(q'^{,n+1}_h,\rho^{n+1}_h) \cl
    & \quad + \chi \nabla \cdot (\rho^{n+1}_h \bv^{*}_h \otimes (\bv^{n+1}_h - \overline{\bv}^{n+1}_h))  - \chi \nabla \cdot (\overline{\mu}_h \nabla \overline{\bv}^{n+1}_h) = \boldsymbol{0}. \el
\end{align}
which can be rewritten as 
\begin{align}
    &\frac{\rho^{n+1}_h \bu^{n+1}_h - \rho^{n}_h\bu^n_h}{\Delta t} +  \nabla \cdot (\rho^{n+1}_h \bv^{*}_h \otimes \bu^{n+1}_h)  + H(q'^{,n+1}_h,\rho^{n+1}_h) - \chi \nabla \cdot (\overline{\mu}_h \nabla {\bu}^{n+1}_h)\cl
    & \quad + \chi \nabla \cdot (\rho^{n+1}_h \bv^{*}_h \otimes (\bv^{n+1}_h - \overline{\bv}^{n+1}_h))  - \chi \nabla \cdot (\overline{\mu}_h \nabla ( \overline{\bv}^{n+1}_h - {\bu}^{n+1}_h)) = \boldsymbol{0}. \label{eq:pert_Eu1}
\end{align}
In \eqref{eq:pert_Eu1}, we have explicitly written a diffusive term involving only the end-of-step velocity $\bu^{n+1}_h$, i.e., the last term in the first line. 
The last term at the left-hand side in \eqref{eq:pert_Eu1} can be rewritten using \eqref{eq:Rv_sd} to get:
\begin{align}
    &\frac{\rho^{n+1}_h \bu^{n+1}_h - \rho^{n}_h\bu^n_h}{\Delta t} +  \nabla \cdot (\rho^{n+1}_h \bv^{*}_h \otimes \bu^{n+1}_h)  + H(q'^{,n+1}_h,\rho^{n+1}_h) - \chi \nabla \cdot (\overline{\mu}_h \nabla {\bu}^{n+1}_h)\cl
    & \quad + \chi \nabla \cdot (\rho^{n+1}_h \bv^{*}_h \otimes (\bv^{n+1}_h - \overline{\bv}^{n+1}_h))  - \chi (1-\chi) \nabla \cdot (\overline{\mu}_h \nabla ( \overline{\bv}^{n+1}_h - {\bv}^{n+1}_h)) = \boldsymbol{0}. \label{eq:pert_Eu2}
\end{align}
Eq.~\eqref{eq:pert_Eu2} shows that the end-of-step velocity $\bu^{n+1}_h$ provided by the EFR algorithm solves a perturbed discrete momentum balance equation. The perturbation consists of three terms, all multiplied by $\chi$: an extra convection term and two extra diffusion terms. 
As mentioned in Sec.~\eqref{LeRayModel}, $\chi$ should be a multiple of $\Delta t$. Thus, as $\Delta t$ tends to zero, the perturbation terms vanish and we recover the discrete momentum balance equation. In addition, we notice that when $\alpha^2/\Delta t$ tends to zero, the artificial viscosity $\overline{\mu}_h$ \eqref{eq:mubar_h} vanishes and $\overline{\bv}^{n+1}_h$ tends to ${\bv}^{n+1}_h$. 
If $\alpha$ is a multiple of $h$ (as it is typically the case), then eq.~\eqref{eq:pert_Eu2} is consistent with eq.~\eqref{eq:mom} so long as the mesh size and time step go to zero at the same rate.

Combining in a similar fashion the space discrete version of eq.~\eqref{eq:ent_td}: 
\begin{align}
    & \frac{\rho^{n+1}_h l^{n+1}_h - \rho^n_h h^n_h}{\Delta t} +  \nabla \cdot (\rho^{n+1}_h \bv^{*}_h l^{n+1}_h) + G\left(\rho^{n+1}_h, K^{n+1}_{\bv,h},q^{n+1}_h, \bv^{n+1}_h\right)  = 0, \cl
& G\left(\rho^{n+1}_h, K^{n+1}_{\bv,h},q^{n+1}_h, \bv^{n+1}_h\right) = \frac{\rho^{n+1}_h K^{n+1}_{\bv,h} - \rho^n_h K^n_h}{\Delta t} +  \nabla \cdot (\rho^{n+1}_h \bv^{*}_h K^{n+1}_{\bv,h}) \cl
&\quad - \frac{q^{n+1}_h - p^n_h}{\Delta t} + \rho^{n+1}_h g \bv^{n+1}_h  \cdot \widehat{\bk}, \el
\end{align}
with the discrete versions of eq.~\eqref{eq:hfil}:
\begin{equation*}
\frac{\rho^{n+1}_h}{\Delta t}
    (\overline{l}^{n+1}_h-l^{n+1}_h)  -\nabla \cdot (\overline{\mu}_h\nabla \overline{l}^{n+1}_h) = 0,
\end{equation*}
and eq.~\eqref{eq:Rh}:
\begin{equation*}
h^{n+1}_h = (1 -\xi){l}^{n+1}_h + \xi \overline{l}^{n+1}_h,
\end{equation*}
we obtain
\begin{align}
    & \frac{\rho^{n+1}_h h^{n+1}_h - \rho^n_h h^n_h}{\Delta t} +  \nabla \cdot (\rho^{n+1}_h \bv^{*}_h h^{n+1}_h) + G\left(\rho^{n+1}_h, K^{n+1}_{\bv,h},q^{n+1}_h, \bv^{n+1}_h\right) -  \xi 
    \nabla \cdot (\overline{\mu}_h \nabla {h}^{n+1}_h)\cl
    & \quad + \xi \nabla \cdot (\rho^{n+1}_h \bv^{*}_h  (l^{n+1}_h - \overline{l}^{n+1}_h))  - \xi (1-\xi) \nabla \cdot (\overline{\mu}_h \nabla ( \overline{l}^{n+1}_h - {l}^{n+1}_h)) = 0. \label{eq:pert_Eu3}
\end{align}
Eq.~\eqref{eq:pert_Eu3} is a perturbed discrete enthalpy balance equation, with the perturbation coming from an extra convection term and two extra diffusion terms. Like in the case of eq.~\eqref{eq:pert_Eu2}, the perturbation terms are multiplied by the relaxation parameter.

\begin{rem}\label{rem:1}
For $\chi = 1$, which corresponds to accepting the filtered velocity as the end-of-step velocity, eq.~\eqref{eq:pert_Eu2} becomes:
\begin{align}
    &\frac{\rho^{n+1}_h \bu^{n+1}_h - \rho^{n}_h\bu^n_h}{\Delta t} +  \nabla \cdot (\rho^{n+1}_h \bv^{*}_h \otimes \bu^{n+1}_h)  + H(q'^{,n+1}_h,\rho^{n+1}_h) - \nabla \cdot (\overline{\mu}_h \nabla {\bu}^{n+1}_h)\cl
    & \quad + \nabla \cdot (\rho^{n+1}_h \bv^{*}_h \otimes (\bv^{n+1}_h - \overline{\bv}^{n+1}_h))  = \boldsymbol{0}, \el
\end{align}
where we clearly see that the artificial diffusion introduced by the EFR algorithm is $\overline{\mu}_h$ \eqref{eq:mubar_h}. 
Similarly, by setting $\xi = 1$ (i.e, we take the filtered enthalpy as the end-of-step enthalpy) eq.~\eqref{eq:pert_Eu3} becomes:
\begin{align}
    & \frac{\rho^{n+1}_h h^{n+1}_h - \rho^n_h h^n_h}{\Delta t} +  \nabla \cdot (\rho^{n+1}_h \bv^{*}_h h^{n+1}_h) + G\left(\rho^{n+1}_h, K^{n+1}_{\bv,h},q^{n+1}_h, \bv^{n+1}_h\right) - 
    \nabla \cdot (\overline{\mu}_h \nabla {h}^{n+1}_h)\cl
    & \quad + \nabla \cdot (\rho^{n+1}_h \bv^{*}_h  (l^{n+1}_h - \overline{l}^{n+1}_h)) = 0. \el
\end{align}
\end{rem}

\begin{rem}
Eddy viscosity models are equivalent to introducing additional terms in eq.~\eqref{eq:mom_split} and \eqref{eq:over_ent} of the form
\begin{align}
&\frac{\partial (\rho \bu)}{\partial t} +  \nabla \cdot (\rho \bu \otimes \bu) + \nabla p' + gz \nabla \rho -  \nabla \cdot (2 \mu_a \boldsymbol{\epsilon}(\bu)) + \nabla \left(\frac{2}{3}\mu_a \nabla \cdot \bu \right)= 0,  \label{eq:mom_LES} \\
&\frac{\partial (\rho h)}{\partial t} +  \nabla \cdot (\rho \bu h) + 
\frac{\partial (\rho K)}{\partial t} +  \nabla \cdot (\rho \bu K) - \dfrac{\partial p}{\partial t}  +  
\rho g \bu \cdot \widehat{\bk}  - \nabla \cdot \left(\frac{\mu_a}{Pr} \nabla h \right) = 0. 
\label{eq:ent_LES}
\end{align}
where $\mu_a$ is an artificial viscosity (defined differently for the different LES models), $\boldsymbol{\epsilon}(\bu) = (\nabla \bu + (\nabla \bu)^T)/2$ is the strain-rate tensor, and $Pr$ is the Prandtl number, i.e., the dimensionless number defined as the ratio of momentum diffusivity to thermal diffusivity. Comparing \eqref{eq:mom_LES}-\eqref{eq:ent_LES} with \eqref{eq:pert_Eu2}-\eqref{eq:pert_Eu3} suggests choosing $\xi = \chi/Pr$.
\end{rem}

\begin{rem}\label{rem:3}
The Smagorinsky model sets $\mu_a$ in \eqref{eq:mom_LES}-\eqref{eq:ent_LES} as
\begin{align}
\mu_a = \rho (C_s \delta)^2 \sqrt{ 2 \boldsymbol{\epsilon} : \boldsymbol{\epsilon}}, \quad C_s^2 = C_k \sqrt{\dfrac{C_k}{C_{\epsilon}}} \label{eq:smago}
\end{align}
where $\delta$ is the filter width (typically comparable with the mesh size), and $C_k$ and $C_{\epsilon}$ are model parameters. In order to obtain the same amount of artificial viscosity with the EFR algorithm and $a_S$ \eqref{eq:a_smago}, one needs $\alpha \approx C_s \delta \sqrt{\Delta t} \|\nabla \bv_h^{n+1}\|_{\infty}$ at time $t^{n+1}$. In practice, one can easily calculate $C_s \delta \sqrt{\Delta t}$ while $\|\nabla \bv_h^{n+1}\|_{\infty}$ would have to be guessed to tune $\alpha$.
\end{rem}

\section{Numerical Results}\label{sec:num_res}

We consider two well-known benchmarks: 
the rising thermal bubble as presented in \cite{ahmadLindeman2007,Feng2021} and the density current \cite{carpenterDroegemeier1990,strakaWilhelmson1993}.
Both test cases involve a perturbation of 
a neutrally stratified atmosphere with uniform background potential temperature over a flat terrain
and the boundaries are treated as if the problem were inviscid (i.e., free-slip boundary conditions are imposed). So, these are not boundary layer flow problems. 
We presents our results for the rising rising thermal bubble and the density current in Sec.~\ref{sec:bubble} and \ref{sec:DC}, respectively, and compare them with other numerical data available in the literature since there is no exact solution for these benchmarks.

\subsection{Rising thermal bubble}\label{sec:bubble}

In computational domain $\Omega = [0,5000] \times [0,10000]$ m$^2$, a neutrally stratified atmosphere with uniform background potential temperature $\theta_0$=300 K is perturbed by a circular bubble of warmer air.
The initial temperature field is
\begin{equation}
\theta^0 = 300 + 2\left[ 1 - \frac{r}{r_0} \right] ~ \textrm{if $r\leq r_0=2000~\mathrm{m}$}, \quad\theta^0 = 300
~ \textrm{otherwise},
\label{warmEqn1}
\end{equation}
where $r = \sqrt[]{(x-x_{c})^{2} + (z-z_{c})^{2}}$, $(x_c,z_c) = (5000,2000)~\mathrm{m}$ is the radius of the circular perturbation \cite{ahmadLindeman2007,ahmad2018}.
The initial density is given by \begin{align}
\rho^0 = \frac{p_g}{R \theta_0} \left(\frac{p}{p_g}\right)^{c_{v}/c_p}, \quad p = p_g \left( 1 - \frac{g z}{c_p \theta^0} \right)^{c_p/R}, \label{eq:rho_wb}
\end{align}
with 
$c_p = R + c_v$, $c_v = 715.5$ J/(Kg K), $R = 287$ J/(Kg K). 
The initial velocity field is zero everywhere. 
Finally, the initial specific enthalpy is given by:
\begin{align}
h^{0} = c_p \theta^0 \left( \frac{p}{p_g} \right)^{\frac{R}{c_{p}}}.
\label{eq:e0}
\end{align}
We let the bubble evolve in the time interval of interest is $(0,1020]$ s. Impenetrable, free-slip boundary conditions are imposed on all walls.

We consider five different meshes with uniform resolution $h = \Delta x = \Delta z = 125, 62.5$, $31.25$, $15.625$~m. The time step is set to $\Delta t = 0.1$ s for all the simulations. In all the cases, we set $\chi = \xi = 1$ so that the artificial diffusion introduced by the EFR algorithm can easily be calculated (see Remark \ref{rem:1}).

\begin{figure}[htb!]
     \centering
     \begin{overpic}[percent,width=0.24\textwidth]{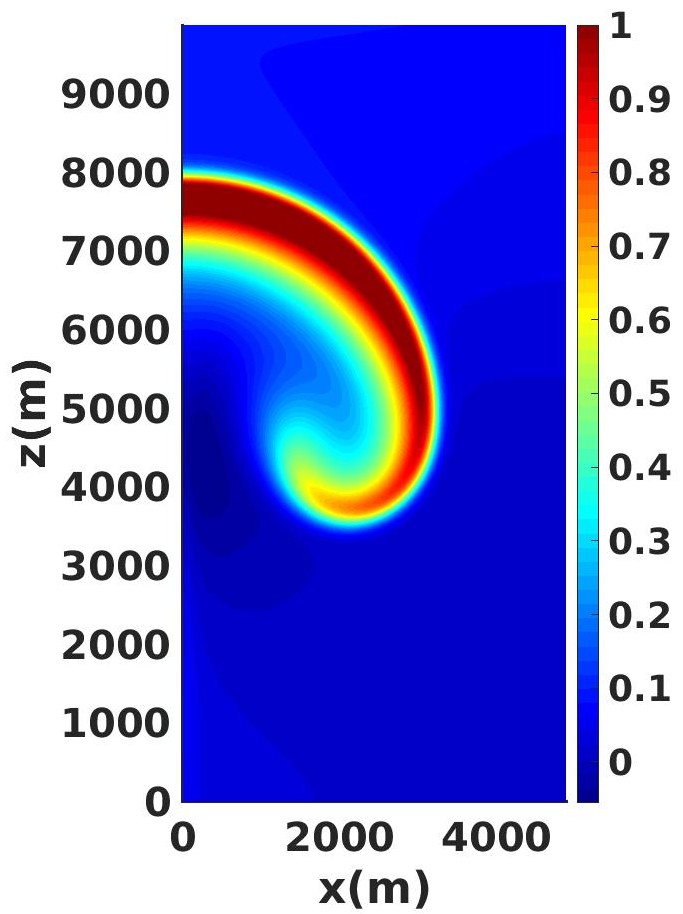}
    \put(22,90){\textcolor{white}{\footnotesize{$h =15.625$ m}}}
    \end{overpic}
    \begin{overpic}[percent,width=0.242\textwidth]{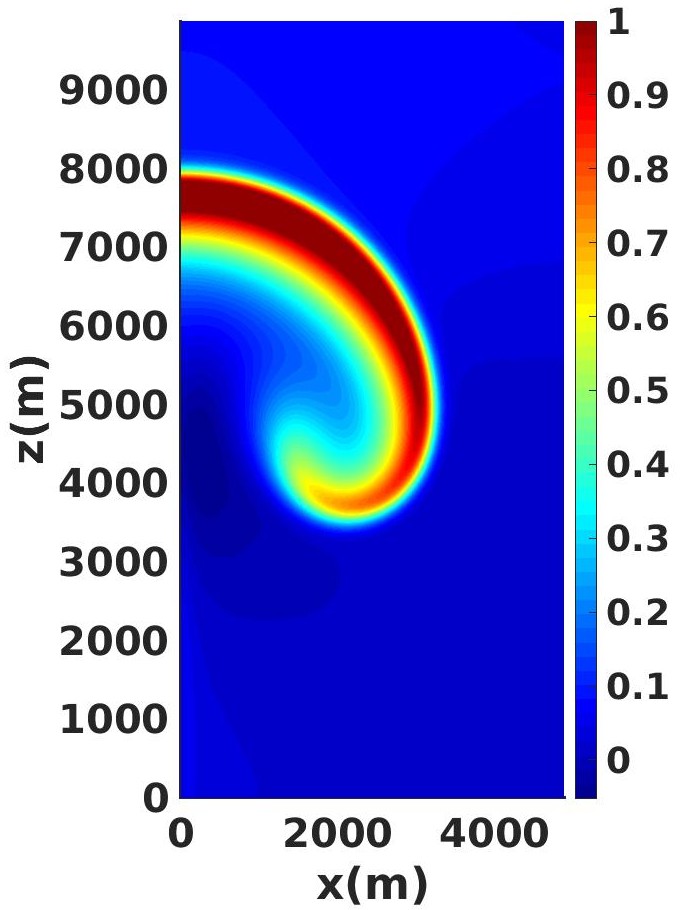}
    \put(22,90){ \textcolor{white}{\footnotesize{$h =31.25$ m}}}
    \end{overpic} 
    \begin{overpic}[percent,width=0.242\textwidth]{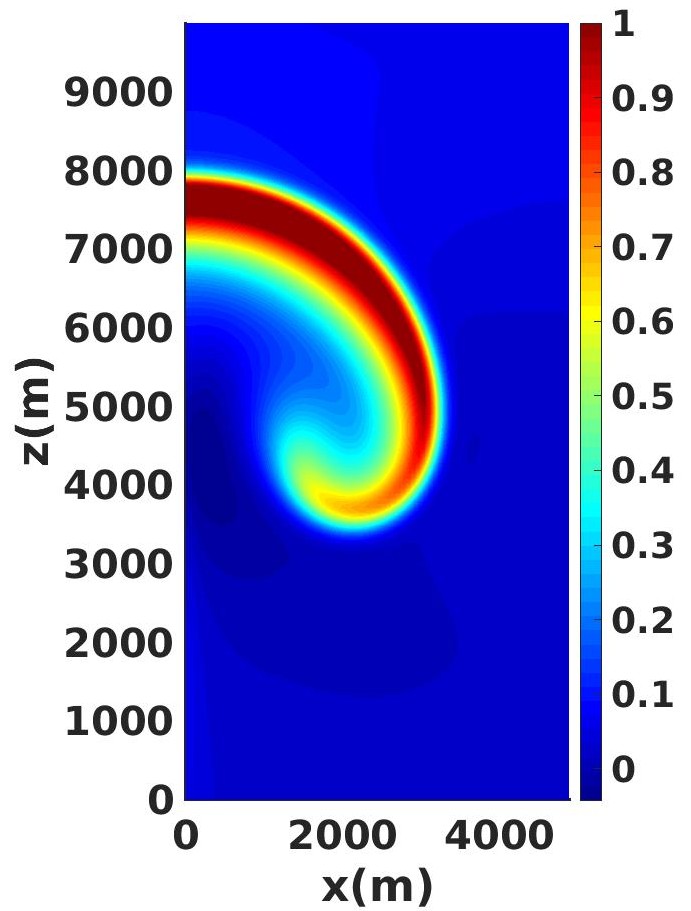}
    \put(23,90){ {\textcolor{white}{\footnotesize{ $h=62.5$ m}}}}
    \end{overpic}
    \begin{overpic}[percent,width=0.242\textwidth]{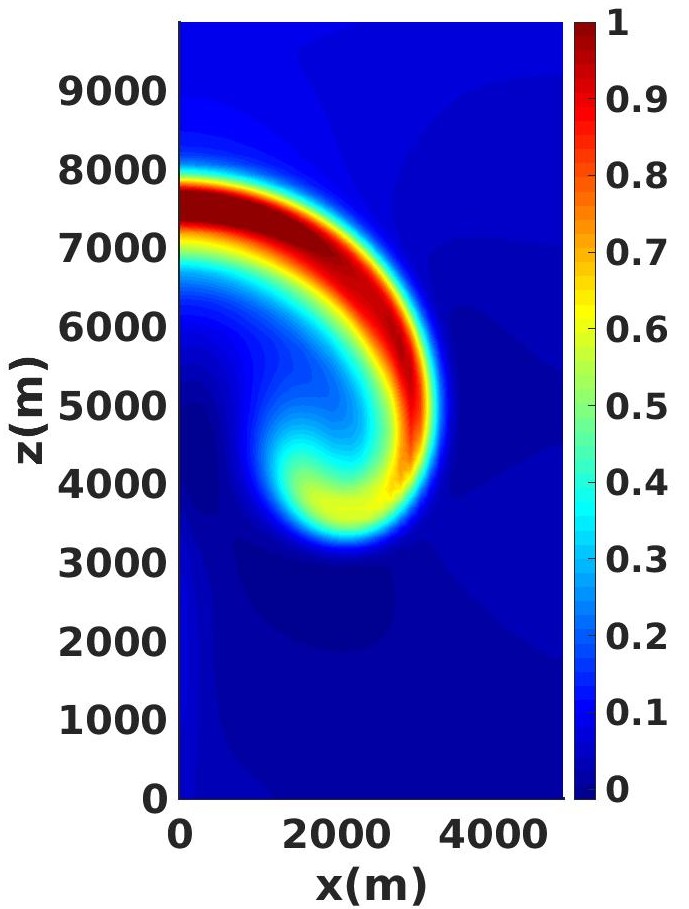}
    \put(23,90){ {\textcolor{white}{\footnotesize{$h=125$ m}}}}
    \end{overpic} 
    \caption{Rising thermal bubble, $a_L$, $\alpha = 1.9$ m: perturbation of potential temperature at $t = 1020$ s computed with four different meshes. The mesh size, specified in each panel, is increasing from left to right. }
    \label{fig:TB_Linear}
\end{figure}

We start with the linear filter, i.e., we take $a_L$ \eqref{eq:a_lin} as indicator function, because it allows us to make a direct comparison with the results obtained by setting $\mu_a = 15$ and $Pr = 1$ in \eqref{eq:mom_LES}-\eqref{eq:ent_LES}
\cite{ahmadLindeman2007,GQR_OF_clima}. We note that both of these are ah-hoc values chosen by the authors of \cite{ahmadLindeman2007} to stabilize the numerical simulations. It is not unusual in benchmarks to set $Pr = 1$ although the air Prandtl number is about 0.71 at 20$^\circ$C (see, e.g., \cite{giraldoRestelli2008b}). Other authors have chosen other arbitrary values, like $Pr = 0.1$ in \cite{nazarovHoffman2013}. For a qualitative analysis of the results for the rising thermal bubble as $Pr$ varies we refer the reader to \cite{marrasNazarovGiraldo2015}.
To introduce the same amount of artificial viscosity with the EFR algorithm and $a_L$, we use \eqref{eq:mubar_h} to get $\alpha \approx \sqrt{\mu_a \Delta t / \rho^{n+1}}$ at time $t^{n+1}$. For simplicity, we keep $\alpha$ constant in space and time and set it to $1.9$, which is obtained by taking the minimum value of density. 
Fig.~\ref{fig:TB_Linear} shows the perturbation of potential temperature $\theta'$ at $t = 1020$ s computed by the EFR algorithm (with $a_L$ and $\alpha = 1.9$) and all the meshes under consideration. 
From Fig.~\ref{fig:TB_Linear}, 
we observe no visible change
in the computed $\theta'$ when the mesh 
is refined past $h = 62.5$ m. In order to facilitate the comparison of the panels in Fig.~\ref{fig:TB_Linear} with data in the literature \cite{ahmadLindeman2007,ahmad2018,marrasNazarovGiraldo2015,GQR_OF_clima}, we have forced the colorbar to range from 0 to 1. Overall, these results are in very good qualitative agreement with the corresponding figures in \cite{ahmadLindeman2007,ahmad2018,GQR_OF_clima,marrasNazarovGiraldo2015}.

\begin{figure}[htb!]
     
     \centering
     
     \begin{overpic}[width=0.19\textwidth]{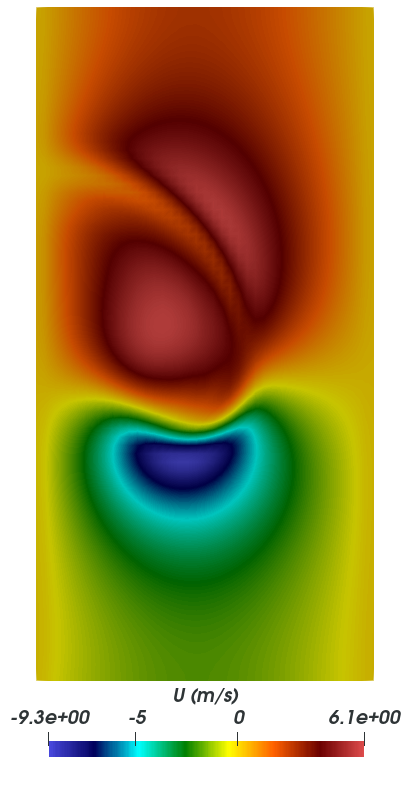}
    \end{overpic} \vspace{0.12cm}
    \begin{overpic}[width=0.195\textwidth]{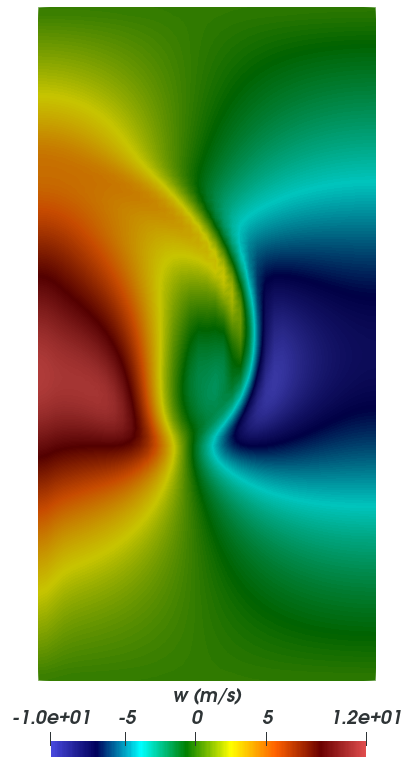}
    \end{overpic}  
    \caption{Rising thermal bubble, $a_L$, $\alpha = 1.9$ m: contour plots of the horizontal velocity component $u$ (left) and the vertical velocity component $w$ (right) at $t = 1020$ s computed with mesh $h = $ 125 m.}
    \label{fig:TB_velocities}
\end{figure}

We obtain good qualitative agreement with data in the literature (e.g., Fig.~7 in \cite{ahmadLindeman2007}) also for 
Fig.~\ref{fig:TB_velocities}, which displays velocity components $u$ and $w$ at $t = 1020$ s computed by the EFR algorithm (with $a_L$ and $\alpha = 1.9$) with mesh $h = $ 125 m. 

\begin{figure}[htb!]
      \centering
     \begin{overpic}[width=0.48\textwidth]{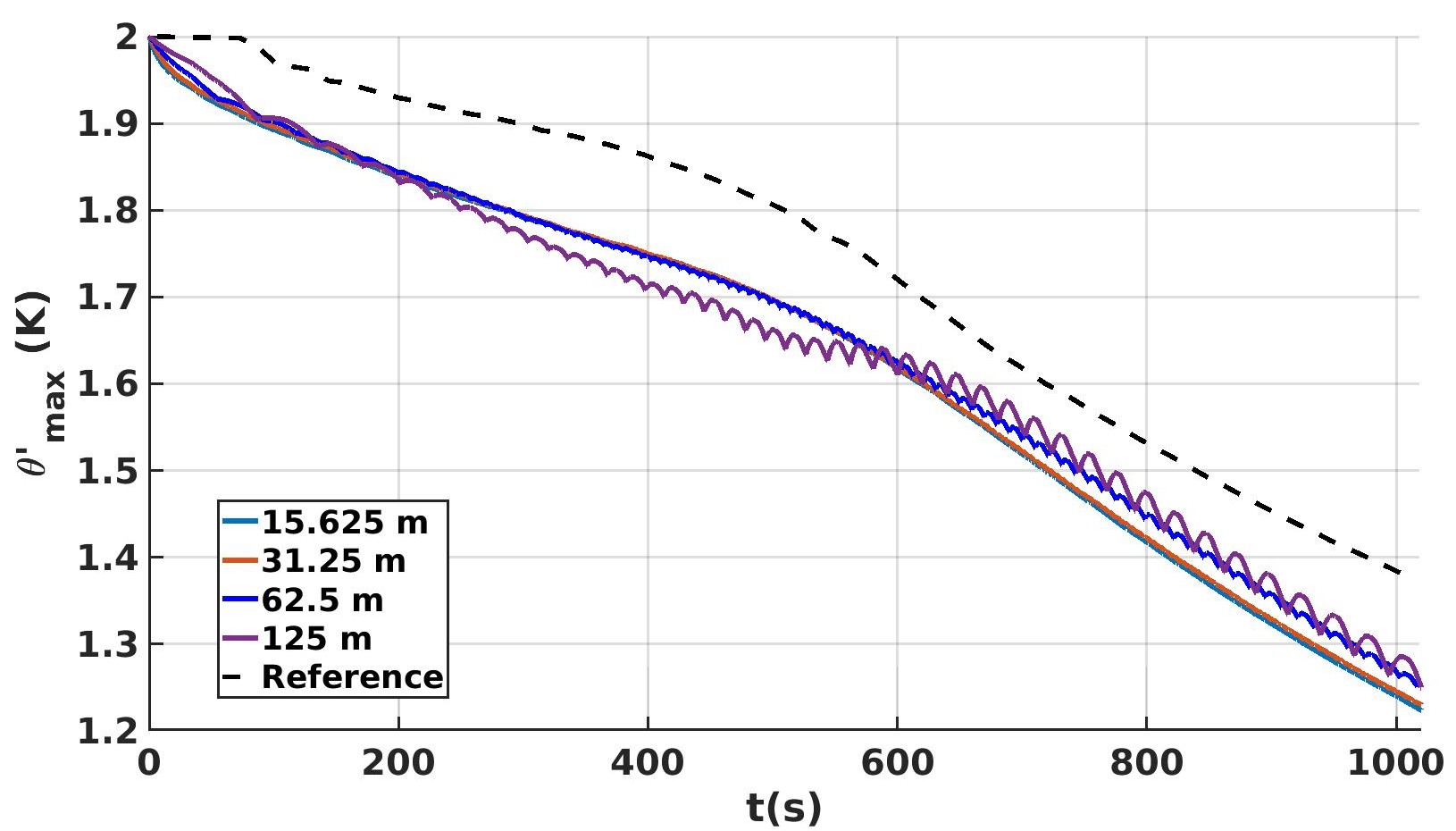}
    \end{overpic}
    \begin{overpic}[width=0.48\textwidth]{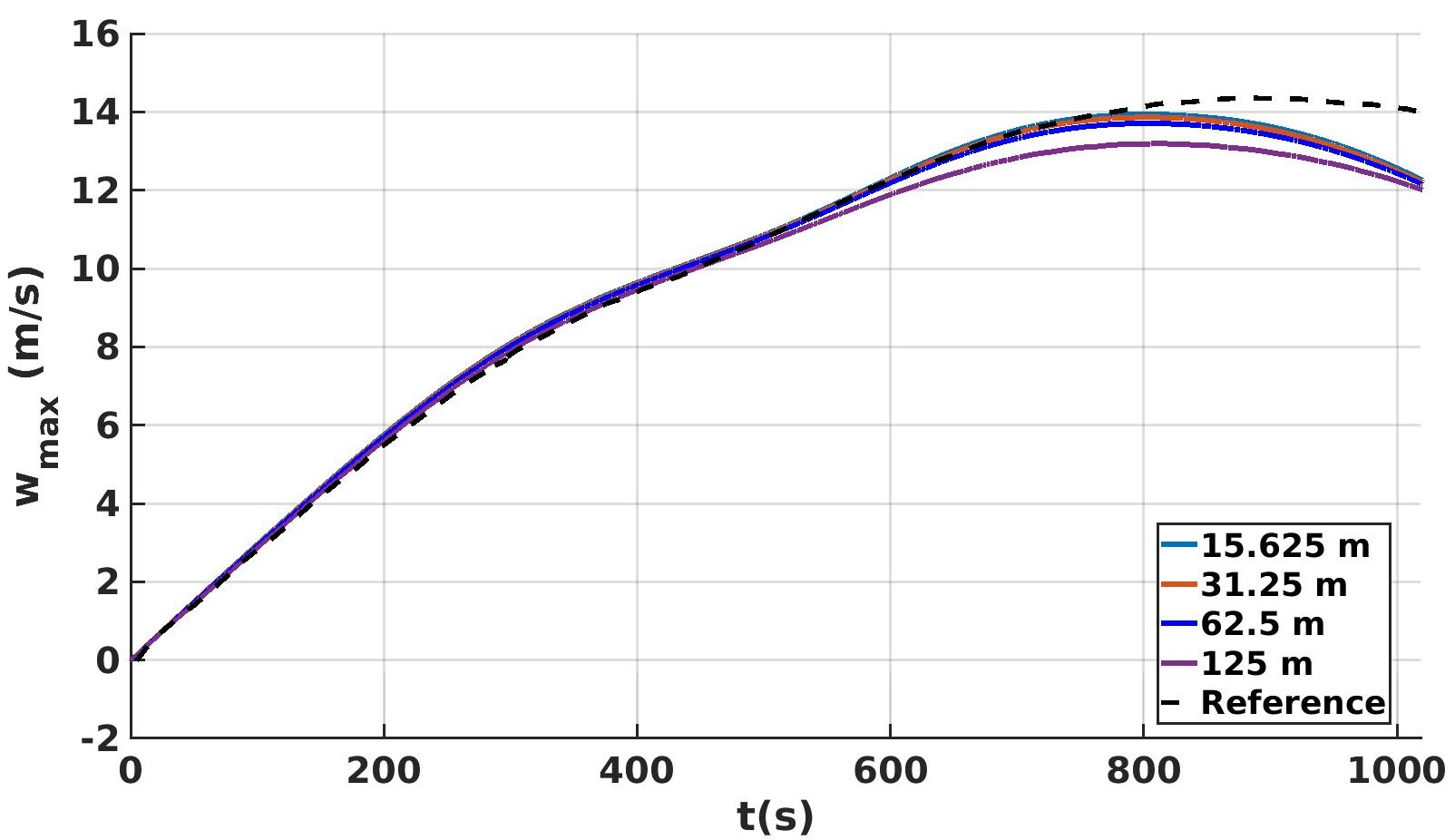}
    \end{overpic} 
     
    \caption{Rising thermal bubble, $a_L$, $\alpha = 1.9$: time evolution of the maximum perturbation of potential temperature $\theta'_{max}$ (left) and the maximum vertical component of the velocity $w_{max}$ (right) computed with all the meshes under consideration. The reference values are taken from \cite{ahmadLindeman2007} and refer to resolution 125 m.}
    \label{fig::TB_Post_max}
\end{figure}

Figure \ref{fig::TB_Post_max} reports a more quantitative comparison. It compares
the time evolution of the maximum perturbation of potential temperature $\theta'_{max}$ and maximum vertical component of the velocity $w_{max}$ computed 
by the EFR algorithm (with $a_L$ and $\alpha = 1.9$) 
against the 
the corresponding results from \cite{ahmadLindeman2007}. We see that the evolution of $\theta'_{max}$ computed with meshes $h = 125$ m is affected by spurious oscillations. Oscillations of small amplitude affect also the $\theta'_{max}$ computed with
mesh $h = 62.5$ m, but they disappear with finer meshes. Since $\theta'_{max}$ and $w_{max}$ computed with meshes $h = 31.25$ m and $h = 15.625$ m are practically overlapped, we chose not to refine the mesh further. 
The ``converged'' $w_{max}$ overlaps with the reference value till about $t = 800$ s, which is a remarkable improvement over our previous results in \cite{GQR_OF_clima}. 
The ``converged'' $\theta'_{max}$ is also closer to the results from \cite{ahmadLindeman2007} than in \cite{GQR_OF_clima}, however there is still some distance between the two curves. 

Table \ref{tab:1} reports the extrema for the vertical velocity $w$ and potential temperature perturbation $\theta'$ at $t = 1020$ s obtained with the EFR algorithm (with $a_L$ and $\alpha = 1.9$), together with the values extracted from the figures in \cite{ahmadLindeman2007}. 
This tables confirms the findings from Fig.~\ref{fig::TB_Post_max}.

\begin{table}[htb!]
\begin{center}
\begin{tabular}{ | c | c | c | c |  c | c | }
\hline
 Type & $h$ (m) &  $w_{min}$ (m/s) & $w_{max}$ (m/s) &  $\theta'_{min} $ (K) & $\theta'_{max} $ (K)\\
 \hline
 Ref.~\cite{ahmadLindeman2007} & 125 & -7.75 & 13.95 & -0.013 & 1.4 \\
 \hline
 $a_L, \alpha = 1.9$ & 125 & -10.35 & 12.01 & -0.012 & 1.23 \\
 \hline
 $a_L, \alpha = 1.9$ & 62.5  & -10.54 & 12.16 & -0.041 & 1.24  \\ 
 \hline
 $a_L, \alpha = 1.9$ & 31.25  & -10.61 & 12.21 & -0.050 & 1.22  \\
 \hline
 $a_L, \alpha = 1.9$ & 15.625 & -10.63 & 12.28 & -0.052 & 1.22 \\
\hline
\end{tabular}
\caption{Rising thermal bubble, $a_L$, $\alpha = 1.9$: minimum and maximum vertical velocity $w$ and potential temperature $\theta'$ at $t=1020$ s compared with the values extracted from the figures in \cite{ahmadLindeman2007}.}\label{tab:1}
\end{center}
\end{table}

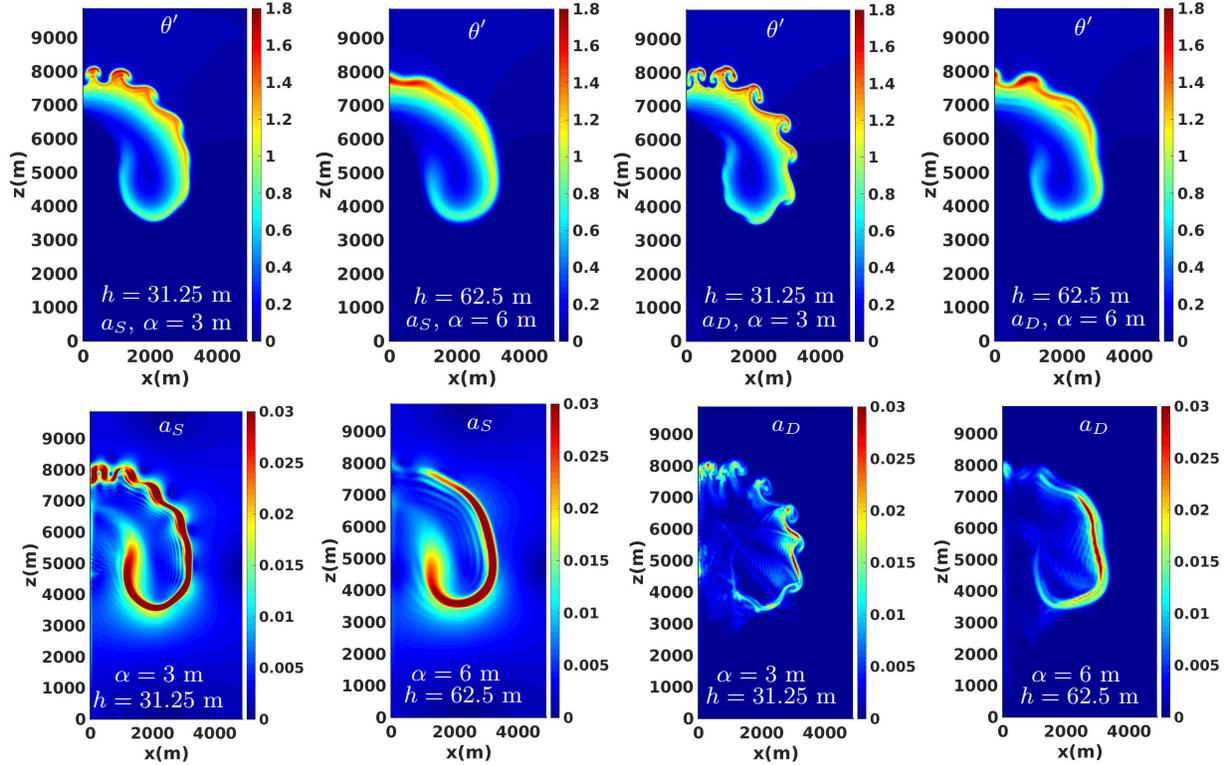
\begin{figure}[htb!]
     \centering
     \hspace*{-0.5cm}
     \begin{overpic}[percent,width=0.245\textwidth]
     {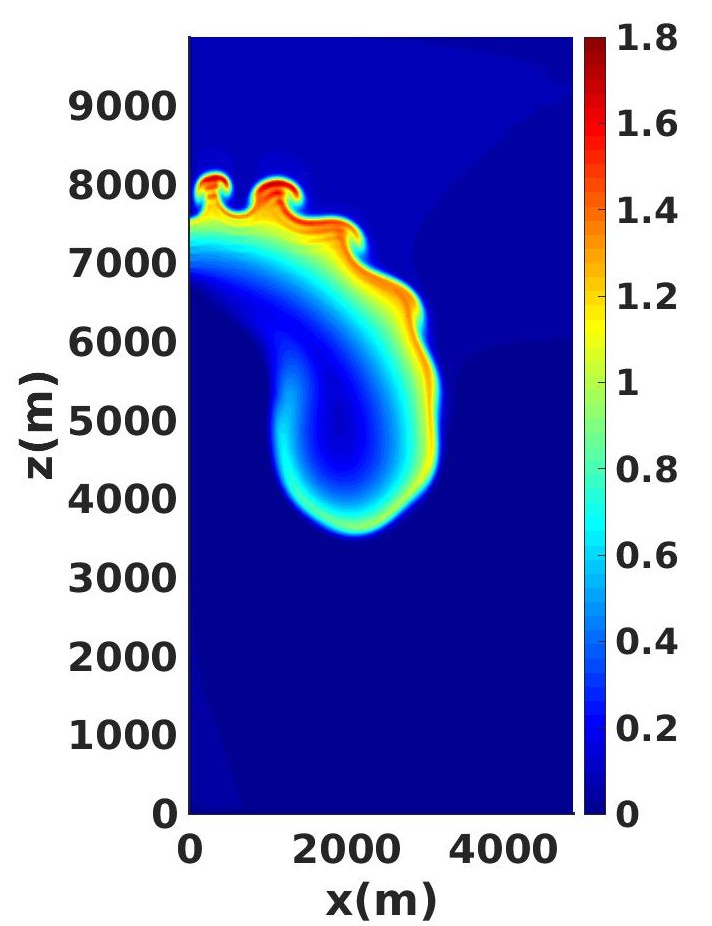}
     \put(40,90){\textcolor{white}{\footnotesize{$\theta'$}}}
     \put(25,22){\textcolor{white}{\footnotesize{$h=31.25$ m}}}
     \put(25,15){\textcolor{white}{\footnotesize{$a_S$, $\alpha=3$ m}}}
    \end{overpic}
    \begin{overpic}[percent,width=0.235\textwidth]{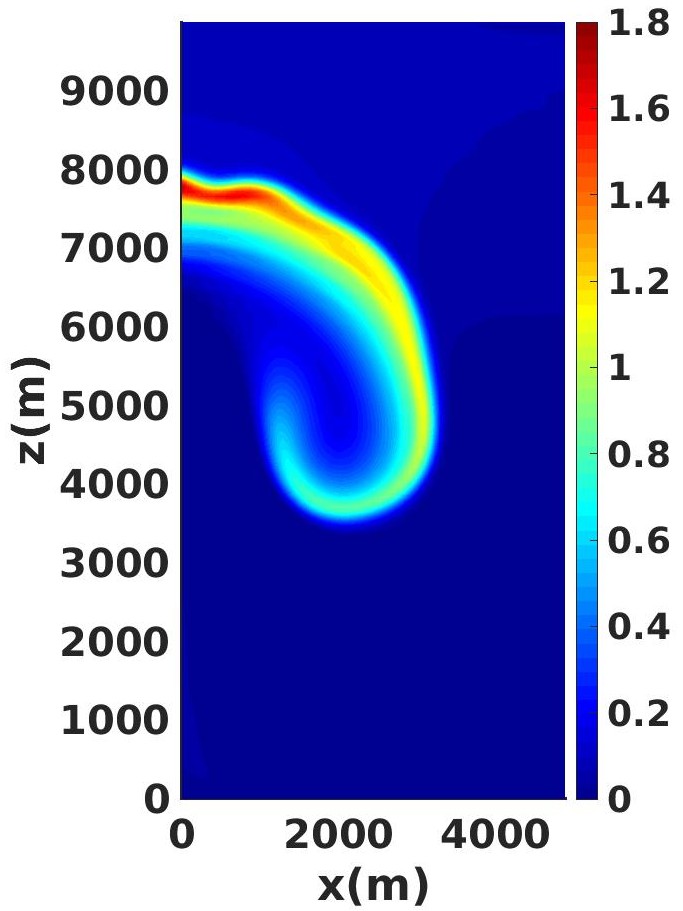}
    \put(40,90){\textcolor{white}{\footnotesize{$\theta'$}}}
     \put(26,22){\textcolor{white}{\footnotesize{$h=62.5$ m}}}
     \put(24,16){\textcolor{white}{\footnotesize{$a_S$, $\alpha=6$ m}}}
    \end{overpic}
     \begin{overpic}[percent,width=0.24\textwidth]
    {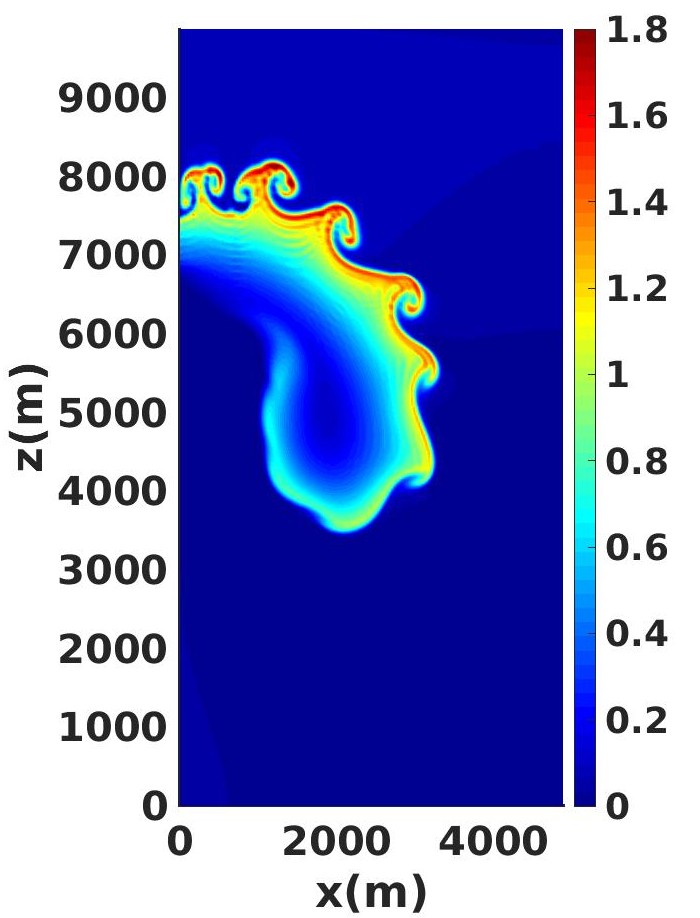}
    \put(40,90){\textcolor{white}{\footnotesize{$\theta'$}}}
     \put(24,22){\textcolor{white}{\footnotesize{$h=31.25$ m}}}
     \put(23.5,16){\textcolor{white}{\footnotesize{$a_D$, $\alpha=3$ m}}}
    \end{overpic} 
    \begin{overpic}[percent,width=0.24\textwidth]{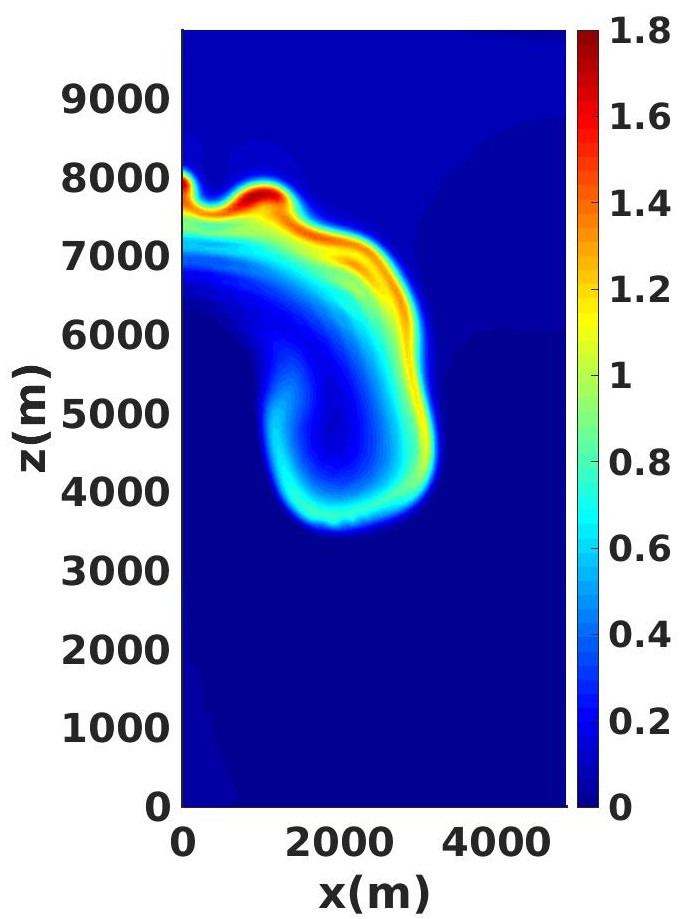}
    \put(40,90){\textcolor{white}{\footnotesize{$\theta'$}}}
     \put(24,22){\textcolor{white}{\footnotesize{$h=62.5$ m}}}
     \put(24,16){\textcolor{white}{\footnotesize{$a_D$, $\alpha=6$ m}}}
    \end{overpic} \\
     \begin{overpic}[percent,width=0.24\textwidth]
     {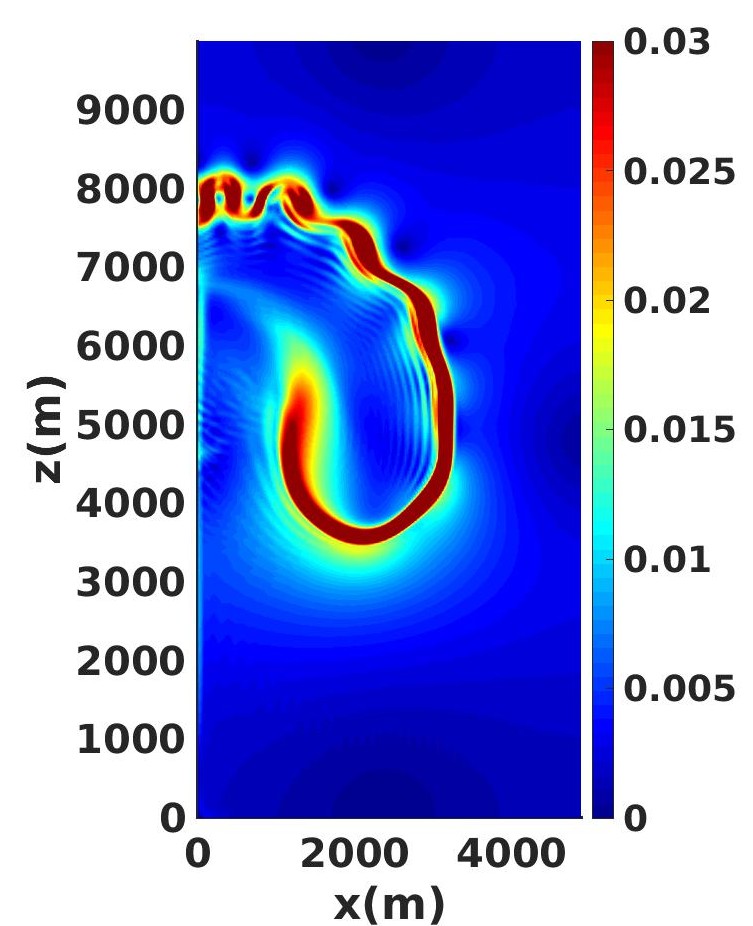}
     \put(40,90){\textcolor{white}{\footnotesize{$a_S$}}}
     \put(22,15){\textcolor{white}{\footnotesize{$h=31.25$ m}}}
     \put(28,22){\textcolor{white}{\footnotesize{$\alpha=3$ m}}}
    \end{overpic}
    \begin{overpic}[percent,width=0.24\textwidth]{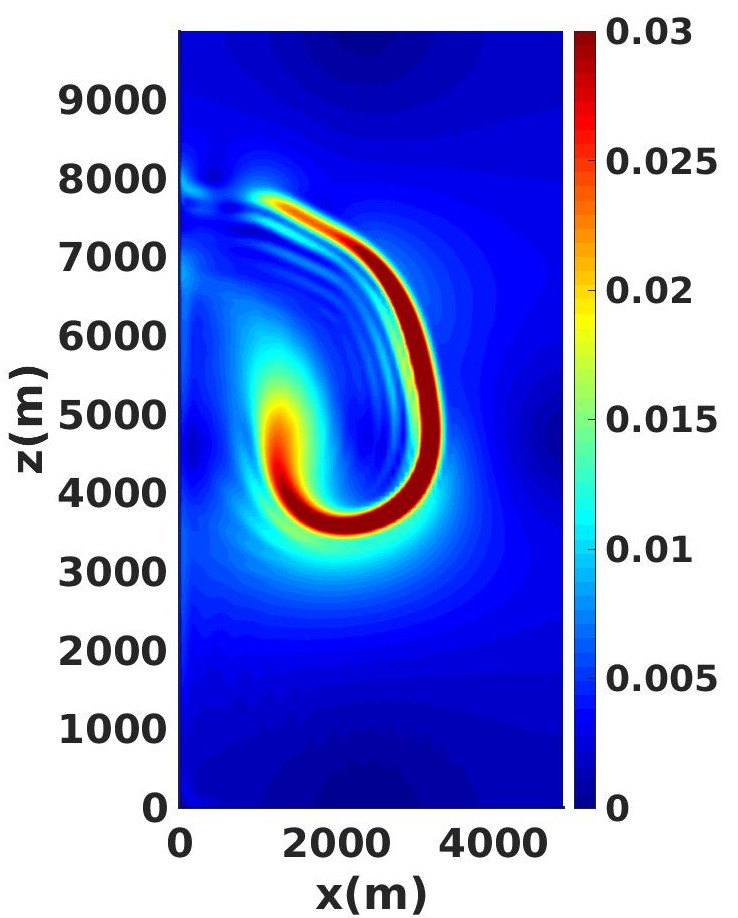}
    \put(40,90){\textcolor{white}{\footnotesize{$a_S$}}}
     \put(23,16){\textcolor{white}{\footnotesize{$h=62.5$ m}}}
     \put(25,22){\textcolor{white}{\footnotesize{$\alpha=6$ m}}}
    \end{overpic}
    \begin{overpic}[percent,width=0.24\textwidth]
    {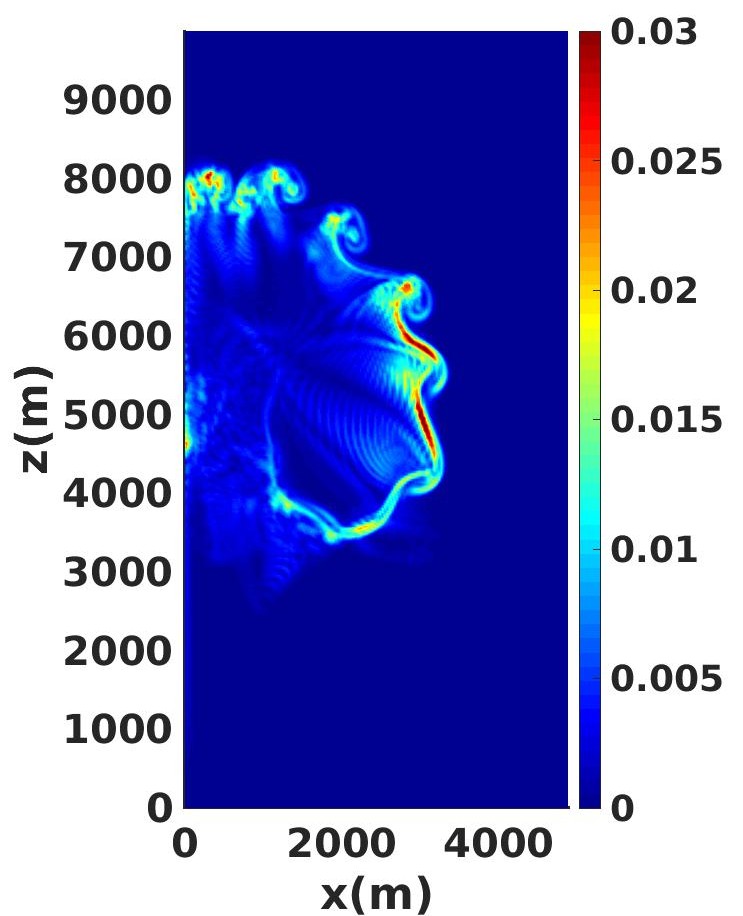}
    \put(40,90){\textcolor{white}{\footnotesize{$a_D$}}}
     \put(22,16){\textcolor{white}{\footnotesize{$h=31.25$ m}}}
     \put(25,22){\textcolor{white}{\footnotesize{$\alpha=3$ m}}}
    \end{overpic}
    \begin{overpic}[percent,width=0.24\textwidth]{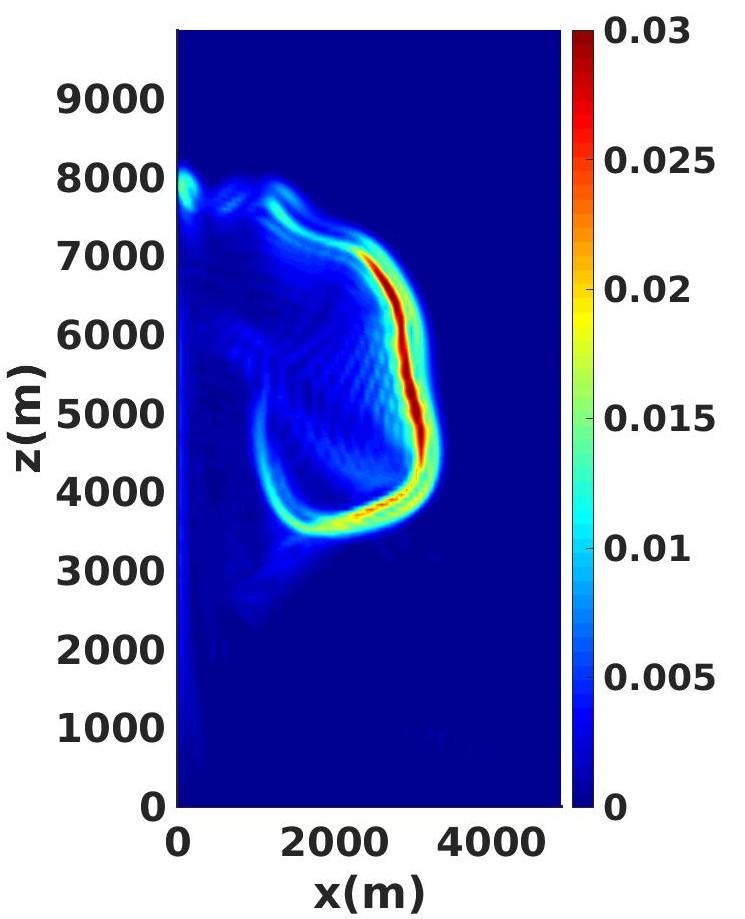}
    \put(40,90){\textcolor{white}{\footnotesize{$a_D$}}}
     \put(24,16){\textcolor{white}{\footnotesize{$h=62.5$ m}}}
     \put(27,22){\textcolor{white}{\footnotesize{$\alpha=6$ m}}}
    \end{overpic}
    \caption{Rising thermal bubble: perturbation of potential temperature (top row) and corresponding indicator function (bottom row) at $t = 1020$ s computed with the EFR and $a_S$ (first two columns) and $a_D$ (last two columns) for different mesh sizes.}
    \label{fig::TB_SM_DB}
\end{figure}

Next, we consider indicator functions $a_S$ and $a_D$ and focus on the two intermediate meshes ($h = 31.25$ m and $h = 62.5$ m).  We set the value of $\alpha$ using $C_s = 0.094$ \cite{GQR_OF_clima} and Remark \ref{rem:3}, 
which suggests an order of magnitude for $\alpha$ rather than a strict value. We take $\alpha=3$ m for mesh $h = 31.25$ m. Since Remark \ref{rem:3} suggests a linear dependence of $\alpha$ on the mesh size, we take $\alpha=6$ m for mesh $h = 62.5$ m. 
We note that these values would apply only for $a_S$ but we will use them for $a_D$ too in order to show the differences in the solutions obtained with the two indicator functions. Fig.~\ref{fig::TB_SM_DB} shows the spatial distribution of $\theta'$ and the indicator function at $t = 1020$ s computed with the EFR algorithm and the nonlinear filters. Note that with the nonlinear filters we can capture a larger amount of
vortical structures than with the linear filter (compare Fig.~\ref{fig::TB_SM_DB} with Fig.~\ref{fig:TB_Linear}).
The results computed with $a_S$ and mesh $h = 31.25$ m (Fig.~\ref{fig::TB_SM_DB}, first panel on the top row) agree very well with those obtained with the Smagorinsky model in \cite{GQR_OF_clima} (Fig.~5, left panel).  
On a given mesh, the Rayleigh-Taylor instability at the edge of the bubble is more developed when using $a_D$ instead of $a_S$, which indicates that $a_D$ introduces less artificial viscosity than $a_S$. Recall that
the artificial viscosity introduced by the EFR algorithm \eqref{eq:mubar_h} is proportional to the indicator function. Indeed, the plots on the bottom row of Fig.~\ref{fig::TB_SM_DB} show that 
$a_S$ at $t = 1020$ s has larger values over wider regions than $a_D$.
This means that indicator function $a_D$
is more selective in identifying the regions of the domain where diffusion is needed.

Table \ref{tab:2} reports the extrema for the vertical velocity $w$ and potential temperature perturbation $\theta'$ at $t = 1020$ s obtained with the EFR algorithm, together with the values from \cite{GQR_OF_clima} for the Smagorinski model. 
The data in Table \ref{tab:2} confirm our observation from Fig.~\ref{fig::TB_SM_DB} about $a_S$ vs $a_D$. Indeed, we see that larger extreme values are found with the EFR algorithm and $a_D$. In addition, we see that the Smagorinsky model from \cite{GQR_OF_clima} gives smaller extreme values than EFR algorithm with $a_S$, which seems to be less diffusive. 

\begin{table}[htb!]
\begin{center}
\begin{tabular}{ | c | c | c | c | c |  c | c | }
\hline
Model & $h$ (m) & $\alpha$ (m) & $w_{min}$ (m/s) & $w_{max}$ (m/s) &  $\theta'_{min} $ (K) & $\theta'_{max} $ (K)\\
 \hline
EFR, $a_D$ & 31.25    & 3  & -13.39 &  15.59  & -0.18 & 1.88 \\
\hline
EFR, $a_S$   & 31.25    & 3  & -11.60 &  15.29  & -0.13 & 1.75 \\
\hline
Smagorisnky \cite{GQR_OF_clima} & 31.25 & - & -11.54 & 15.04 & -0.072 & 1.89\\
\hline
EFR, $a_D$ & 62.5  & 6  & -11.3 &  14.88 & -0.12 & 1.78\\
\hline
EFR, $a_S$ & 62.5  & 6  & -10.76 &  13.38  & -0.14 & 1.72\\
\hline 
\end{tabular}
\caption{Rising thermal bubble: minimum and maximum vertical velocity $w$ and potential temperature $\theta'$ at $t=1020$ s computed by the EFR algorithm with different meshes, indicator functions, and values of $\alpha$. For comparison, the table reports the values from \cite{GQR_OF_clima} obtained with the Smagorisnky model.}\label{tab:2}
\end{center}
\end{table}

\begin{figure}[htb!]
\centering
     \begin{overpic}[percent,width=0.235\textwidth]{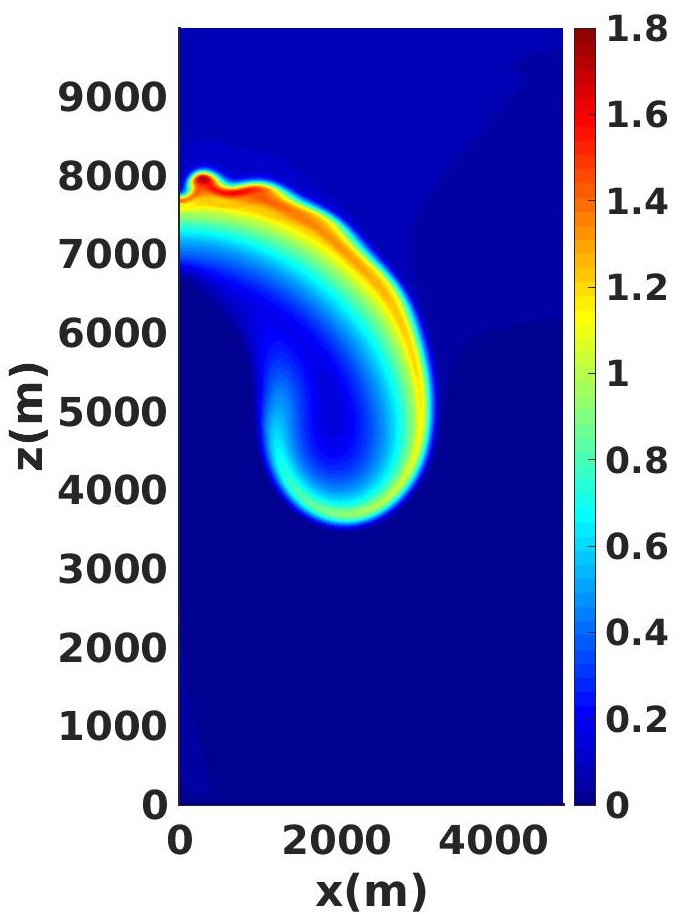}
     \put(25,16){\textcolor{white}{\footnotesize{$\alpha=5$ m}}}
    \end{overpic}
\begin{overpic}[percent,width=0.233\textwidth]{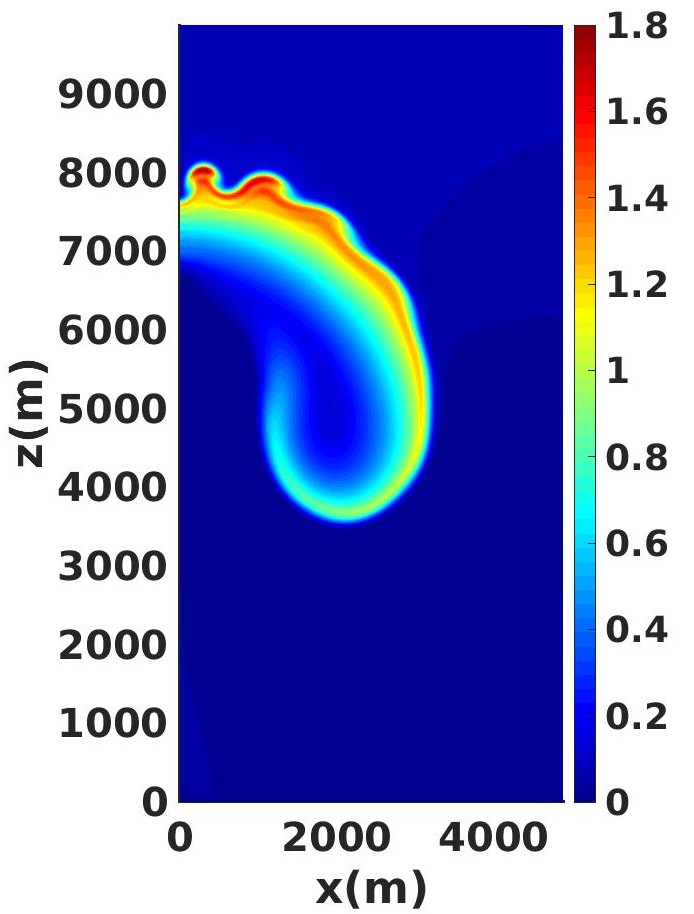}
     \put(32,16){\textcolor{white}{\footnotesize{$\alpha=4$ m}}}
    \end{overpic} 
    \begin{overpic}[percent,width=0.24\textwidth]{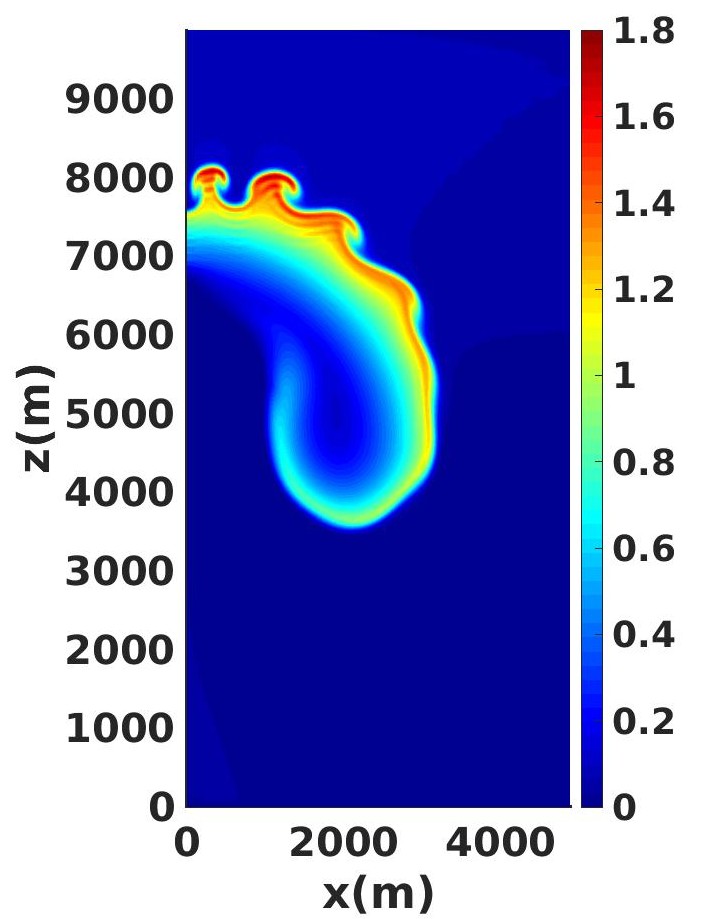}
    \put(32,16){\textcolor{white}{\footnotesize{$\alpha=3$ m}}}
    \end{overpic} 
    \begin{overpic}[percent,width=0.237\textwidth]{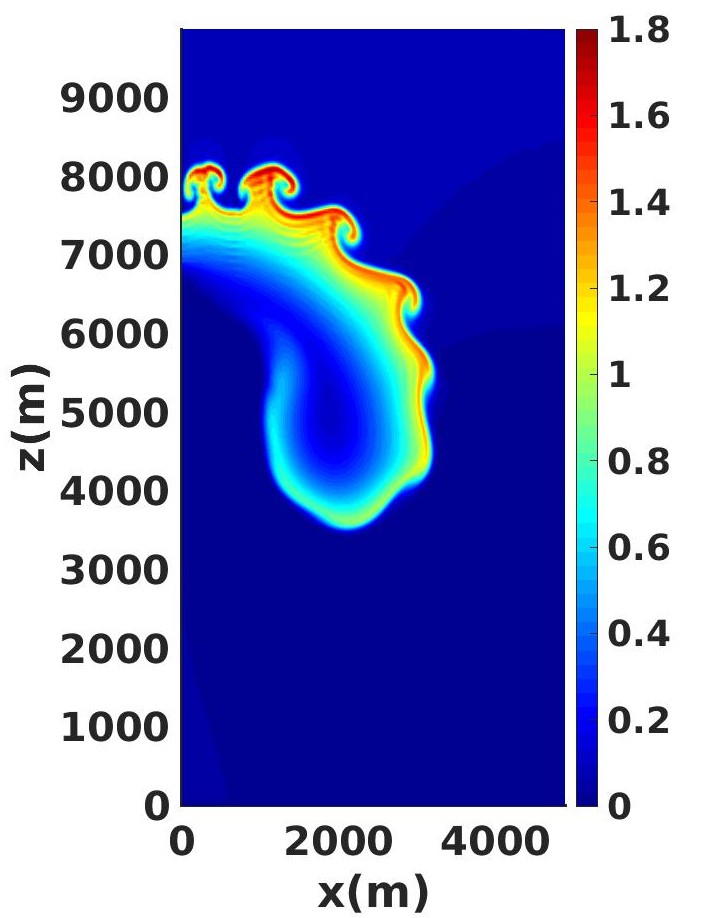}
    \put(32,16){\textcolor{white}{\footnotesize{$\alpha=2$ m}}}
    \end{overpic}
    \caption{Rising thermal bubble, $a_S$: perturbation of potential temperature at $t = 1020$ s computed
    by the EFR algorithm with mesh $h=31.25$ m and (from left to right) $\alpha = 5,4,3,2$ m.}
    \label{fig::TB_31.25_alpha_comp}
\end{figure}

We conclude this section by highlighting the important role played by the filtering radius. Fig.~\ref{fig::TB_31.25_alpha_comp} shows $\theta'$ at $t = 1020$ s computed by the EFR algorithm with $a_S$
and different values of $\alpha$ for mesh $h=31.25$ m. Although the three values of $\alpha$ are all of the same order of magnitude, we see a big difference in the solution. This sensitivity to $\alpha$ can be mitigated by choosing $\chi \neq 1$ (see, e.g., \cite{BQV}) and $\xi \neq 1$.  

\subsection{Density current}\label{sec:DC}

The computational domain in the $xz$-plane is $\Omega=[0,25600]\times[0,6400]~\mathrm{m}^2$ and the time interval of interest is $(0,900]$ s. 
Impenetrable, free-slip boundary conditions are imposed on all the walls. 
The initial density is given by \eqref{eq:rho_wb} with initial potential temperature:
\begin{equation}
\theta^0 = 300 - \frac{15}{2}\left[  1 + \cos(\pi r)\right] ~ \textrm{if $r\leq 1$},\quad\theta^0 = 300
~ \textrm{otherwise},
\label{dcEqn1}
\end{equation}
where $r = \sqrt[]{\left(\frac{x-x_{c}}{x_r}\right)^{2} + \left(\frac{z-z_{c}}{z_r}\right)^{2}}$, with $(x_r,z_r)=(4000, 2000)~{\rm m}$ and $(x_c,z_c) = (0,3000)~\mathrm{m}$. The initial bubble in this test is of cold air. The initial velocity field is zero everywhere and the initial specific enthalpy is given by \eqref{eq:e0}.


We consider uniform, orthogonal meshes with mesh sizes $h =\Delta x$ = $\Delta z  = 200,100,50,25$ m. The time step is set to $\Delta t = 0.1$ s. 
Just like in the case of the warm bubble, we set $\chi = \xi = 1$ in the EFR algorithm. 

\begin{figure}
    \centering
    \begin{overpic}[width=0.49\textwidth]{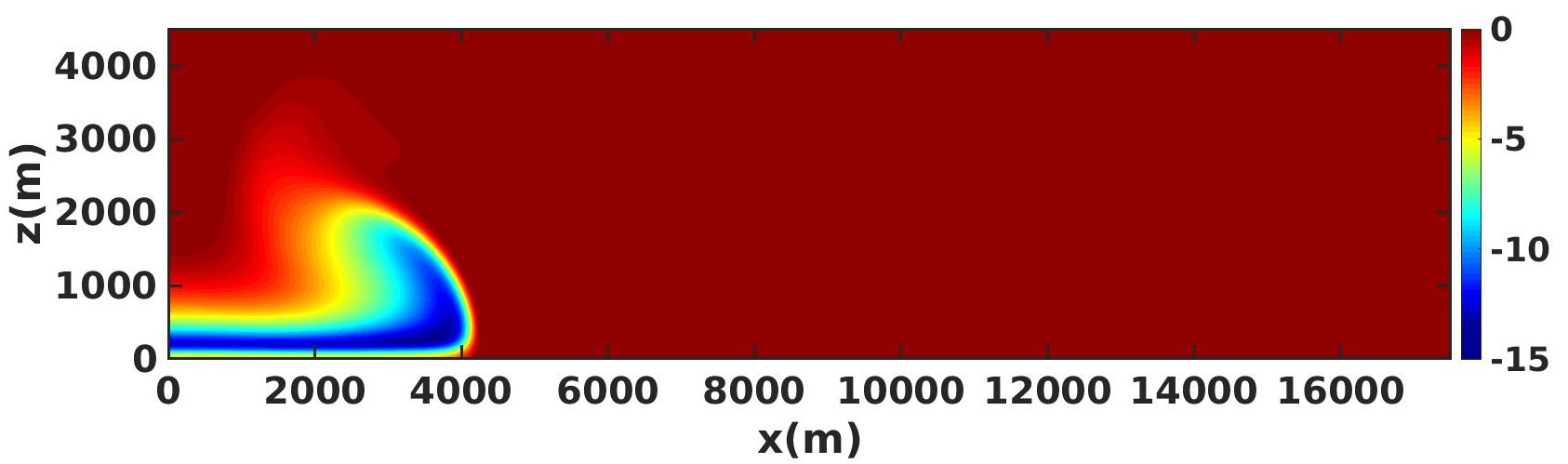}
    \put(160,50){\textcolor{white}{$t = 300$ s}}
    \label{fig:refined_300}
    \end{overpic}
    \begin{overpic}[width=0.49\textwidth]{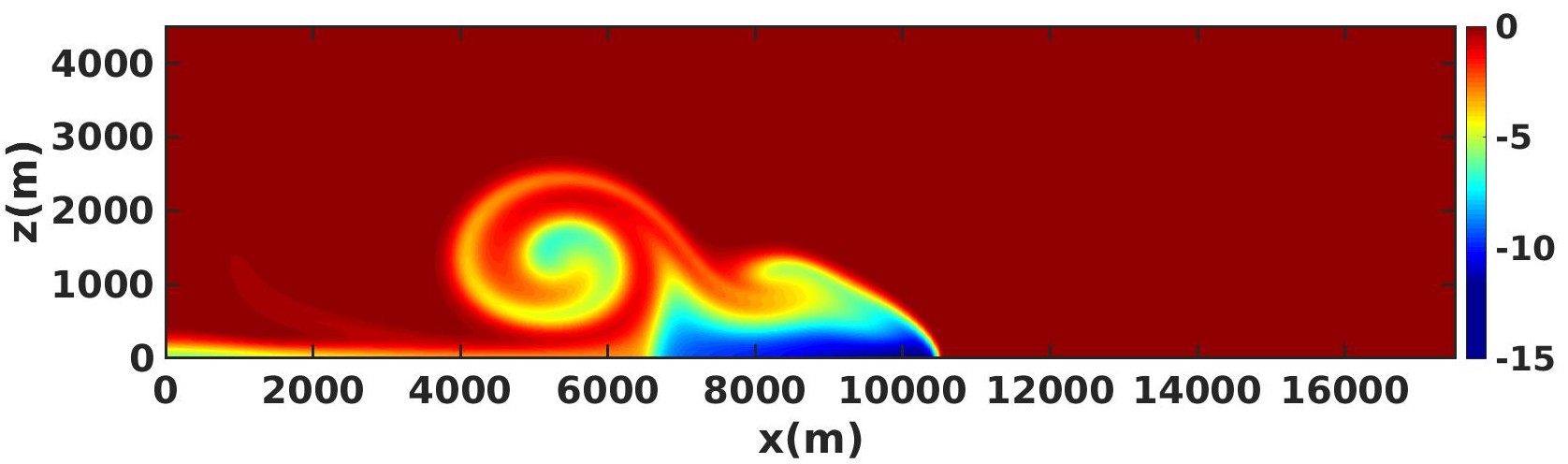}
    \put(160,50){\textcolor{white}{$t = 600$ s}}
    \label{fig:refined_600}
    \end{overpic} \\ \vspace{0.2cm}
    \begin{overpic}[width=0.49\textwidth]{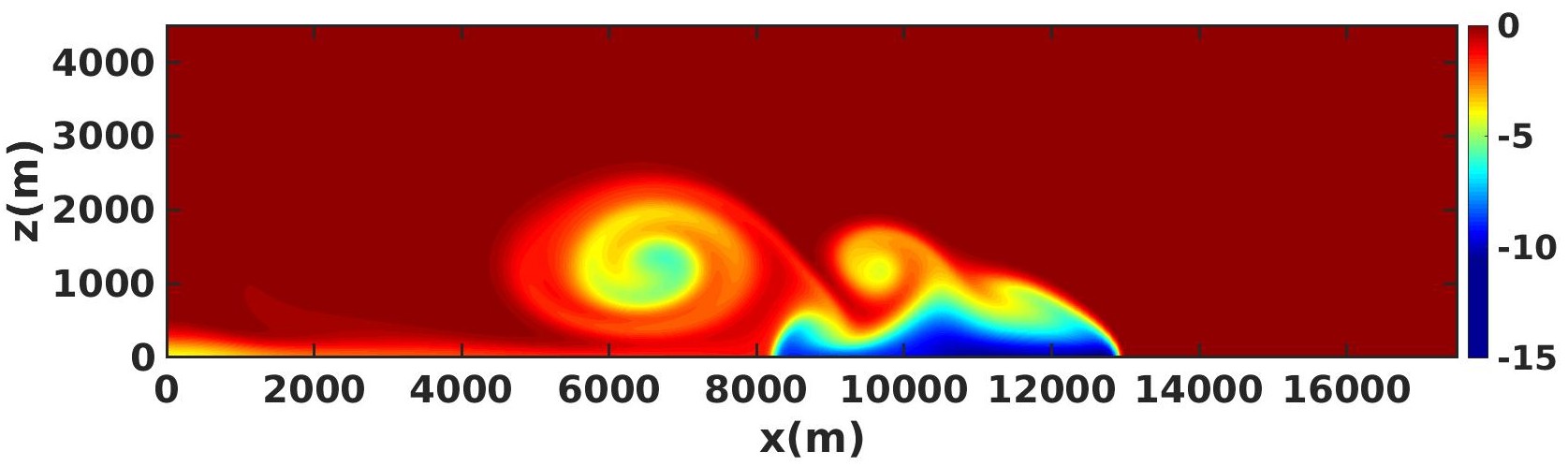}
    \put(160,50){\textcolor{white}{$t = 750$ s}}
    \label{fig:refined_750}
    \end{overpic}
    \begin{overpic}[width=0.49\textwidth]{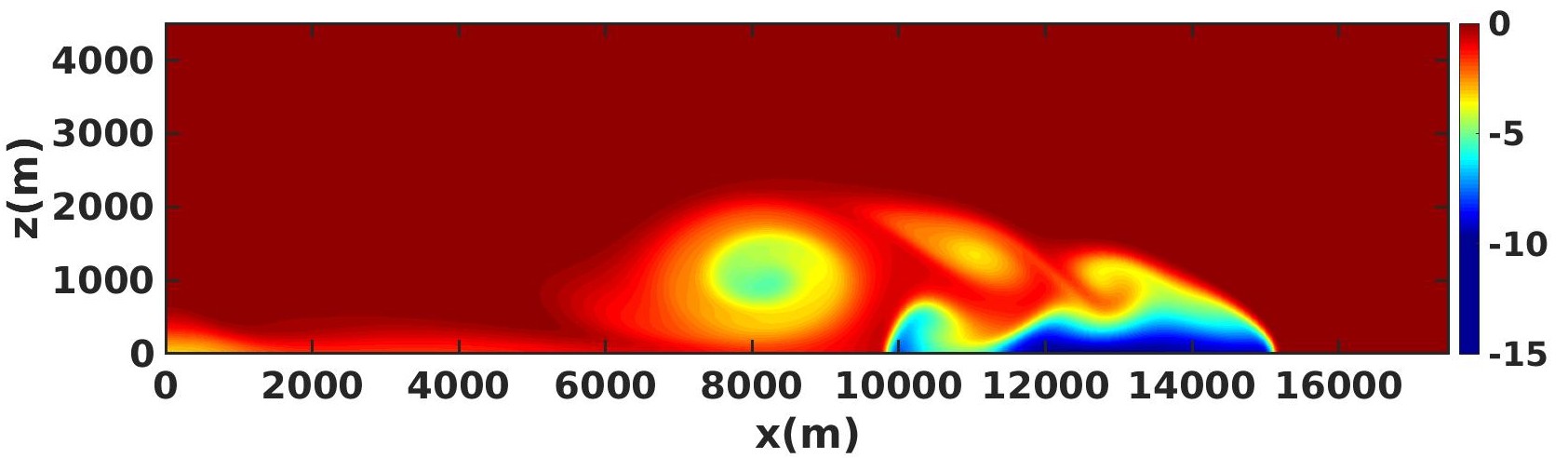}
    \put(160,50){\textcolor{white}{$t = 900$ s}}
    \label{fig:refined_900}
    \end{overpic} 
    \caption{Density current, $a_L$, $\alpha = 2.7$: time evolution of potential temperature fluctuation $\theta'$ computed with mesh $h = 25$ m.}
        \label{fig::linear_25}
\end{figure}

\begin{figure}
    \centering
    \begin{overpic}[width=.49\textwidth,grid=false]{Linear-Results/AV75_COMP/DC25_2.74_LIN_900.jpg}
    \put(150,50){\textcolor{white}{$h =25$ m}}
    \end{overpic} \\ \vspace{0.2cm}
    \begin{overpic}[width=.49\textwidth]{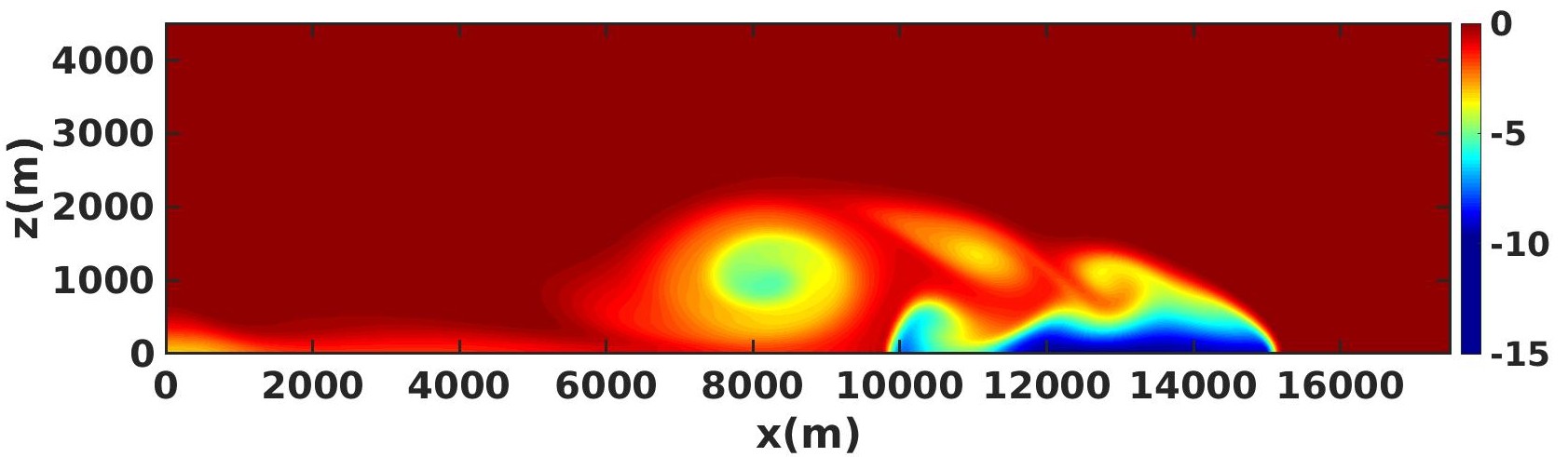}
    \put(150,50){\textcolor{white}{$h=50$ m}}
    \end{overpic} \\ \vspace{0.2cm}
    \begin{overpic}[width=.49\textwidth]{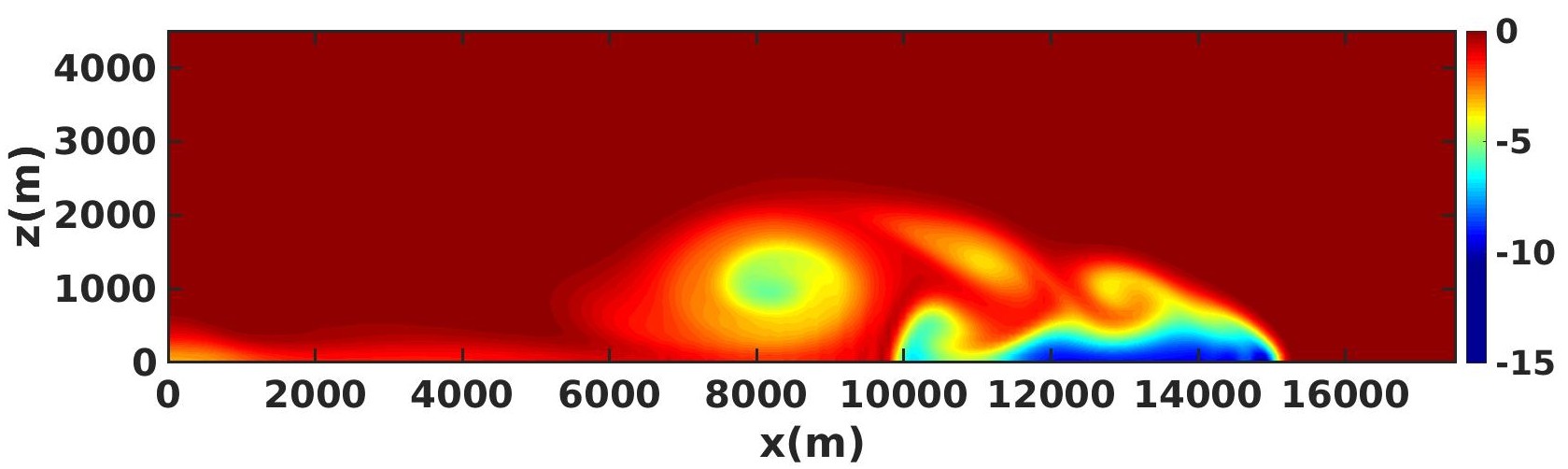}
    \put(150,50){\textcolor{white}{$h=100$ m}}
    \end{overpic} \\ \vspace{0.2cm}
    \begin{overpic}[width=.49\textwidth]{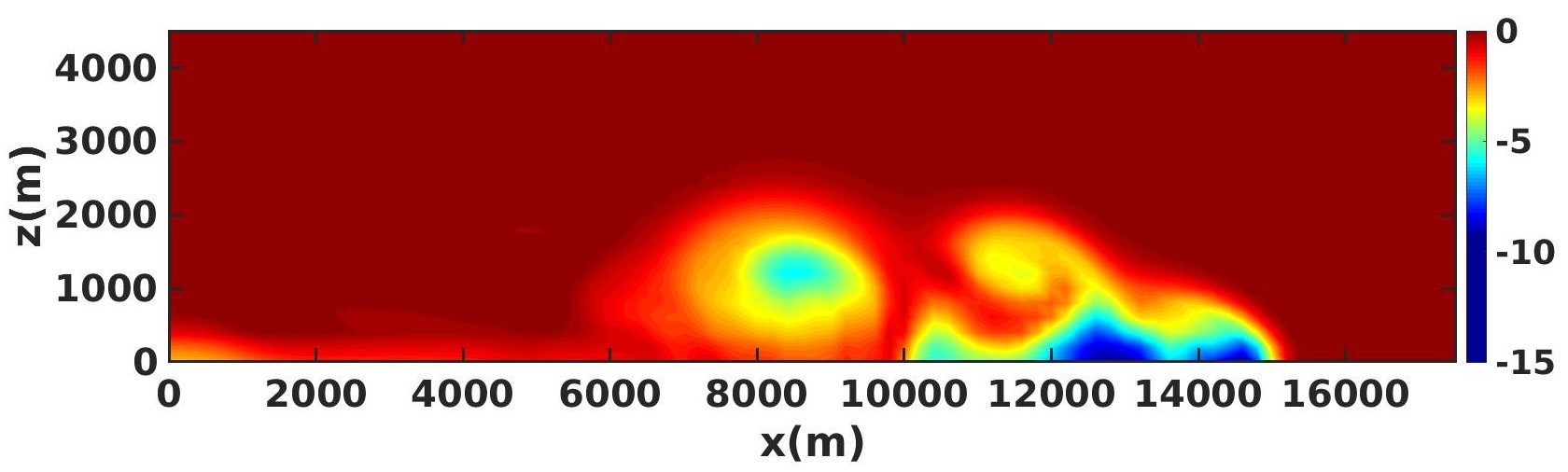}
    \put(150,50){\textcolor{white}{$h=200$ m}}
    \end{overpic} 
\caption{Density current, $a_L$, $\alpha = 2.7$: potential temperature fluctuation $\theta'$ computed at $t = 900$ s with meshes $h = 25, 50, 100, 200, 400~$m. The mesh size is increasing from top to bottom.}
\label{fig:linear_sharing}
\end{figure}

We start again with the linear filter, i.e., we take $a_L$ \eqref{eq:a_lin} as indicator function, because it allows us to make a direct comparison with the results obtained by setting $\mu_a = 75$ and $Pr = 1$ in \eqref{eq:mom_LES}-\eqref{eq:ent_LES}
\cite{GQR_OF_clima,strakaWilhelmson1993}.
To introduce the same amount of artificial viscosity with the EFR algorithm and $a_L$, we use \eqref{eq:mubar_h} and get $\alpha = 2.7$ m by using the minimum density in the computational domain. 
Figure \ref{fig::linear_25} shows $\theta'$ computed with this value of $\alpha$ and mesh 
$h = 25$ m (i.e., the finest mesh among those considered) at $t = 300, 600, 750, 900$ s. 
We observe very good agreement with the results reported in Fig.~1 of \cite{strakaWilhelmson1993}, which were obtained with the same resolution.
In order to understand the behavior of the liner filter as the mesh size is varied, we report $\theta'$ computed at $t = 900$ s with all the meshes mentioned above in Fig.~\ref{fig:linear_sharing}.  
We observe that the EFR algorithm with $a_L$ and $\alpha = 2.7$ m does not introduce sufficient artificial diffusion to stabilize the solution with the coarsest mesh we consider (i.e., $h = 200$ m). For all the other meshes though, we see the emergence of a clear three-rotor structure when the mesh is refined.
Also the results in Fig.~\ref{fig:linear_sharing} are in very good agreement with those reported in the literature. See, e.g., \cite{ahmadLindeman2007,giraldo_2008,GQR_OF_clima,marrasEtAl2013a,marrasNazarovGiraldo2015,strakaWilhelmson1993}.

\begin{figure}[htb!]
    \centering
    \begin{subfigure}{0.5\textwidth}
        \includegraphics[width=\linewidth]{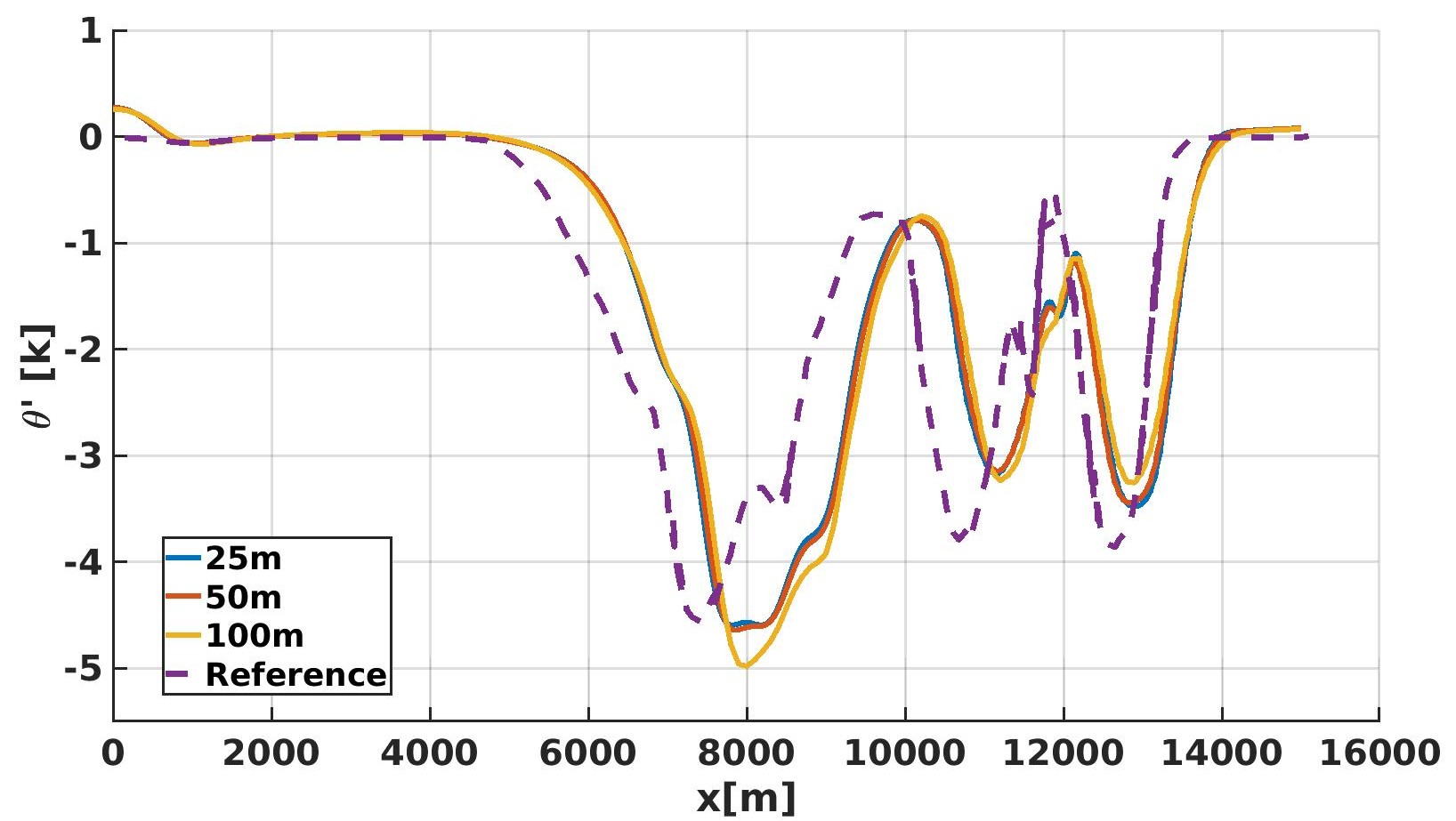}
        \label{fig:theta_Linear_comparison}
    \end{subfigure}
    \caption{Density current, $a_L$, $\alpha = 2.7$: potential temperature perturbation $\theta'$ at $t = 900$ s
along the horizontal direction at a height of $z = 1200$ m for meshes $h = 100, 50, 25$ m compared against data from \cite{giraldo_2008} (denoted as ``Reference'' and referred to resolution 25 m).}
    \label{fig:thetaLine_}
\end{figure}

For a quantitative comparison, we consider the 
potential temperature perturbation $\theta'$ at $t = 900$ s
along the horizontal direction at height $z = 1200$ m.
Fig.~\ref{fig:thetaLine_} displays a comparison between the results given by the EFR model with $a_L$ and $\alpha = 2.7$ for meshes $h = 100, 50, 25$ m. 
We see that the curves associated to meshes $h = 50, 25$ m are practically superimposed. 
In Fig.~\ref{fig:thetaLine_}, 
we report also the results from \cite{giraldo_2008}, which were obtained by setting a constant artificial viscosity (i.e., $\mu_a = 75$) and using a spectral element method. Such results are labeled as ``Reference'' and refer to resolution 25 m. We observe that our results are slightly out of phase with respect to the reference data. Each dip in $\theta'$  
in  Fig.~\ref{fig:thetaLine_} corresponds to a recirculation in Fig.~\ref{fig:linear_sharing}, top three panels. So, from Fig.~\ref{fig:thetaLine_} we learn that in our simulations the front is faster than in the simulations from \cite{giraldo_2008}. For this reason, Table \ref{tab:3} reports the front location (defined as the location on the ground where $\theta'$ = -1 K) at $t = 900$ s
obtained with EFR and $a_L$ and compares it with 
the data in Table 4 of \cite{strakaWilhelmson1993}. The data in \cite{strakaWilhelmson1993} refer to the model with constant artificial viscosity (i.e., $\mu_a = 75$) and 14 different numerical approaches. We note that our results fall well within the values from
\cite{strakaWilhelmson1993}. Hence, we attribute the difference with the Reference in Fig.~\ref{fig:thetaLine_} to the use of different numerical methods.

\begin{table}[htb!]
\begin{center}
\begin{tabular}{ | c | c |  c | }
\hline
 Method & $h$ (m) & Front Location (m)  \\
 \hline
  EFR, $a_L$ & 25 &  15170  \\
 \hline
 EFR, $a_L$ & 50 & 15190  \\ 
 \hline
 EFR, $a_L$ & 100 & 15210  \\
 \hline
  Ref.~\cite{strakaWilhelmson1993} & (25, 200) & (14533,17070)
 \\
  \hline
\end{tabular}
\caption{Density current, $a_L$, $\alpha = 2.7$: our results for the front location at $t = 900$ s obtained with different meshes compared against results reported in 
\cite{strakaWilhelmson1993}.
For reference \cite{strakaWilhelmson1993}, we provide the range of mesh sizes and front location values obtained with different methods.}\label{tab:3}
\end{center}
\end{table}

\begin{figure}[htb!]
     \centering
     \begin{overpic}[width=0.48\textwidth]{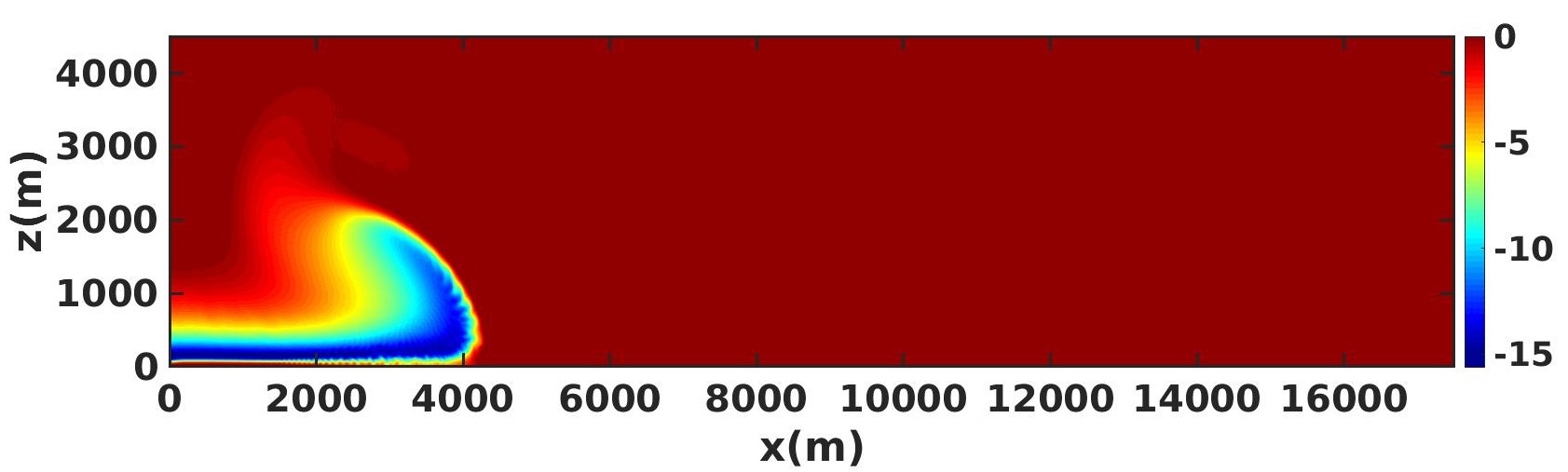}
     \label{Smago_125_300}
    \put(90,50){\textcolor{white}{\footnotesize{$h = 12.5$ m, $\alpha = 4$, $t = 300$ s}}}
    \end{overpic} 
    \begin{overpic}[width=0.48\textwidth]{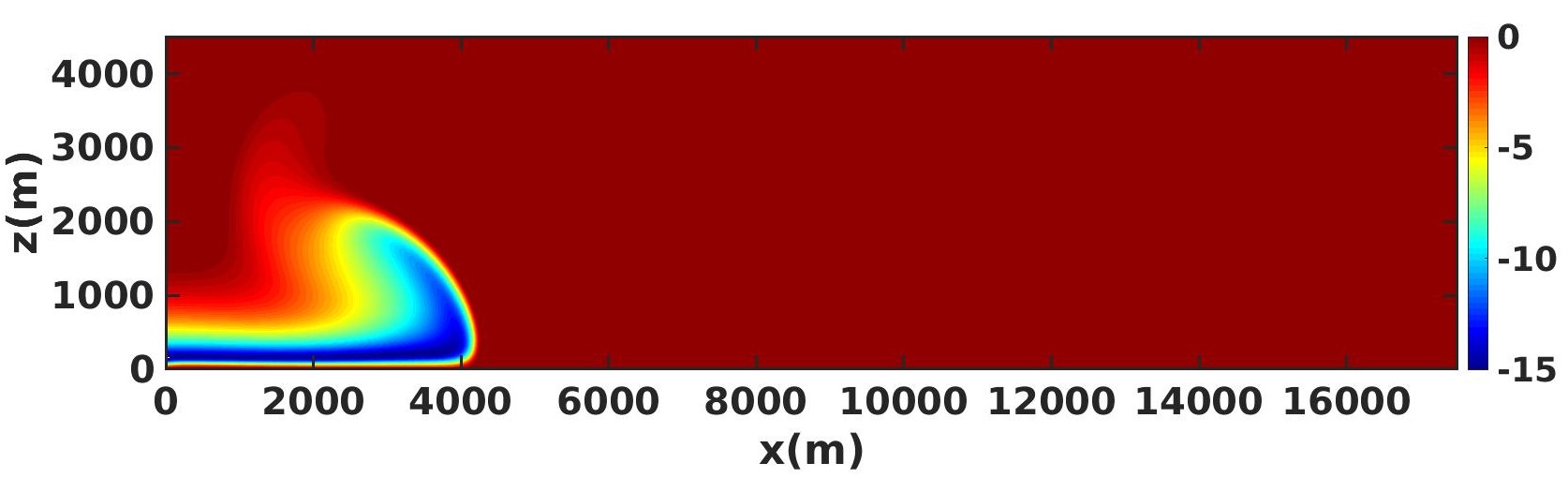}
     \label{fig:Smago_25_300}
    \put(90,50){\textcolor{white}{\footnotesize{$h = 25$ m, $\alpha = 8$, $t = 300$ s}}}
    \end{overpic}
    \\
    \begin{overpic}[width=0.48\textwidth]{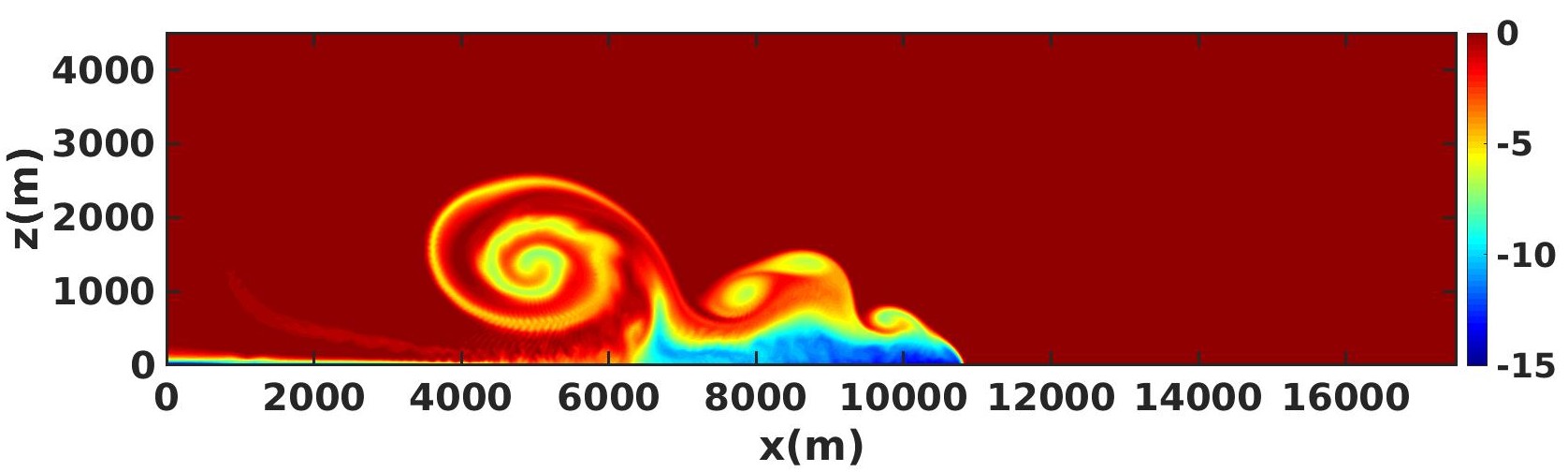}
    \label{Smago_125_600}
    \put(90,50){\textcolor{white}{\footnotesize{$h = 12.5$ m, $\alpha = 4$, $t = 600$ s}}}
    \end{overpic} 
        \begin{overpic}[width=0.48\textwidth]{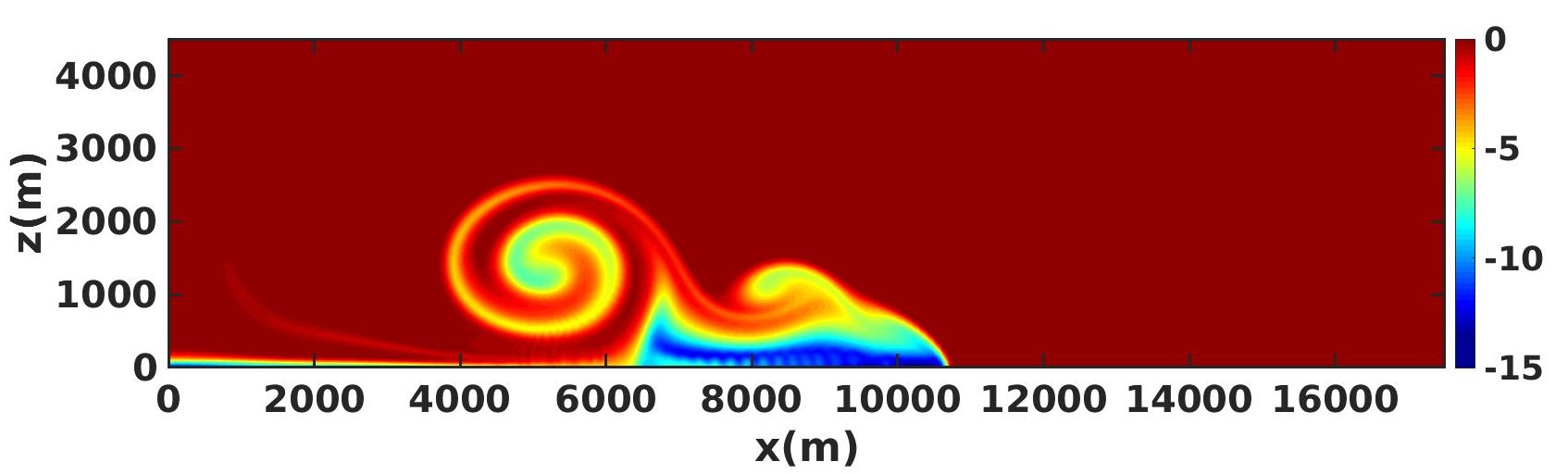}
    \label{fig:Smago_25_600}
    \put(90,50){\textcolor{white}{\footnotesize{$h = 25$ m, $\alpha = 8$, $t = 600$ s}}}
    \end{overpic}
    \\ 
    \begin{overpic}[width=0.48\textwidth]{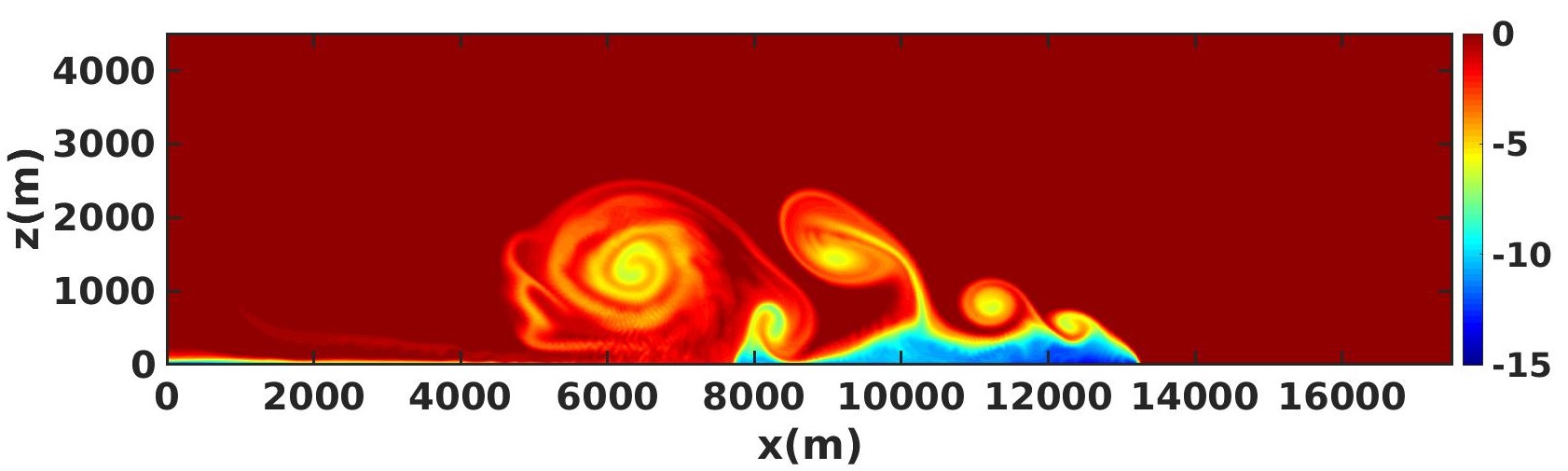}
    \label{Smago_125_750}
    \put(90,50){\textcolor{white}{\footnotesize{$h = 12.5$ m, $\alpha = 4$, $t = 750$ s}}}
    \end{overpic} 
        \begin{overpic}[width=0.48\textwidth]{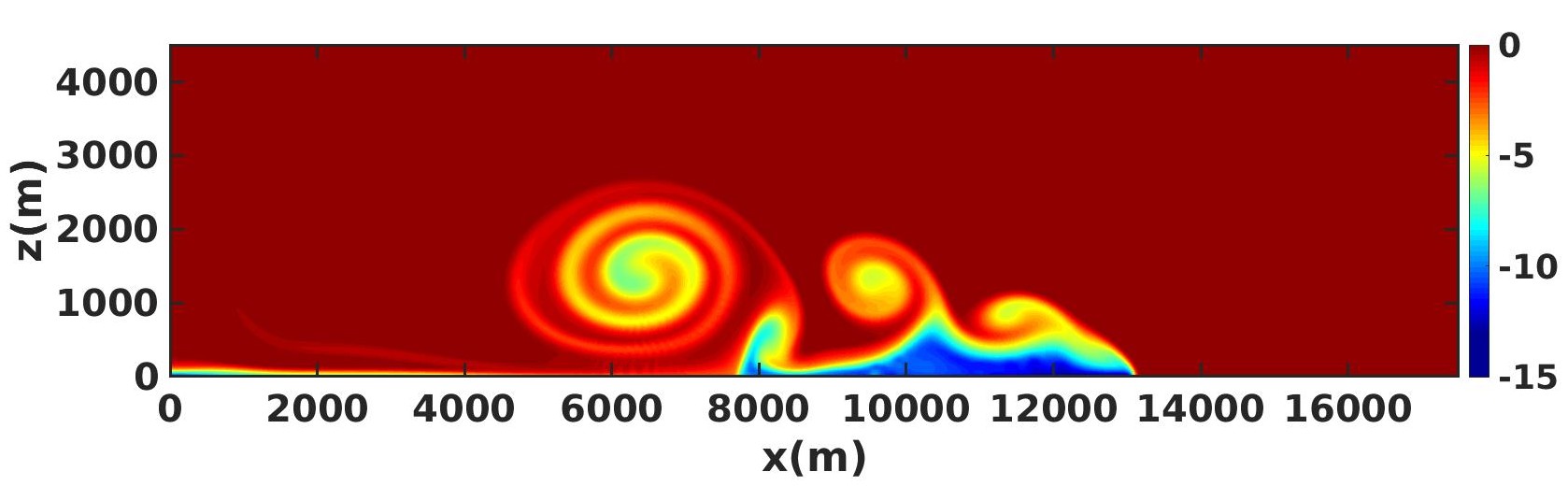}
    \label{fig:Smago_25_750}
    \put(90,50){\textcolor{white}{\footnotesize{$h = 25$ m, $\alpha = 8$, $t = 750$ s}}}
    \end{overpic}
    \\
    \begin{overpic}[width=0.48\textwidth]{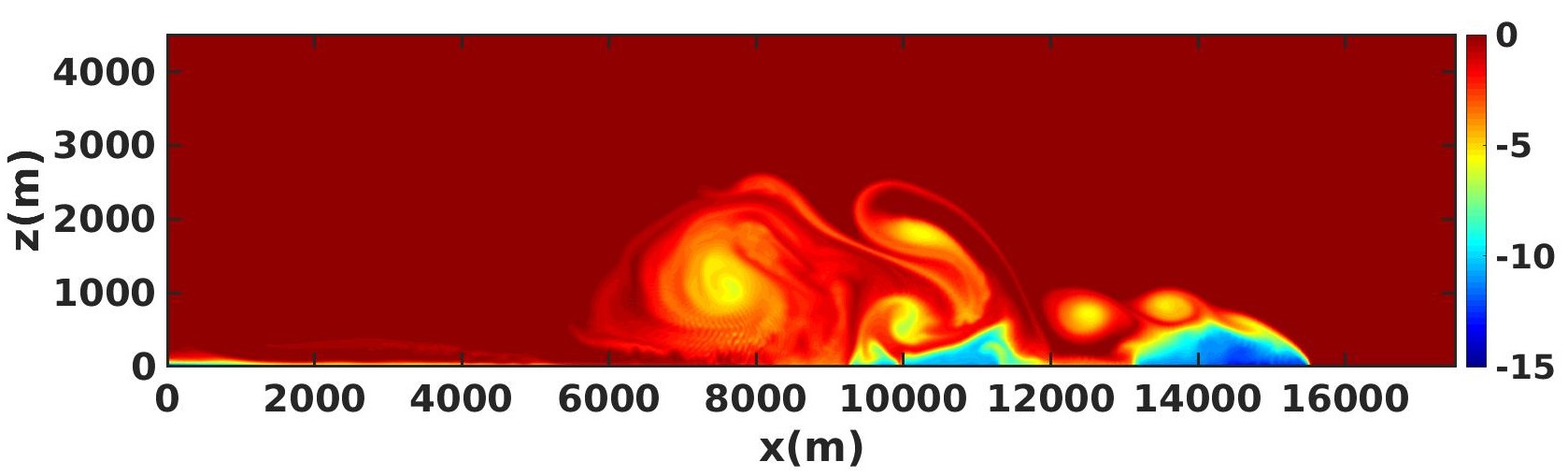}
    \put(90,50){\textcolor{white}{\footnotesize{$h = 12.5$ m, $\alpha = 4$, $t = 900$ s}}}
    \end{overpic} 
        \begin{overpic}[width=0.48\textwidth]{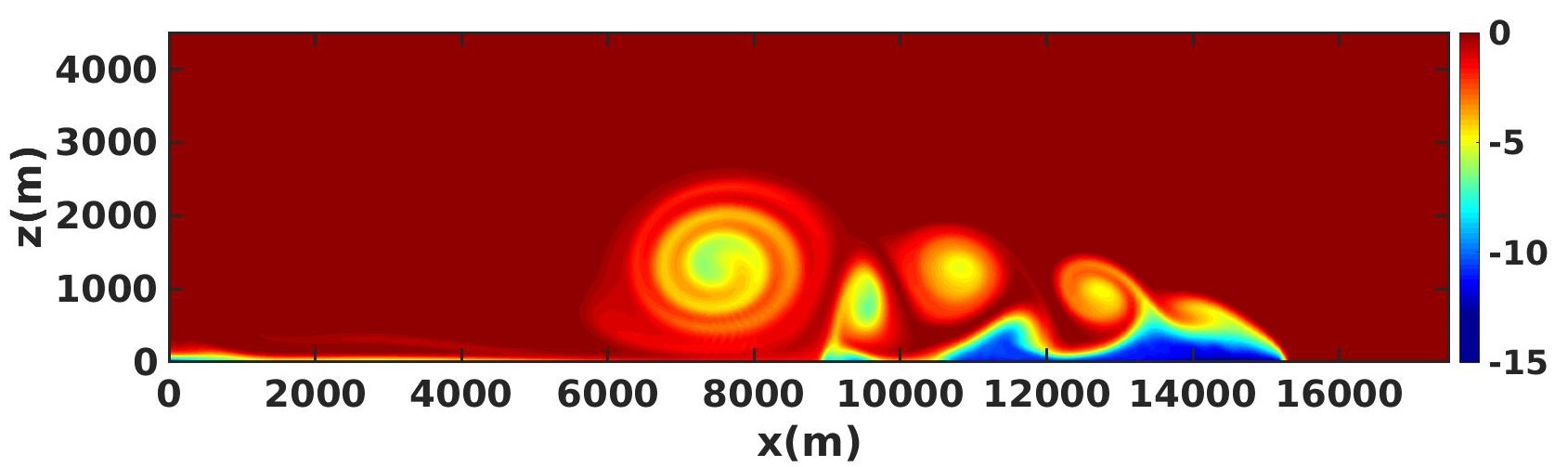}
    \label{fig:Smago_25_900}
    \put(90,50){\textcolor{white}{\footnotesize{$h = 25$ m, $\alpha = 8$, $t = 900$ s}}}
    \end{overpic} 
    \caption{
    Density current, $a_S$:
time evolution of potential temperature fluctuation $\theta'$ computed with mesh $h = 12.5$ m and $\alpha = 4$ (left) and with mesh $h = 25$ m and $\alpha = 8$ (right).}
    \label{fig::SmagoRayModel125_25}
\end{figure}

\begin{figure}[htb!]
    \centering
     \begin{overpic}[width=0.48\textwidth]{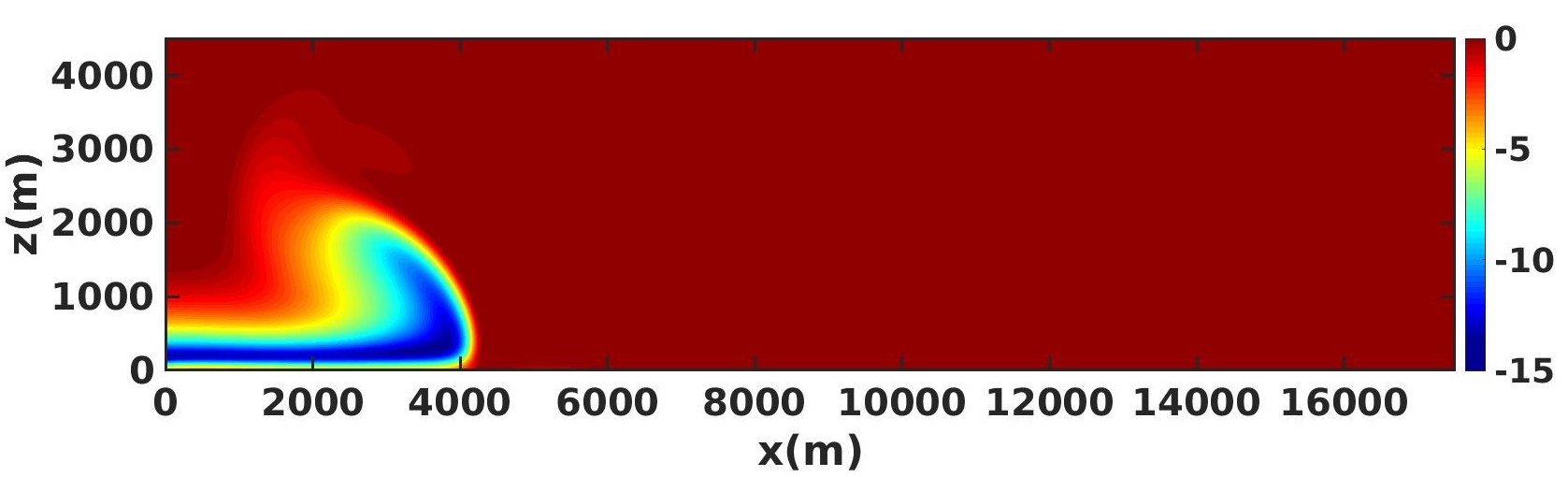}
     \label{fig:Smago_50_300}
    \put(90,50){\textcolor{white}{\footnotesize{$h = 50$ m, $\alpha = 16$, $t = 300$ s}}}
    \end{overpic} 
         \begin{overpic}[width=0.48\textwidth]{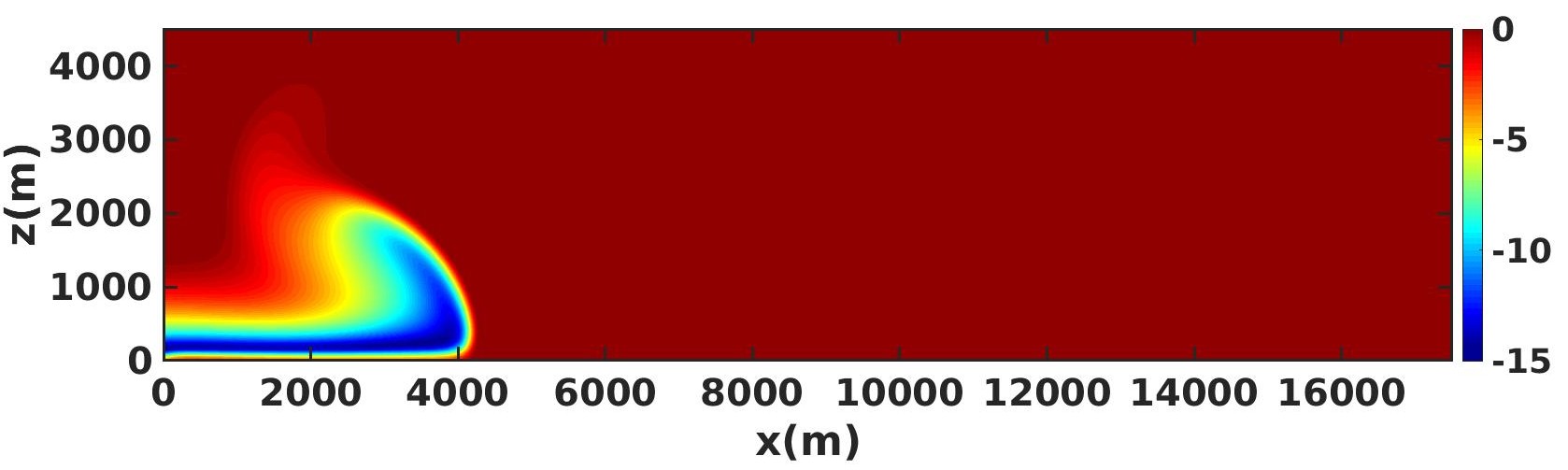}
     \label{fig:Smago_50_300}
    \put(90,50){\textcolor{white}{\footnotesize{$h = 50$ m, $\alpha = 11$, $t = 300$ s}}}
    \end{overpic}
    \\
    \begin{overpic}[width=0.48\textwidth]{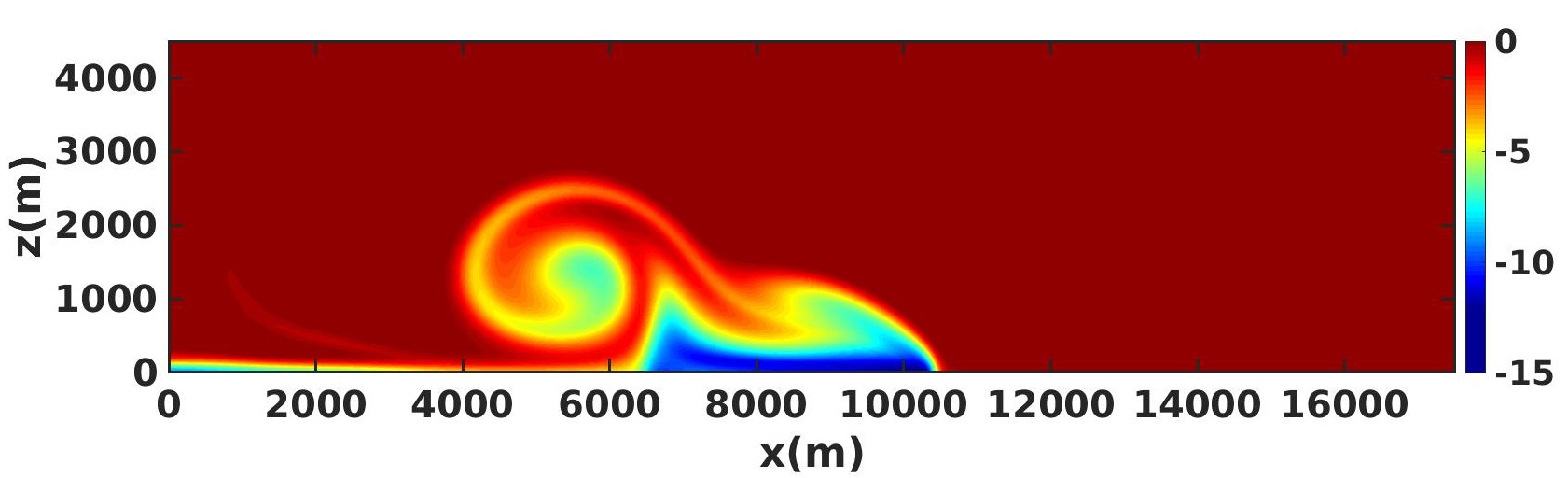}
    \label{fig:Smago_50_600}
    \put(90,50){\textcolor{white}{\footnotesize{$h = 50$ m, $\alpha = 16$, $t = 600$ s}}}
    \end{overpic} 
        \begin{overpic}[width=0.48\textwidth]{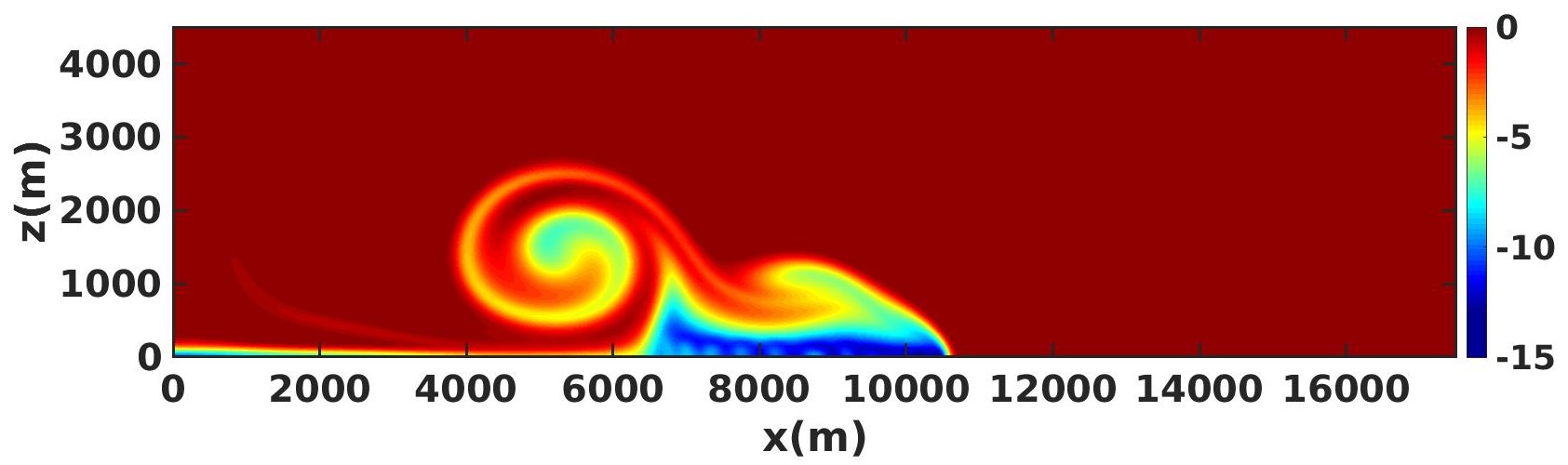}
    \label{fig:Smago_50_600}
    \put(90,50){\textcolor{white}{\footnotesize{$h = 50$ m, $\alpha = 11$, $t = 600$ s}}}
    \end{overpic} 
    \\ 
    \begin{overpic}[width=0.48\textwidth]{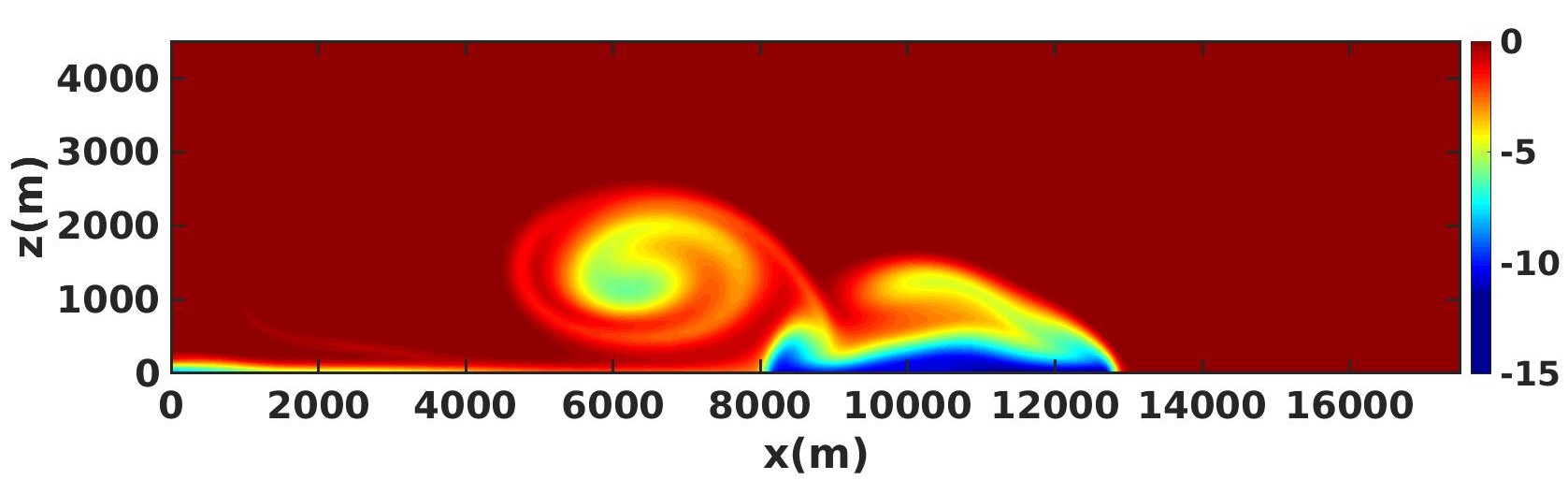}
    \label{fig:Smago_50_750}
    \put(90,50){\textcolor{white}{\footnotesize{$h = 50$ m, $\alpha = 16$, $t = 750$ s}}}
    \end{overpic} 
    \begin{overpic}[width=0.48\textwidth]{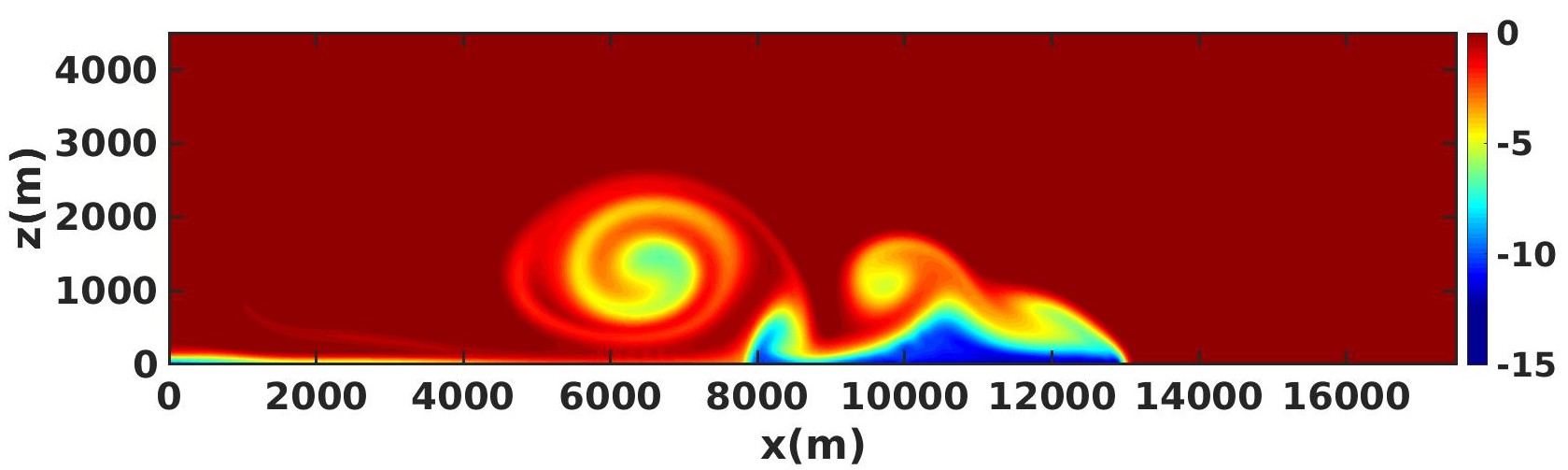}
    \label{fig:Smago_50_750}
    \put(90,50){\textcolor{white}{\footnotesize{$h = 50$ m, $\alpha = 11$, $t = 750$ s}}}
    \end{overpic}
    \\
    \begin{overpic}[width=0.48\textwidth]{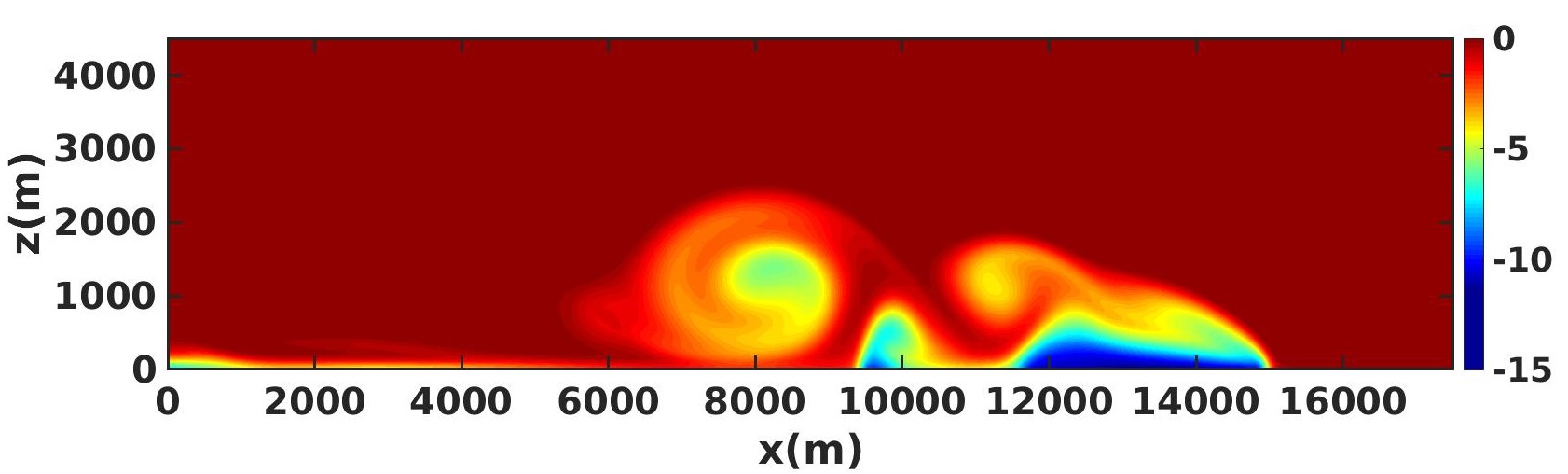}
    \label{fig:Smago_50_900}
    \put(90,50){\textcolor{white}{\footnotesize{$h = 50$ m, $\alpha = 16$, $t = 900$ s}}}
    \end{overpic} 
        \begin{overpic}[width=0.48\textwidth]{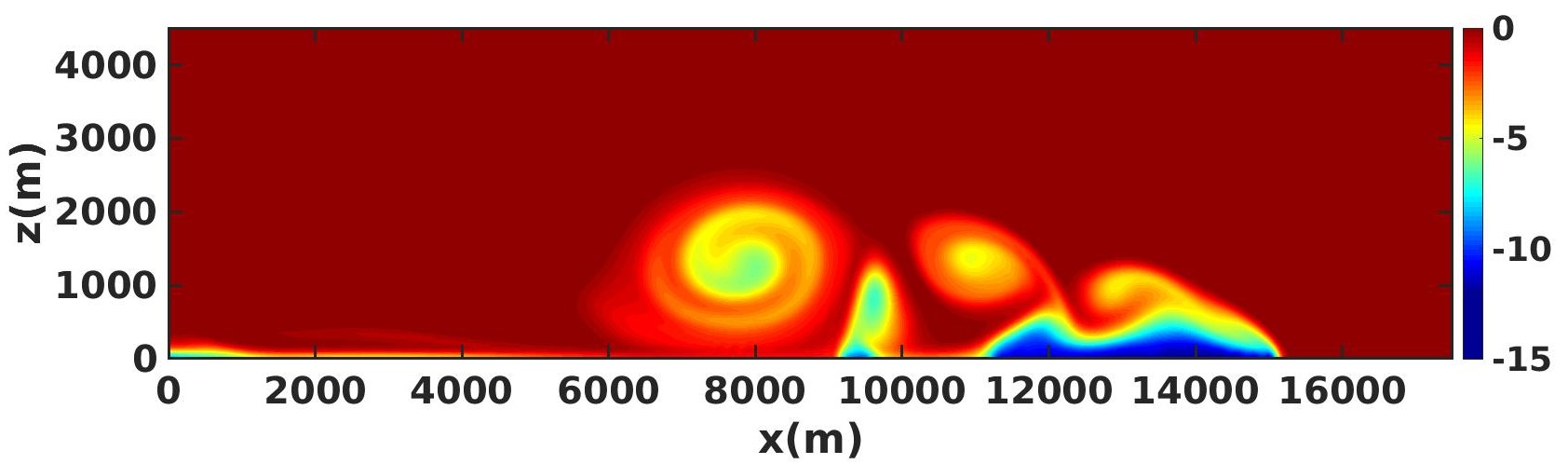}
    \label{fig:Smago_50_900}
    \put(90,50){\textcolor{white}{\footnotesize{$h = 50$ m, $\alpha = 11$, $t = 900$ s}}}
    \end{overpic} 
    \caption{Density current, $a_S$:
time evolution of potential temperature fluctuation $\theta'$ computed with mesh $h = 50$ m and two values of $\alpha$: $\alpha = 16$ (left) and $\alpha = 11$ (right).}
        \label{fig::SmagoRayModel50}
\end{figure}

Next, we focus on the EFR algorithm with indicator function $a_S$. We restrict our attention to meshes $h = 25, 50$ m and, following \cite{marrasNazarovGiraldo2015}, we further refine the mesh to get $h = 12.5$ m. Like in Sec.~\ref{sec:bubble}, we set the value of $\alpha$ using the rule of thumb in
Remark \ref{rem:3} and $C_s = 0.454$ \cite{GQR_OF_clima}. We take $\alpha=4$ m for mesh $h = 12.5$ m because smaller values would lead to instabilities. Then, we use the linear dependence of $\alpha$ on the mesh size to set $\alpha=8$ m for mesh $h = 25$ m and $\alpha=16$ m for mesh $h = 50$ m.
Fig.~\ref{fig::SmagoRayModel125_25} (left), \ref{fig::SmagoRayModel125_25} (right), and \ref{fig::SmagoRayModel50} (left) show the time evolution of the potential temperature fluctuation computed with meshes $h = 12.5, 25, 50$ m, respectively. As expected, more vortical structures appear when we reduce the mesh size. 
The EFR algorithm with $a_S$ produces very similar results to a standard implementation of the Smagorinsky model with mesh $h = 25$ m: compare Fig.~\ref{fig::SmagoRayModel125_25} (right) with Fig.~10 in \cite{GQR_OF_clima}. However, our method does a better job at stabilizing the larger eddies with mesh $h = 12.5$ m: Fig.~\ref{fig::SmagoRayModel125_25} (left) with Fig.~9 in \cite{GQR_OF_clima}. For mesh $h = 50$ m, Fig.~\ref{fig::SmagoRayModel50} (left) indicates that the $\alpha = 16$ leads to overdiffusion. In fact, it provides a smoothed out solution even when compared to the linear filter (see Fig.~\ref{fig:linear_sharing}, second panel from the top). A less dissipative solution 
can be found by lowering the values of $\alpha$. Fig.~\ref{fig::SmagoRayModel50} (right) shows $\theta'$ computed with mesh $h = 50$ m and $\alpha = 11$.
We this new value of $\alpha$ the solution obtained with mesh $h = 50$ m looks similar to the solution given by 
mesh $h = 25$ m. This is confirmed by Table \ref{tab:4}, 
which reports the front locations at $t = 900$ s
obtained with EFR and $a_S$. The location computed with 
mesh $h = 25$ m and $\alpha = 8$ is very close to the location given by mesh $h = 50$ m and $\alpha = 11$. This is exactly what expected from the EFR algorithm: with a proper tuning of $\alpha$, one can use coarser meshes without compromising accuracy. 

\begin{table}[htb!]
\begin{center}
\begin{tabular}{ | c | c | c |  c | }
\hline
Method & $h$ (m) & $\alpha$ (m) & Front Location (m)  \\
 \hline 
 Ref.~\cite{marrasNazarovGiraldo2015} & 12.5 &-& 15056 \\
 \hline
 EFR, $a_S$ & 12.5 & 4 &  15550 \\
 \hline
 EFR, $a_S$ & 25 & 8 & 15300   \\
 \hline
  EFR, $a_S$ & 50 & 11 & 15220 \\
 \hline
 EFR, $a_S$ & 50 & 16 & 15090 \\
 \hline
    Ref.~\cite{strakaWilhelmson1993} & (25, 200) & -& (14533,17070)
 \\
  \hline
\end{tabular}
\caption{Density current, $a_S$: front location at $t = 900$ s obtained with the EFR algorithm and different meshes. Our results are compared against results from 
\cite{marrasNazarovGiraldo2015,strakaWilhelmson1993}. For reference \cite{marrasNazarovGiraldo2015}, we report only the front location computed with the finest resolution.
For reference \cite{strakaWilhelmson1993}, we provide the range of mesh sizes and front location values obtained with different methods.
}\label{tab:4}
\end{center}
\end{table}

\begin{figure}[htb!]
    \centering
    \begin{overpic}[width=0.5\textwidth]{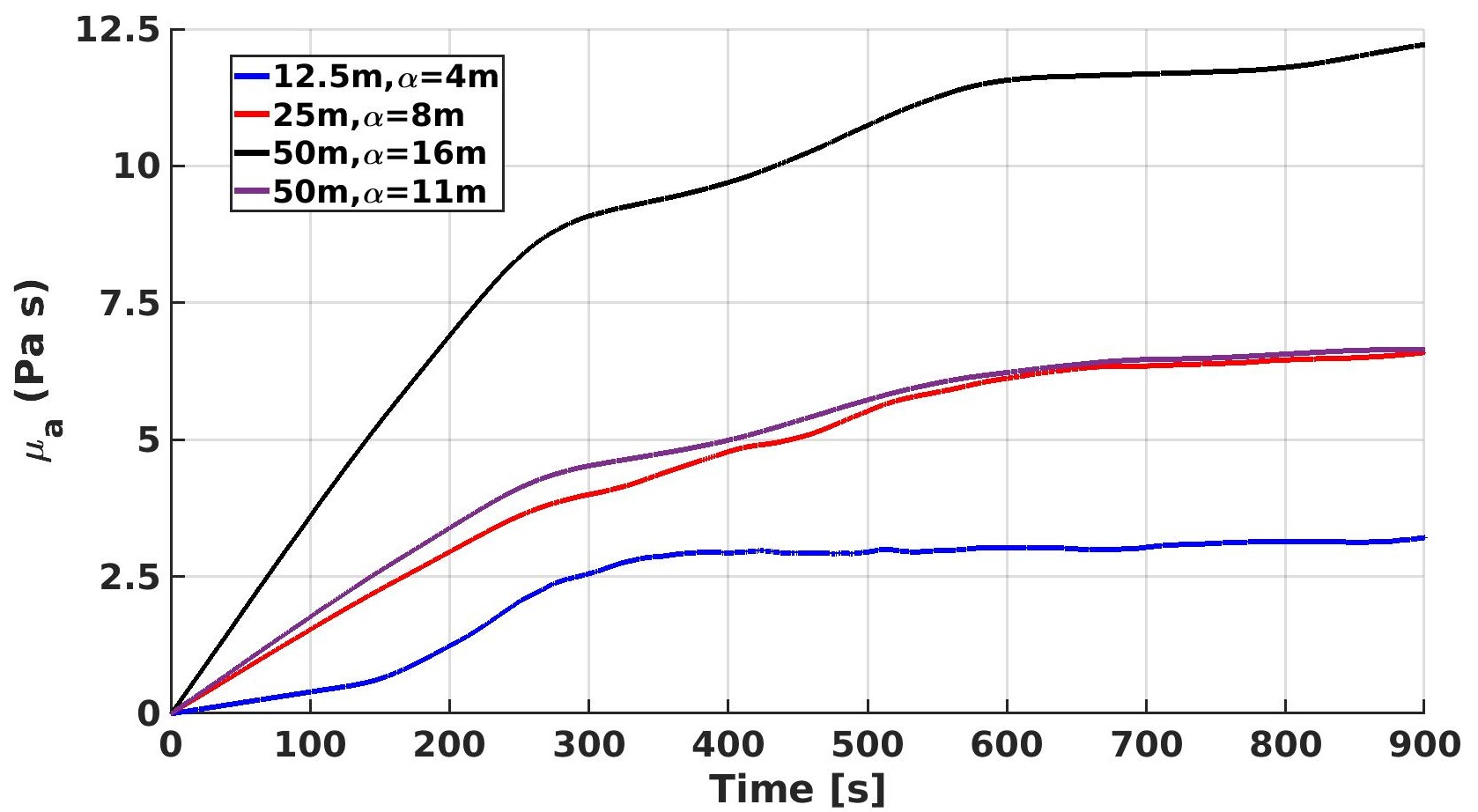}
    \put(120,132){{$a_S$ }}
    \end{overpic}
    \begin{overpic}[width=0.49\textwidth]{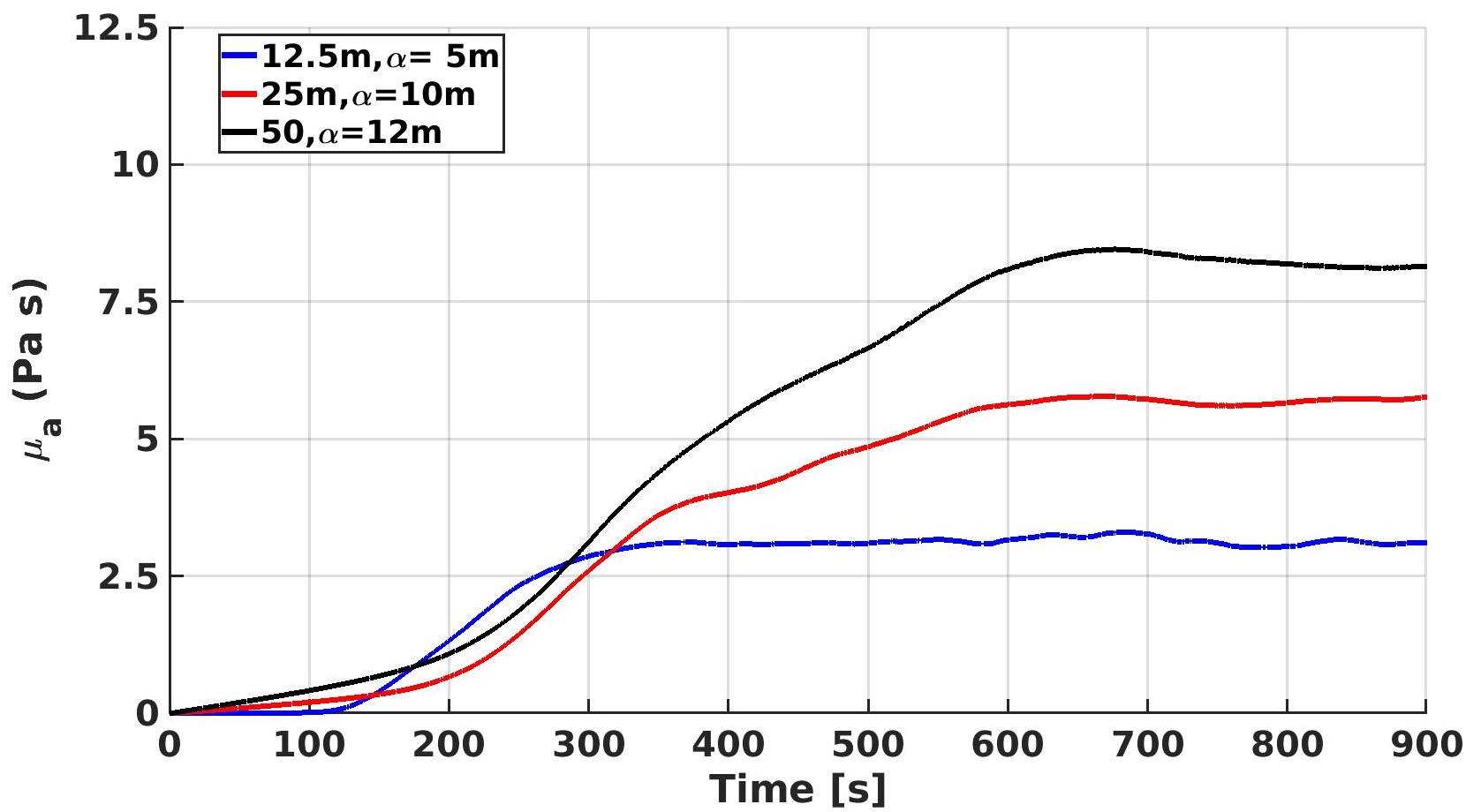} 
    \put(120,132){{$a_D$}}
    \end{overpic}
    \caption{Density current: time evolution of the average eddy viscosity \eqref{eq:mu_av} for $a_S$ (left) and $a_D$ (right) with meshes $h = 12.5, 25, 50$ m.
    }
    \label{fig:CompViscosities}
\end{figure}

Now, let us turn our attention to $a_D$. In Sec.~\ref{sec:bubble}, we have shown that $a_D$ is a more selective indicator function than $a_S$. Thus, we slightly increase the values of $\alpha$ used for $a_S$ since at the moment we do not have a better criterion to set $\alpha$ for $a_D$. We take $\alpha = 5$ for mesh $h = 12.5$ m,  $\alpha = 10$ for mesh $h = 25$ m, $\alpha = 12$ for mesh $h = 50$ m. Before showing the solutions obtained with 
$a_D$ and these values of $\alpha$, in Fig.~\ref{fig:CompViscosities} we compare the time evolution of the space-averaged artificial viscosity
\begin{equation}\label{eq:mu_av}
\mu_{av} = \dfrac{1}{\Omega} \int_{\Omega} \overline{\mu}_h d \Omega
\end{equation}
obtained with $a_S$ and $a_D$ for meshes $h = 12.5, 25, 50$ m.
We recall that $\overline{\mu}_h$ is defined in \eqref{eq:mubar_h}. 
From Fig.~\ref{fig:CompViscosities}, we see that $a_S$ with $\alpha = 4$ and $a_D$ with $\alpha = 5$ introduce roughly the same amount of space-averaged artificial viscosity for most of the time interval under consideration in the case of mesh $h = 12.5$ m. The main difference for this mesh is that $a_D$ introduces almost no artificial diffusion till about 100 s and then ramps it up faster than $a_S$. A rather fast ramp is observed for $a_D$ also in the case of meshes $h = 25, 50$ m: $\mu_{av}$ remains small while the cold perturbation falls due to negative buoyancy and it increases as the cold front propagates horizontally. For mesh $h = 25$ m, $\mu_{av}$ given by $a_S$ with $\alpha = 8$ grows almost linearly till about 300 s and then around 600 s it flattens. As already evident from Fig.~\ref{fig::SmagoRayModel50}, $a_S$ with $\alpha = 16$ introduces too much artificial viscosity for mesh $h = 50$ m. This can be fixed by decreasing the value to $\alpha = 11$, which introduces a similar amount of $\mu_{av}$ as $\alpha = 8$ does for mesh 
$h = 25$ m.

\begin{figure}[htb!]
    \centering
    \begin{overpic}[width=0.48\textwidth]{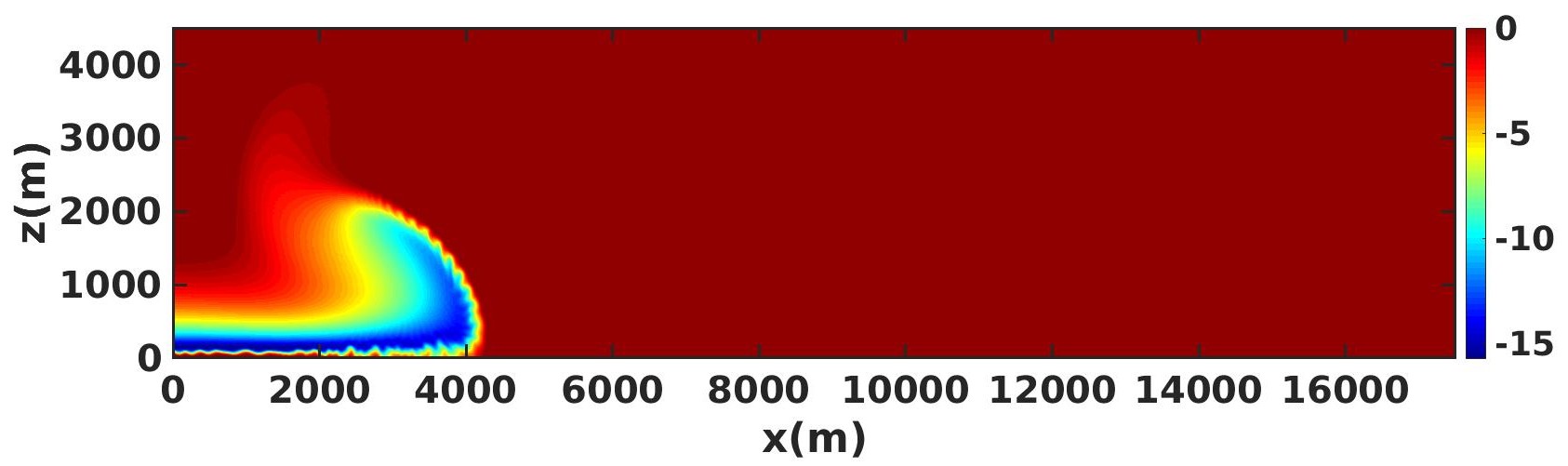}
    \label{fig:decon_125_300}
    \put(90,50){\textcolor{white}{\footnotesize{$h = 12.5$ m, $\alpha = 5$, $t = 300$ s}}}
    \end{overpic} 
    \begin{overpic}[width=0.48\textwidth]{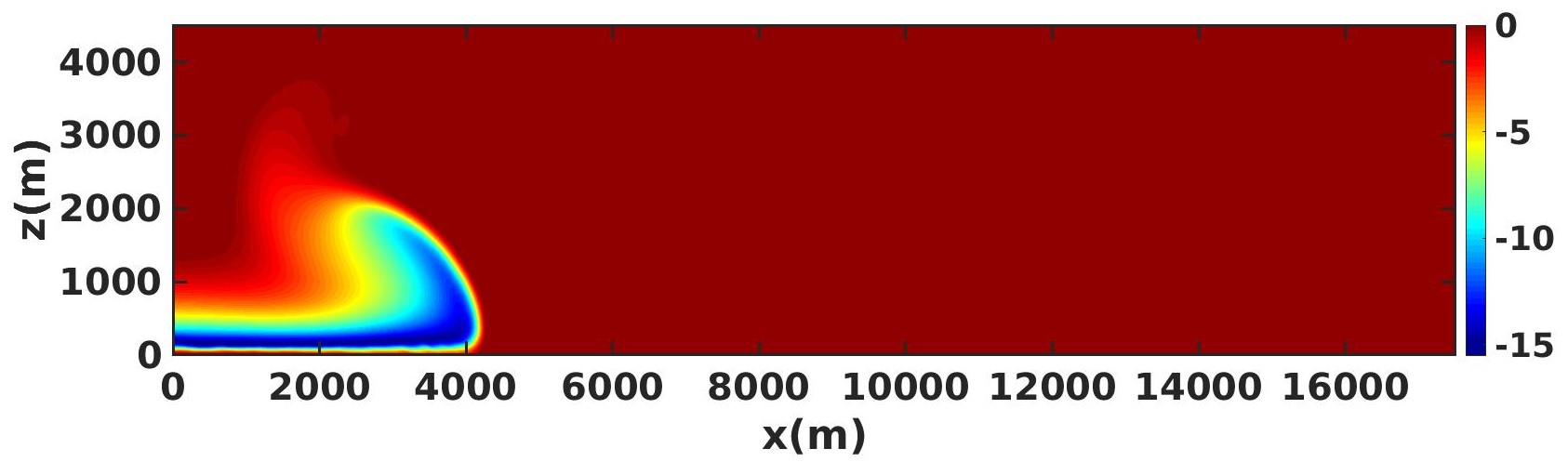}
    \label{fig:decon_125_300}
    \put(90,50){\textcolor{white}{\footnotesize{$h = 25$ m, $\alpha = 10$, $t = 300$ s}}}
    \end{overpic}
    \\
    \begin{overpic}[width=0.48\textwidth]{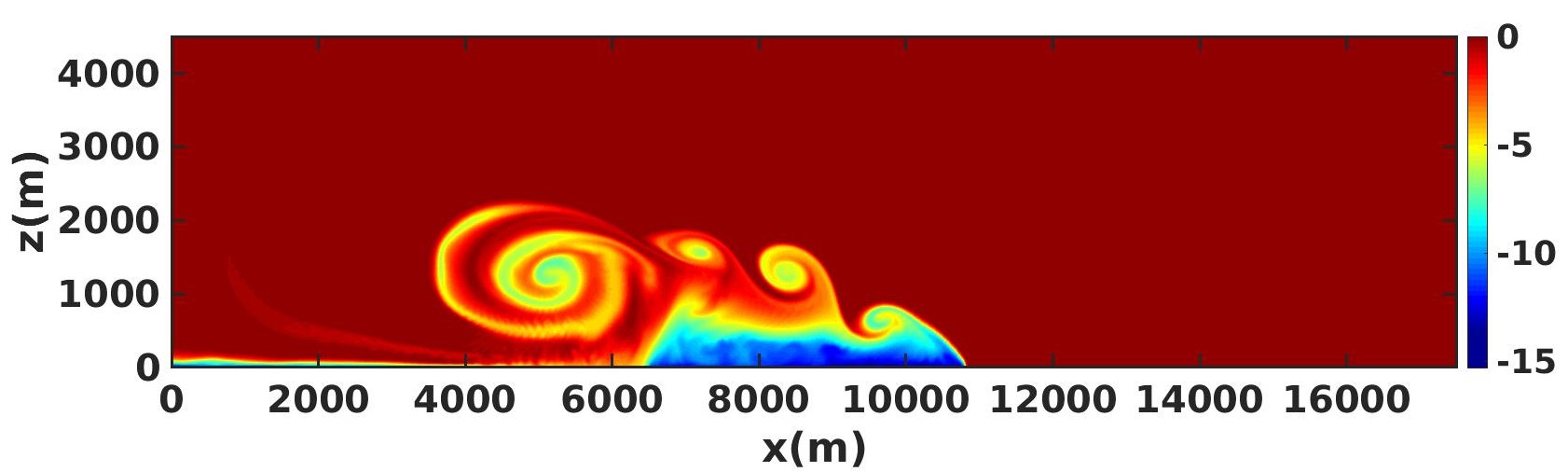}
    \label{fig:decon_125_600}
    \put(90,50){\textcolor{white}{\footnotesize{$h = 12.5$ m, $\alpha = 5$, $t = 600$ s}}}
    \end{overpic} 
    \begin{overpic}[width=0.48\textwidth]{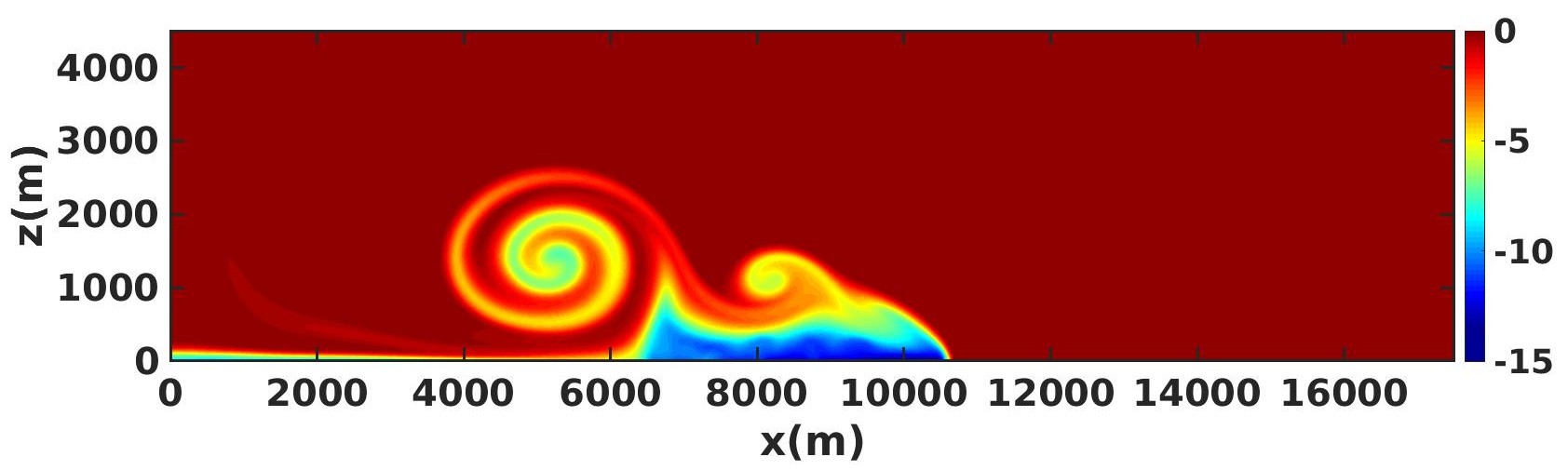}
    \label{fig:decon_125_600}
    \put(90,50){\textcolor{white}{\footnotesize{$h = 25$ m, $\alpha = 10$, $t = 600$ s}}}
    \end{overpic} 
    \\ 
    \begin{overpic}[width=0.48\textwidth]{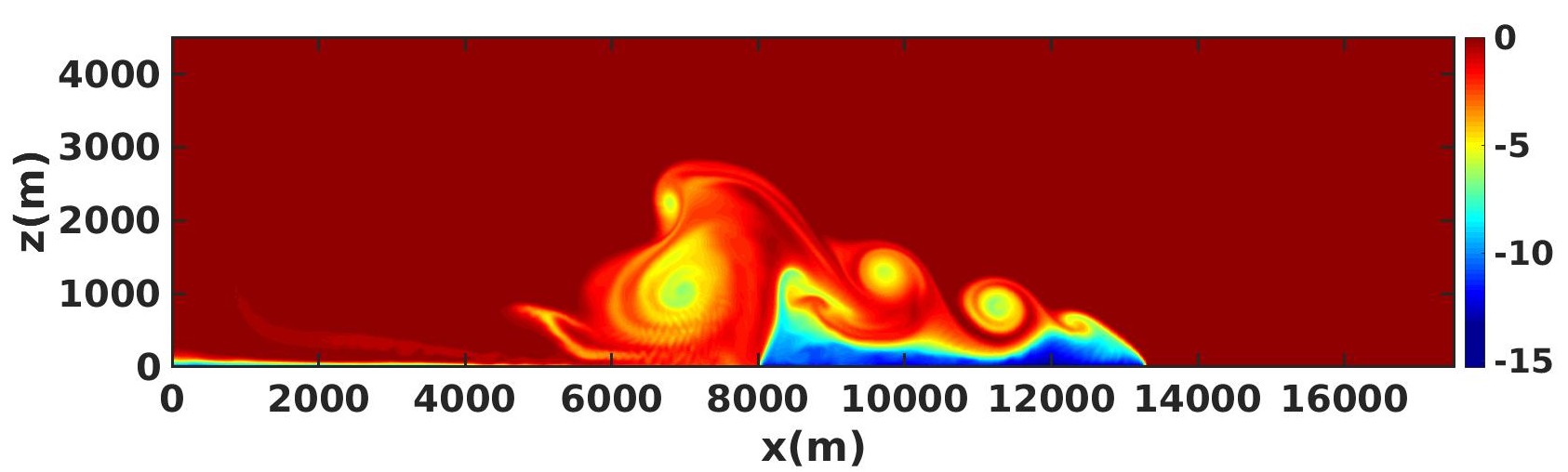}
    \label{fig:decon_125_750}
    \put(90,50){\textcolor{white}{\footnotesize{$h = 12.5$ m, $\alpha = 5$, $t = 750$ s}}}
    \end{overpic}
    \begin{overpic}[width=0.48\textwidth]{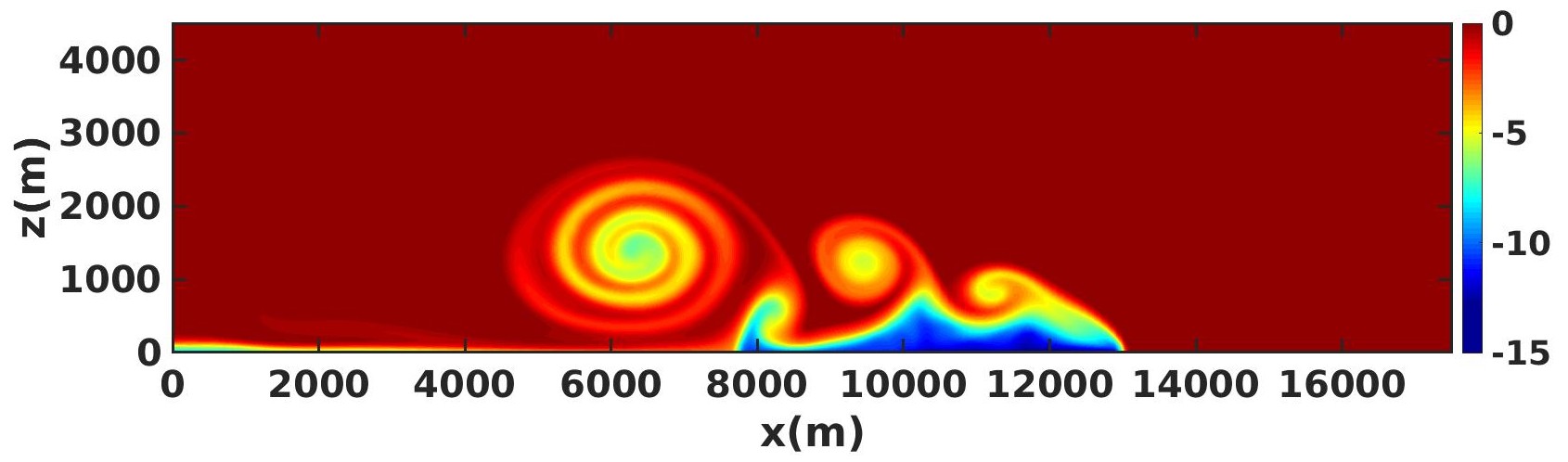}
    \label{fig:decon_125_750}
    \put(90,50){\textcolor{white}{\footnotesize{$h = 25$ m, $\alpha = 10$, $t = 750$ s}}}
    \end{overpic}
    \\
    \begin{overpic}[width=0.48\textwidth]{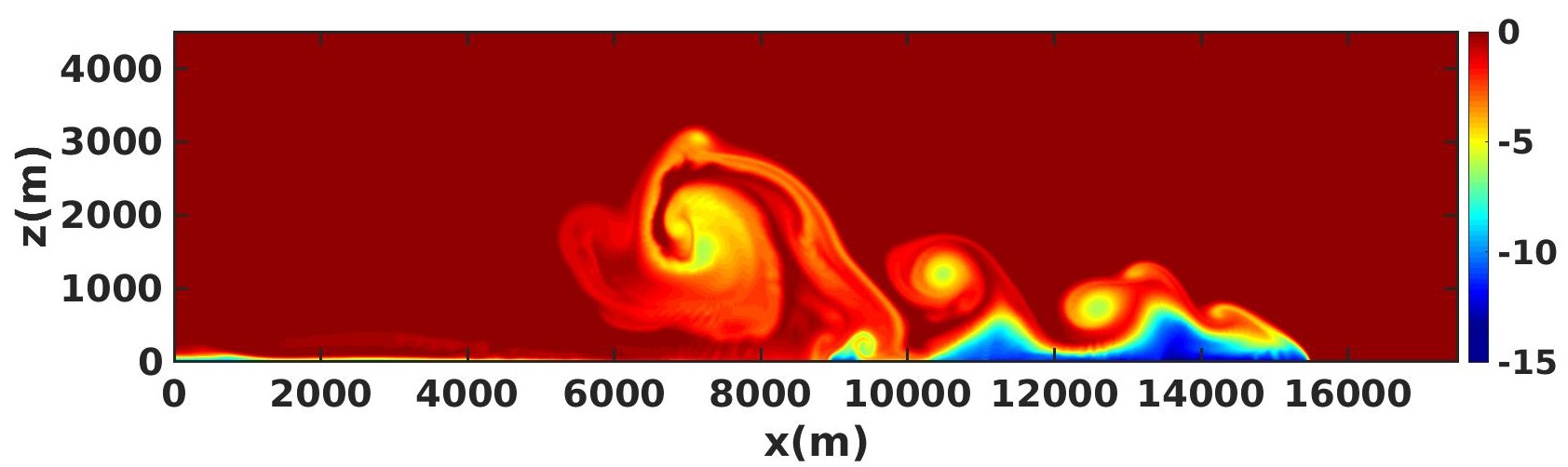}
    \label{fig:decon_125_900}
    \put(90,50){\textcolor{white}{\footnotesize{$h = 12.5$ m, $\alpha = 5$, $t = 900$ s}}}
    \end{overpic} 
    \begin{overpic}[width=0.48\textwidth]{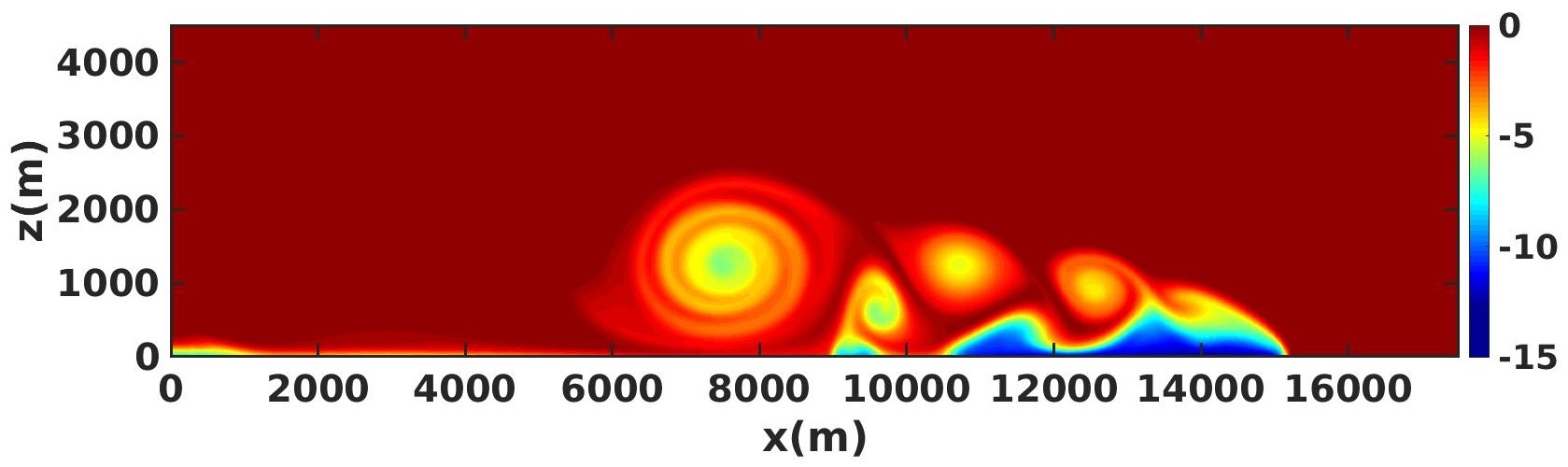}
    \label{fig:decon_125_900}
    \put(90,50){\textcolor{white}{\footnotesize{$h = 25$ m, $\alpha = 10$, $t = 900$ s}}}
    \end{overpic} 
    \caption{
    Density current, $a_D$:
time evolution of potential temperature fluctuation $\theta'$ computed with mesh $h = 12.5$ m and $\alpha = 5$ (left) and with mesh $h = 25$ m and $\alpha = 10$ (right).}
        \label{fig:deconRayModel125_25}
\end{figure}

\begin{figure}[htb!]
    \centering
    \begin{overpic}[width=0.48\textwidth]{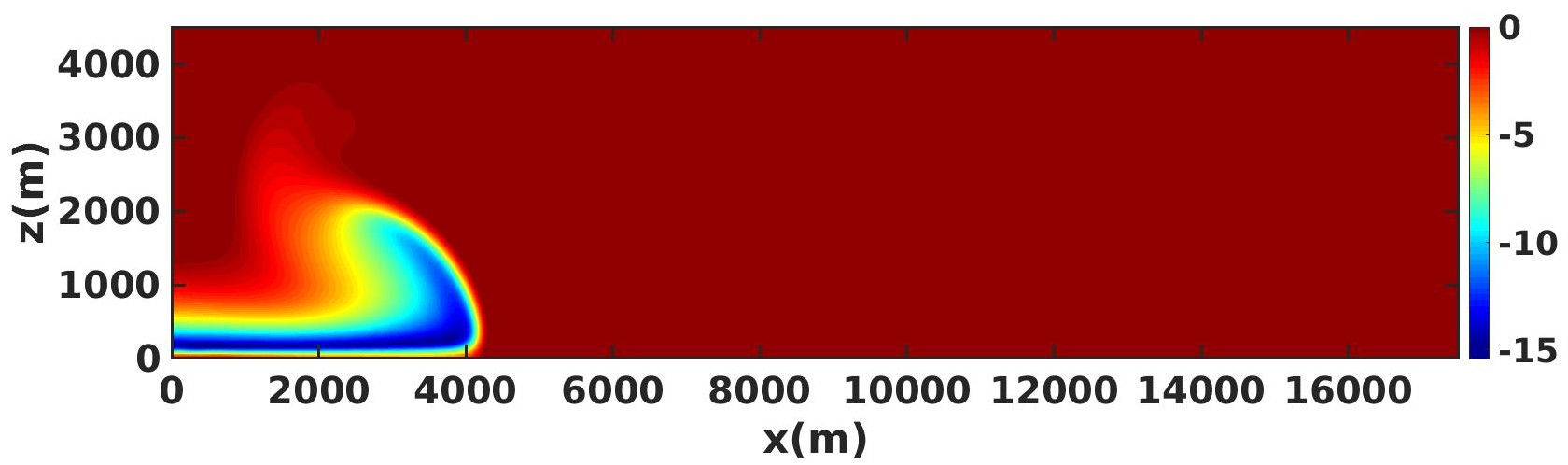}
    \label{fig:decon_125_300}
    \put(150,50){\textcolor{white}{$t = 300$ s}}
    \end{overpic}\\
    \begin{overpic}[width=0.48\textwidth]{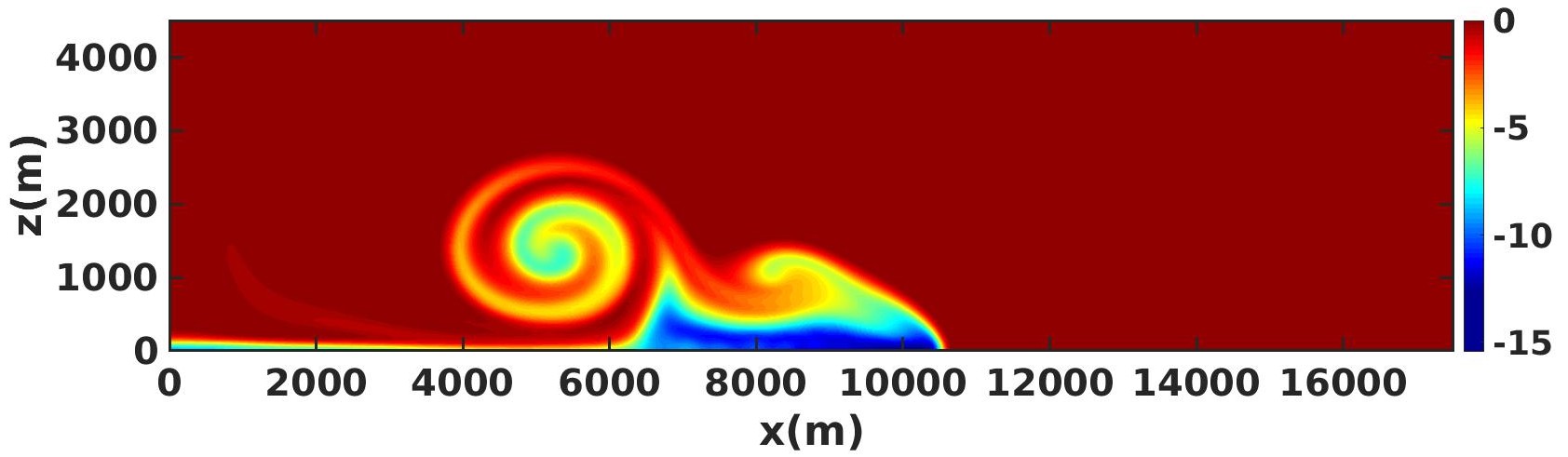}
    \label{fig:decon_125_600}
    \put(150,50){\textcolor{white}{$t = 600$ s}}
    \end{overpic} \\ 
    \begin{overpic}[width=0.48\textwidth]{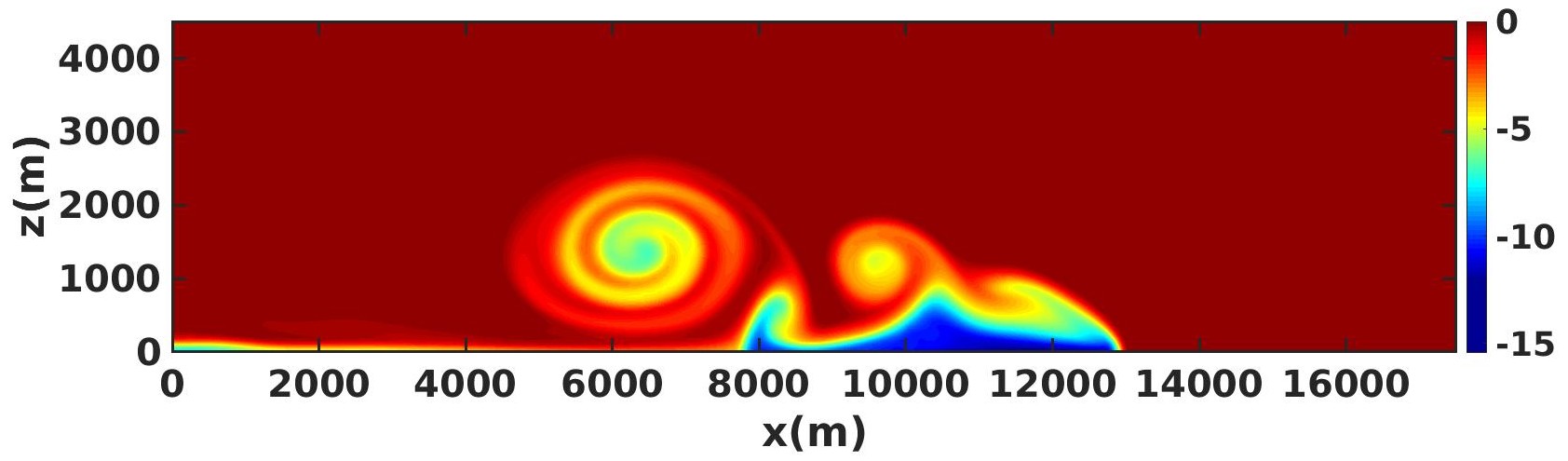}
    \label{fig:decon_125_750}
    \put(150,50){\textcolor{white}{$t = 750$ s}}
    \end{overpic}\\
    \begin{overpic}[width=0.48\textwidth]{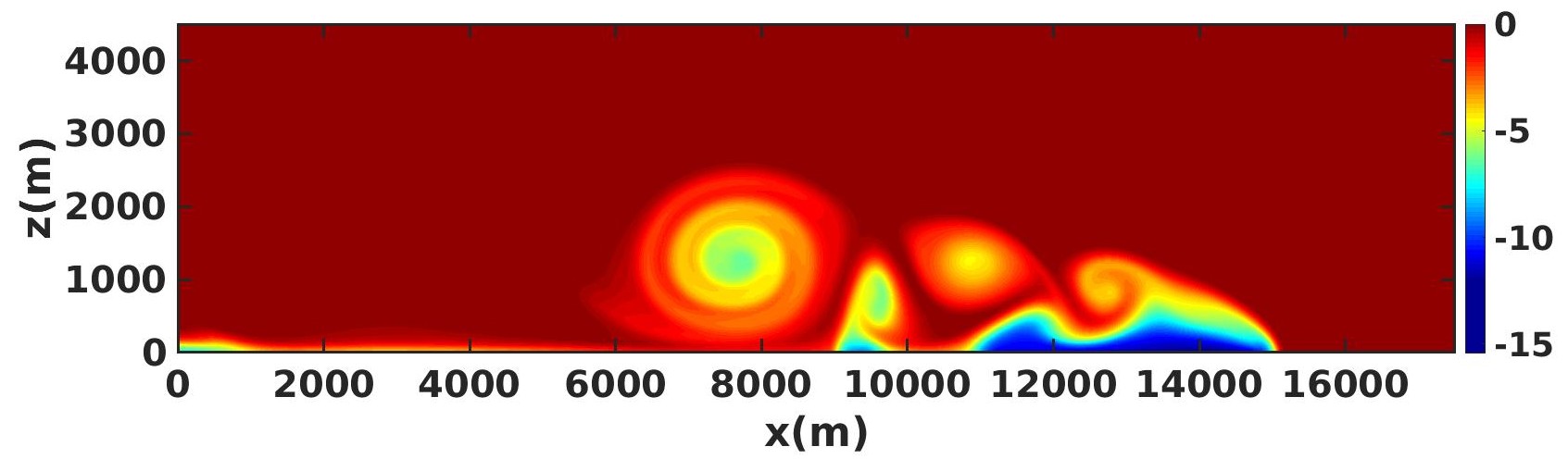}
    \label{fig:decon_125_900}
    \put(150,50){\textcolor{white}{$t = 900$ s}}
    \end{overpic} 
    \caption{
    Density current, $a_D$, $\alpha = 12$:
time evolution of potential temperature fluctuation $\theta'$ computed with mesh $h = 50$ m.}
        \label{fig:deconRayModel50_12}
\end{figure}

Fig.~\ref{fig:deconRayModel125_25} (left), \ref{fig:deconRayModel125_25}  (right), and \ref{fig:deconRayModel50_12} display the time evolution of the potential temperature fluctuation computed with $a_D$ and the chosen values of $\alpha$ for meshes $h = 12.5, 25, 50$ m, respectively. When compared to the respective counterparts obtained with
$a_S$ (namely Fig.~\ref{fig::SmagoRayModel125_25} (left), \ref{fig::SmagoRayModel125_25} (right), and \ref{fig::SmagoRayModel50} (right)), all the observations made about Fig.~\ref{fig:CompViscosities} are confirmed: the solutions obtained with the finer mesh are initially comparable early and then some differences are observed for the larger recirculations, the solutions for the intermediate mesh are remarkably similar. While Fig.~\ref{fig:CompViscosities} suggests that EFR with $a_D$ and $\alpha = 12$ is more diffusive (in average) than with $a_S$ and $\alpha = 11$, the respective solutions in Fig.~\ref{fig:deconRayModel50_12} and Fig.~\ref{fig::SmagoRayModel50} (right) are surprisingly similar.

Table \ref{tab:5} reports the front locations at $t = 900$ s
obtained with EFR and $a_D$ for the three meshes under consideration. The three locations are within about 400 m of each other, with the front becoming faster as the mesh is refined. This was the case also for $a_S$ (see Table \ref{tab:4}). The opposite trend is observed for $a_L$ (see Table \ref{tab:3}), i.e., the front slows down as the mesh is refined, although the locations are only roughly 40 m apart. In any case, our results fall well within the results from
\cite{strakaWilhelmson1993} and are close to the results from \cite{marrasNazarovGiraldo2015}.

\begin{table}[htb!]
\begin{center}
\begin{tabular}{ | c | c | c |  c | }
\hline
Method & $h$ (m) & $\alpha$ (m) & Front Location (m)  \\

 \hline
  Ref.~\cite{marrasNazarovGiraldo2015} & 12.5 &-& 15056 \\
 \hline 
 EFR, $a_D$ & 12.5 & 5 & 15560  \\
 \hline
 EFR, $a_D$ & 25 & 10 & 15215  \\
 \hline
 EFR, $a_D$ & 50 & 12 & 15120 \\
    \hline
    Ref.~\cite{strakaWilhelmson1993} & (25, 200) & -& (14533,17070)
 \\
  \hline
\end{tabular}
\caption{Density current, $a_D$: front location at $t = 900$ s obtained with the EFR algorithm and different meshes. Our results are compared against results from 
\cite{strakaWilhelmson1993,marrasNazarovGiraldo2015}.
For reference \cite{strakaWilhelmson1993}, we provide the range of mesh sizes and front location values obtained with different methods.
For reference \cite{marrasNazarovGiraldo2015}, we report only the front location computed with the finest resolution.}\label{tab:5}
\end{center}
\end{table}

\begin{figure}[htb!]
    \centering
    \begin{overpic}[width=0.48\textwidth,grid=false]{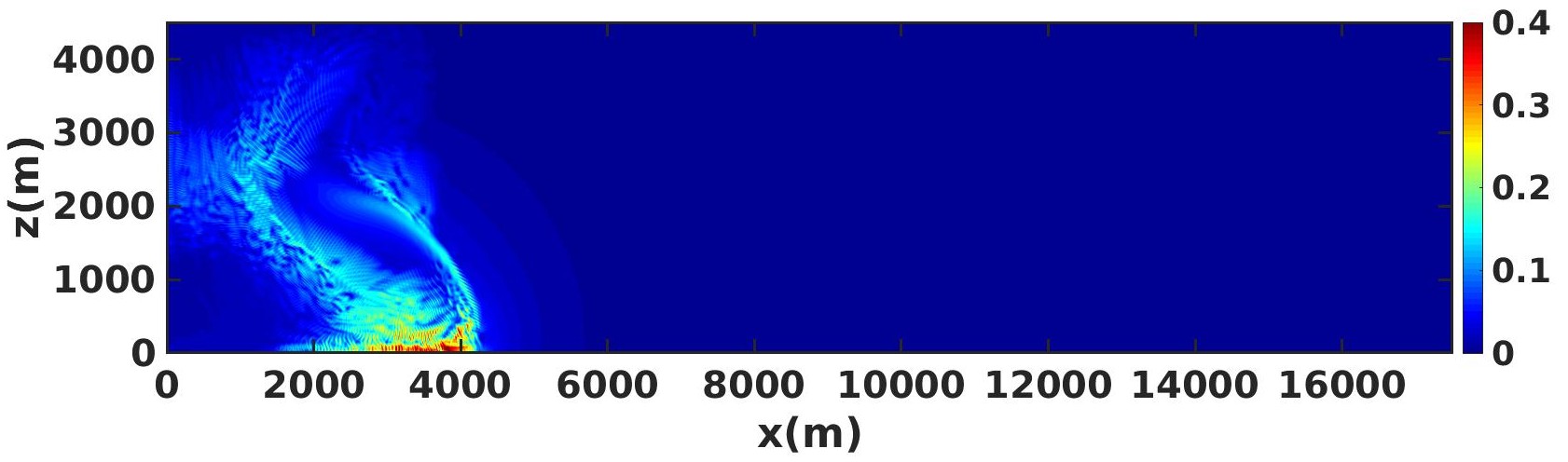}
    \put(125,50){\textcolor{white}{\footnotesize{$a_S$, $\alpha = 4$, $t = 300$ s}}}
    \end{overpic}
    \begin{overpic}[width=0.48\textwidth]{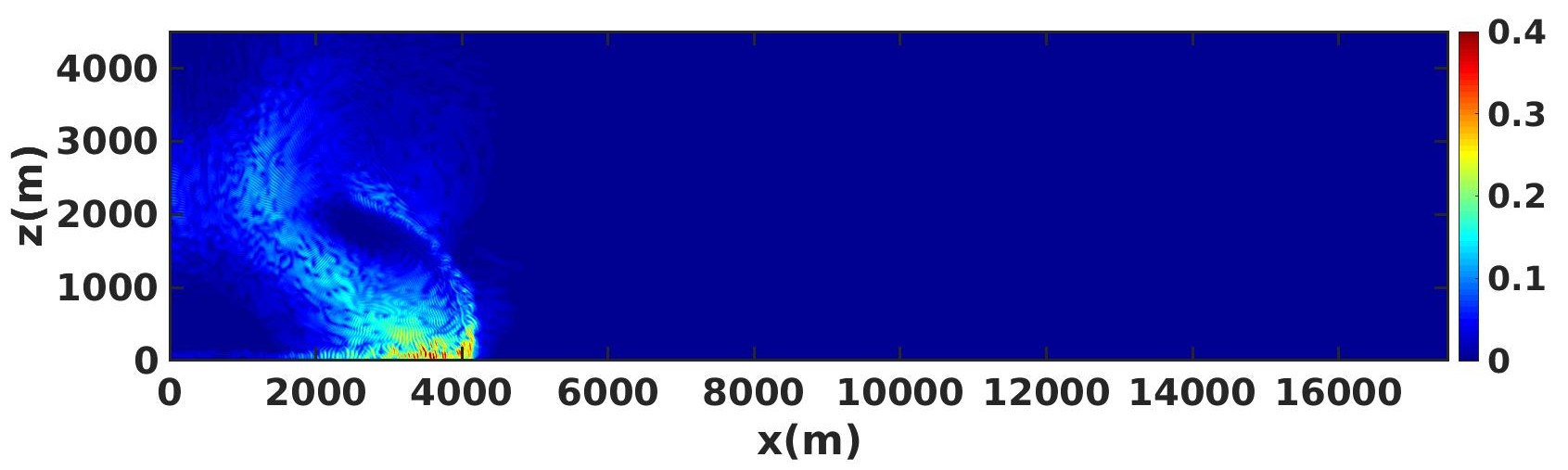}
    \put(125,50){\textcolor{white}{\footnotesize{$a_D$, $\alpha = 5$, $t = 300$ s}}}
    \end{overpic} \\ \vspace{0.2cm}
    \begin{overpic}[width=0.48\textwidth]{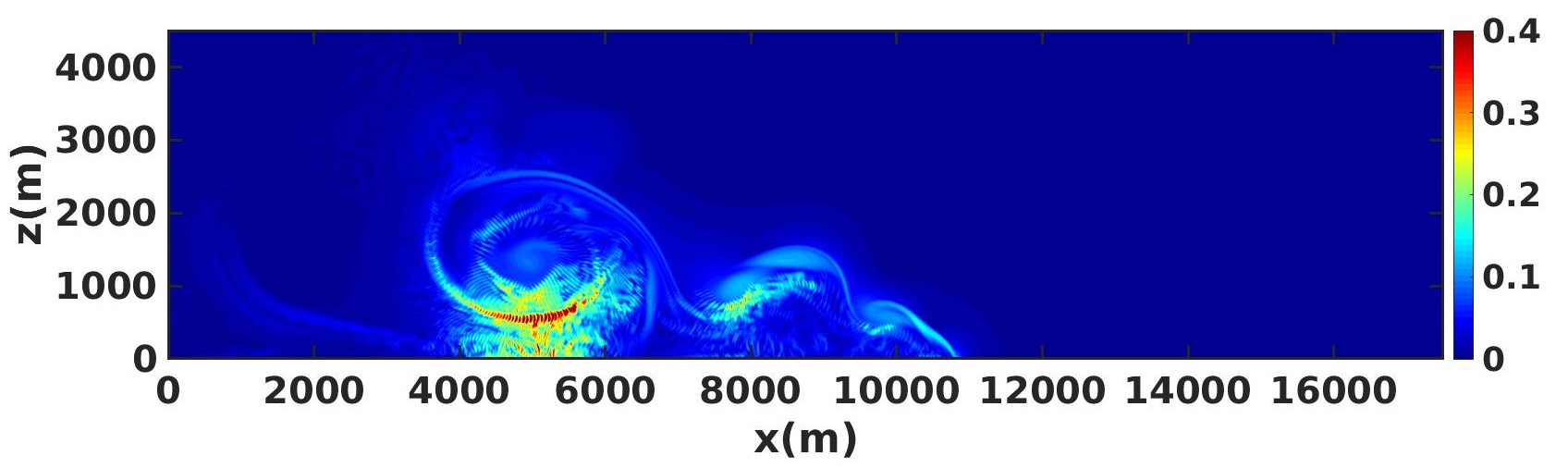}
    \put(125,50){\textcolor{white}{\footnotesize{$a_S$, $\alpha = 4$, $t = 600$ s}}}
    \end{overpic}
    \begin{overpic}[width=0.48\textwidth]{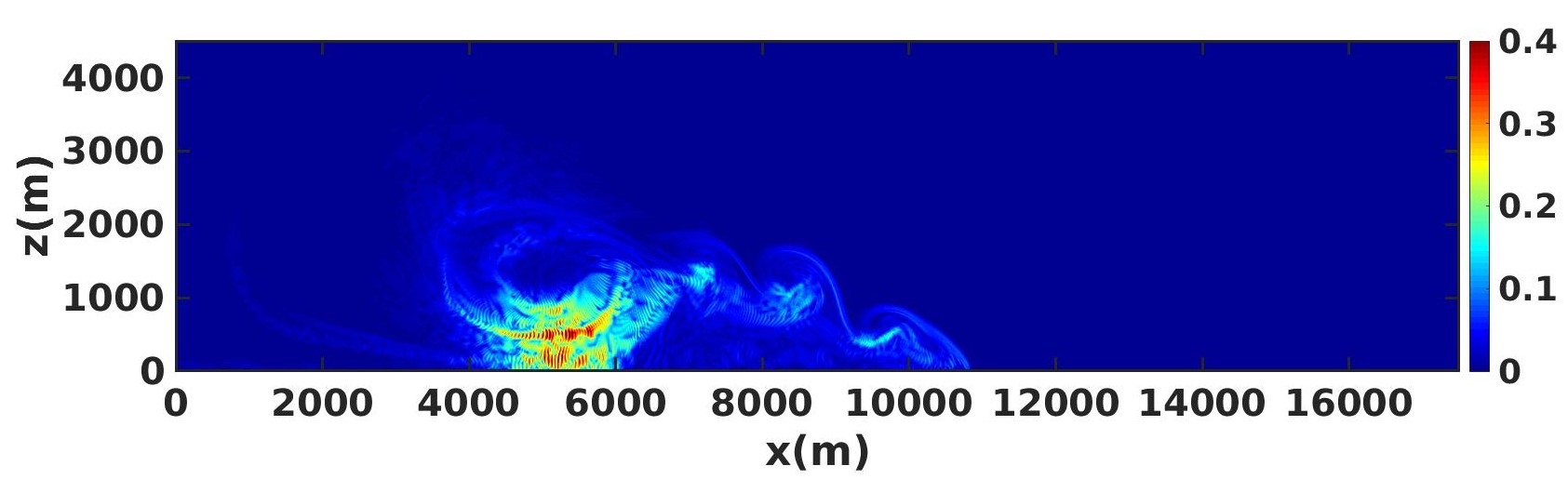}
    \put(125,50){\textcolor{white}{\footnotesize{$a_D$, $\alpha = 5$, $t = 600$ s}}}
    \end{overpic} \\ \vspace{0.2cm}
    \begin{overpic}[width=0.48\textwidth]{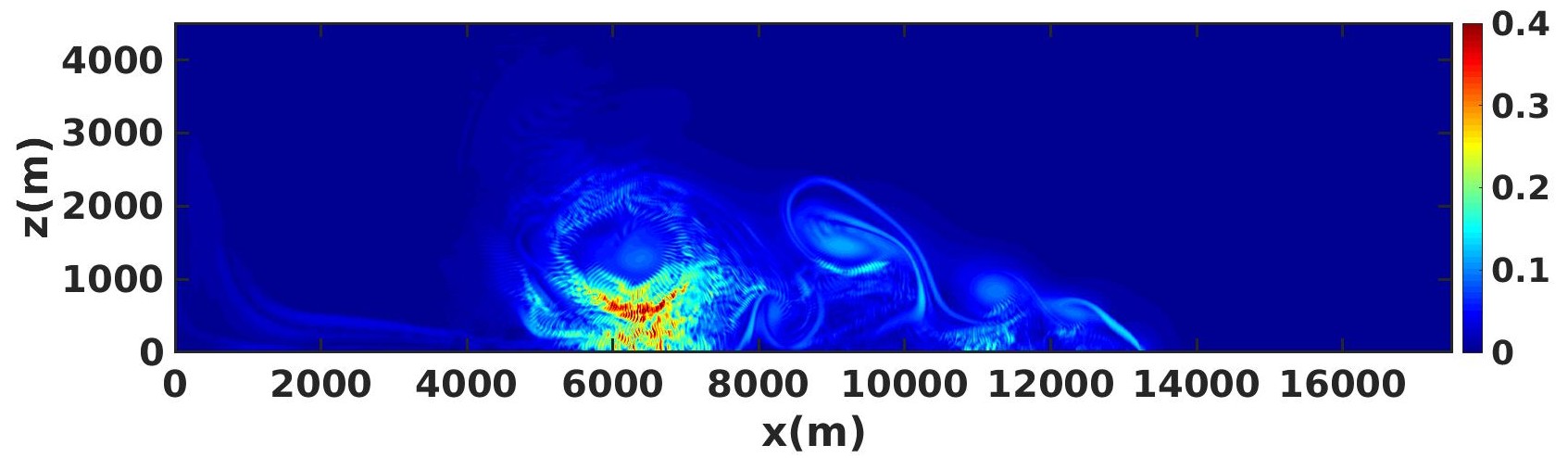}
    \put(125,50){\textcolor{white}{\footnotesize{$a_S$, $\alpha = 4$, $t = 750$ s}}}
    \end{overpic}
    \begin{overpic}[width=0.48\textwidth]{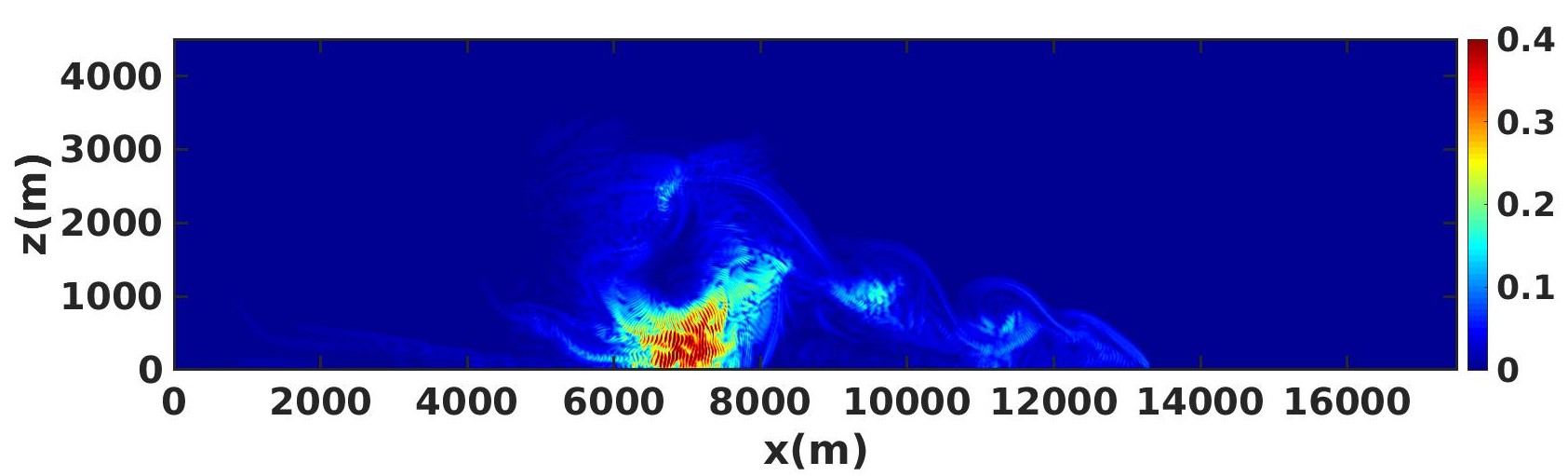}
    \put(125,50){\textcolor{white}{\footnotesize{$a_D$, $\alpha = 5$, $t = 750$ s}}}
    \end{overpic} \\ \vspace{0.2cm}
    \begin{overpic}[width=0.48\textwidth]{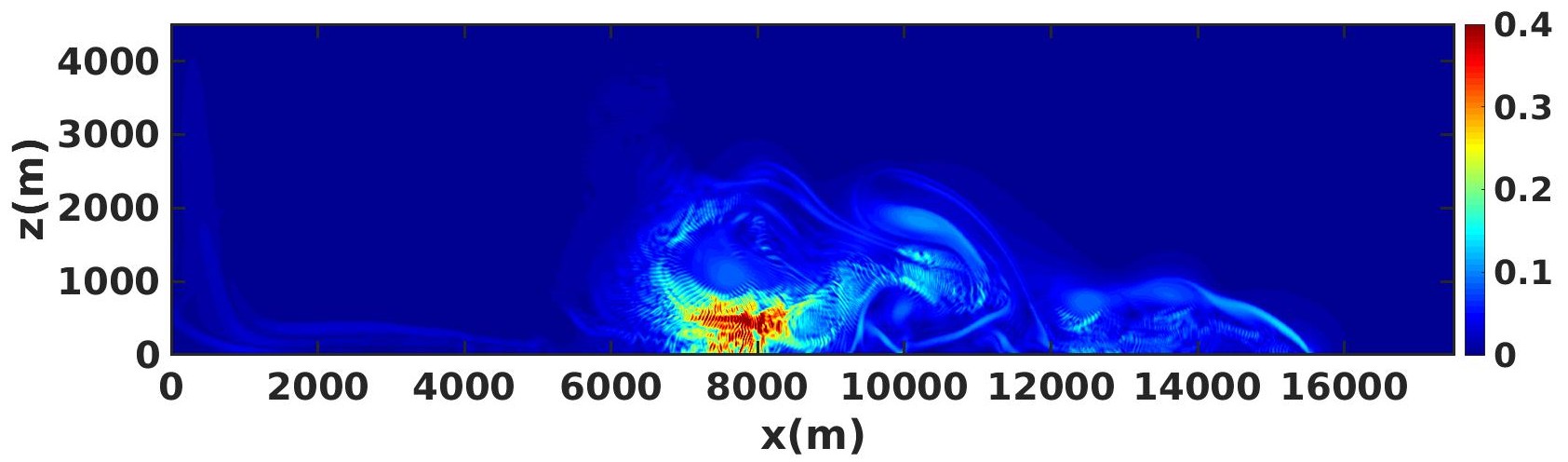}
    \put(125,50){\textcolor{white}{\footnotesize{$a_S$, $\alpha = 4$, $t = 900$ s}}}
    \end{overpic}
    \begin{overpic}[width=0.48\textwidth]{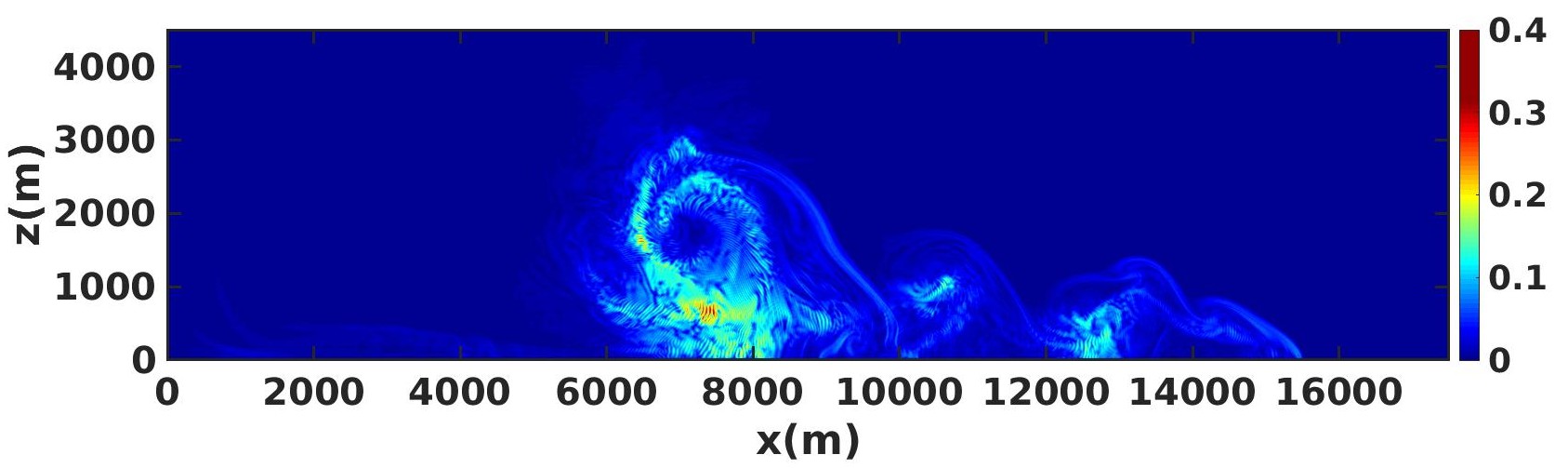}
    \put(125,50){\textcolor{white}{\footnotesize{$a_D$, $\alpha = 5$, $t = 900$ s}}}
    \end{overpic} 
    \caption{Density current, mesh $h = 12.5$ m:
time evolution of $a_S$ with $\alpha=4$ m (left) and $a_D$ with $\alpha=5$ m (right).
    }\label{fig::dc_12.5_aU_1}
\end{figure}

\begin{figure}[htb!]
    \centering
    \begin{overpic}[width=0.48\textwidth]{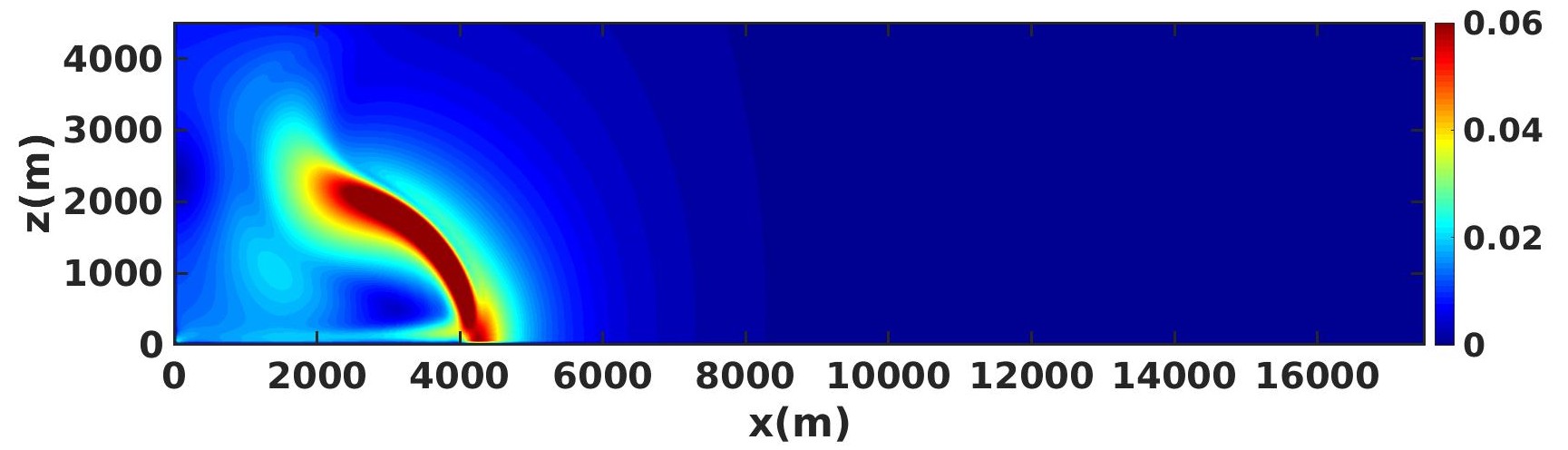}
    \put(115,50){\textcolor{white}{\footnotesize{$a_S$, $\alpha = 11$, $t = 300$ s}}}
    \end{overpic}
    \begin{overpic}[width=0.48\textwidth]{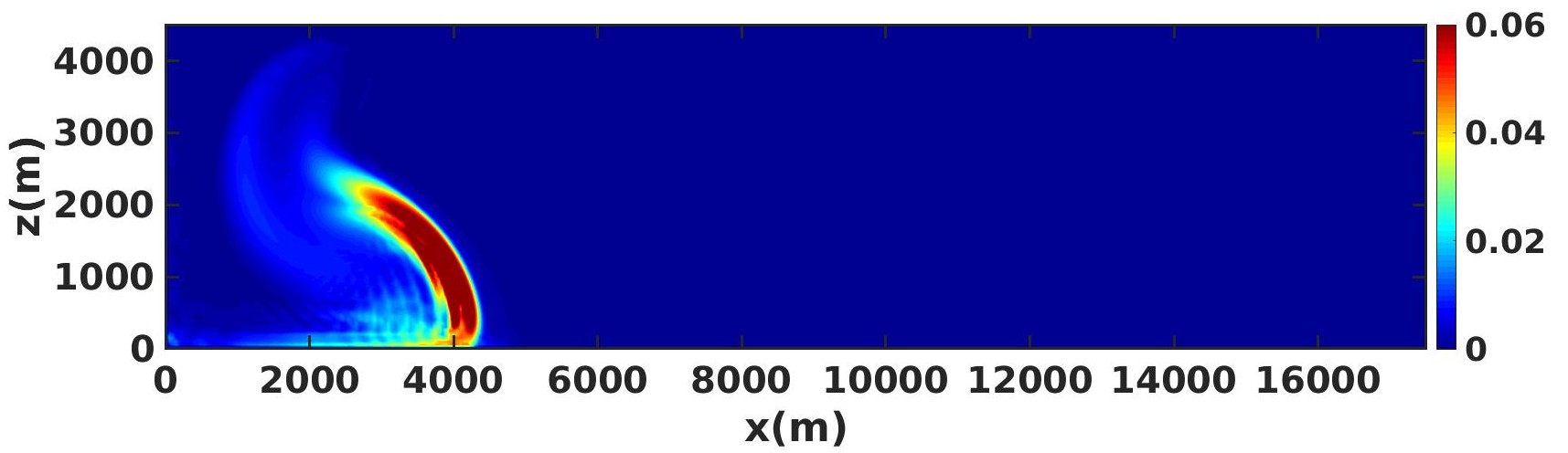}
    \put(115,50){\textcolor{white}{\footnotesize{$a_D$, $\alpha = 12$, $t = 300$ s}}}
    \end{overpic} \\ \vspace{0.2cm}
    \begin{overpic}[width=0.48\textwidth]{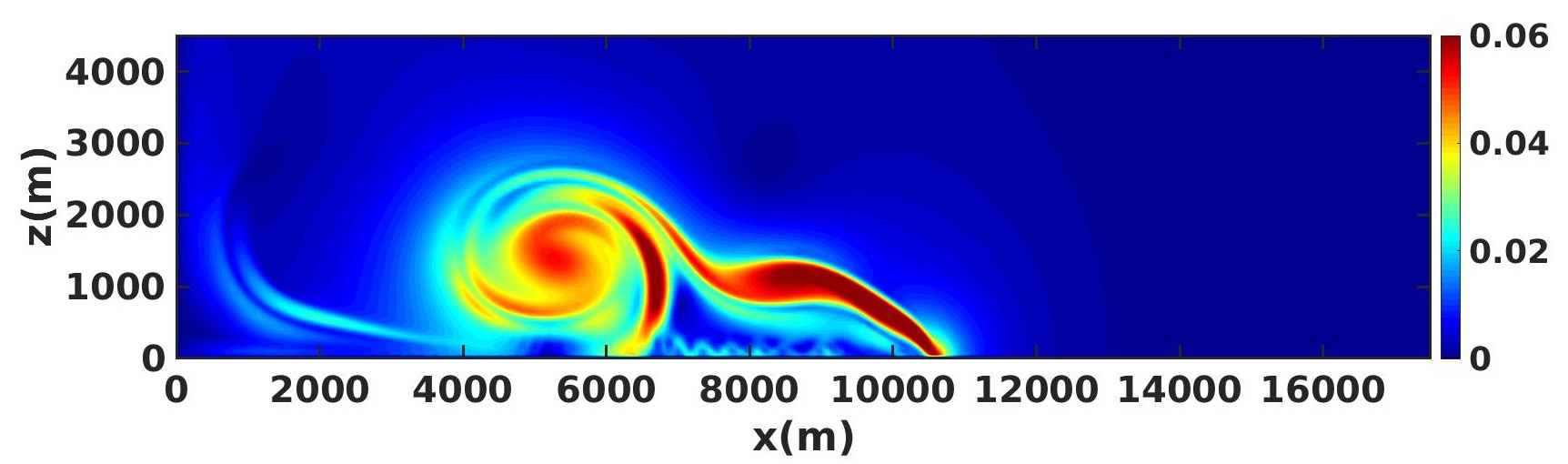}
    \put(115,50){\textcolor{white}{\footnotesize{$a_S$, $\alpha = 11$, $t = 600$ s}}}
    \end{overpic}
    \begin{overpic}[width=0.48\textwidth]{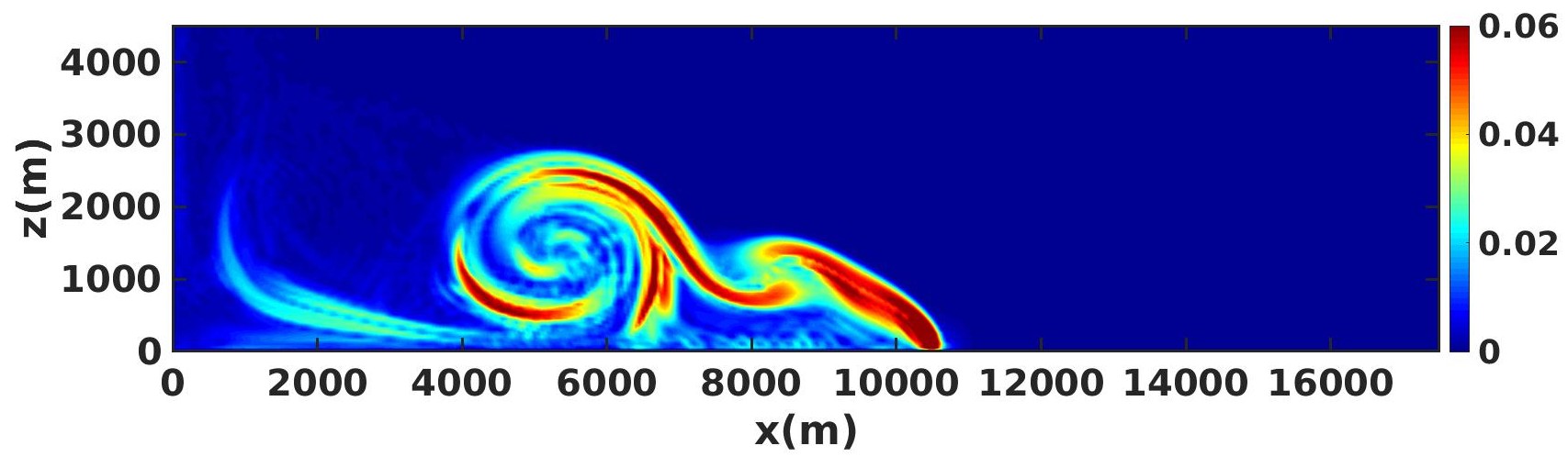}
    \put(115,50){\textcolor{white}{\footnotesize{$a_D$, $\alpha = 12$, $t = 600$ s}}}
    \end{overpic} \\ \vspace{0.2cm}
    \begin{overpic}[width=0.48\textwidth]{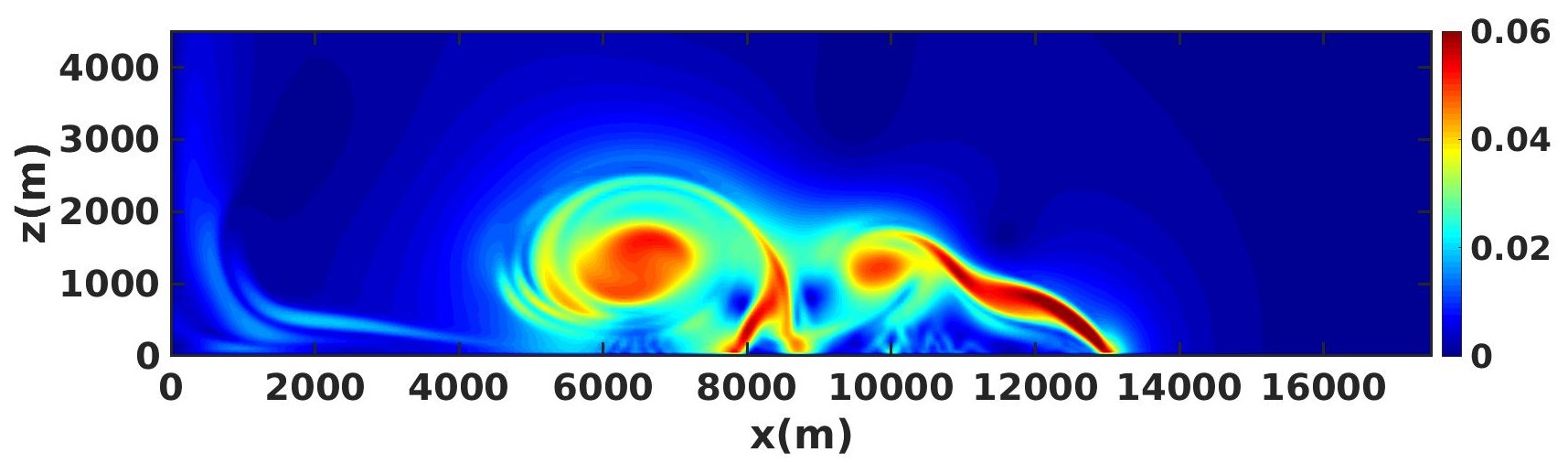}
    \put(115,50){\textcolor{white}{\footnotesize{$a_S$, $\alpha = 11$, $t = 750$ s}}}
    \end{overpic}
    \begin{overpic}[width=0.48\textwidth]{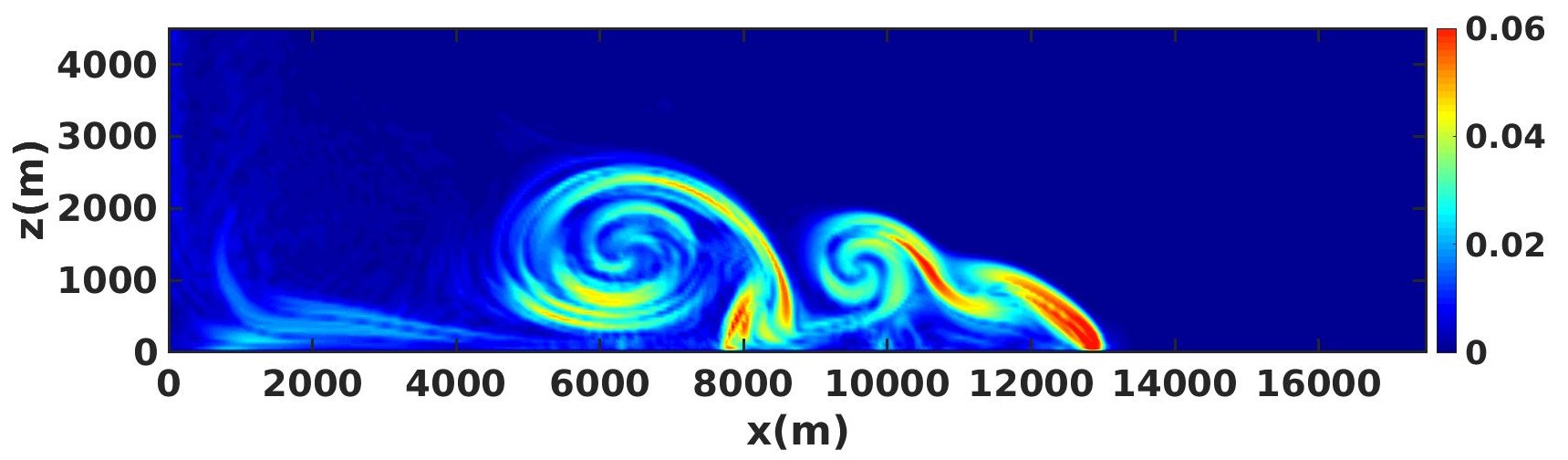}
    \put(115,50){\textcolor{white}{\footnotesize{$a_D$, $\alpha = 12$, $t = 750$ s}}}
    \end{overpic} \\ \vspace{0.2cm}
    \begin{overpic}[width=0.48\textwidth]{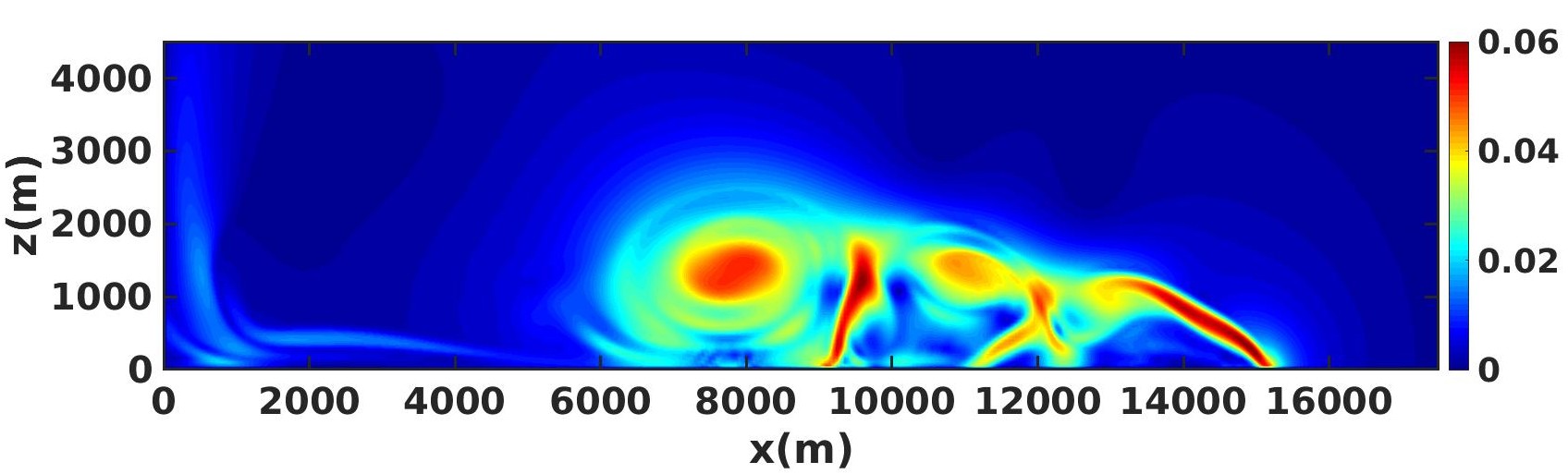}
    \put(115,50){\textcolor{white}{\footnotesize{$a_S$, $\alpha = 11$, $t = 900$ s}}}
    \end{overpic}
    \begin{overpic}[width=0.48\textwidth]{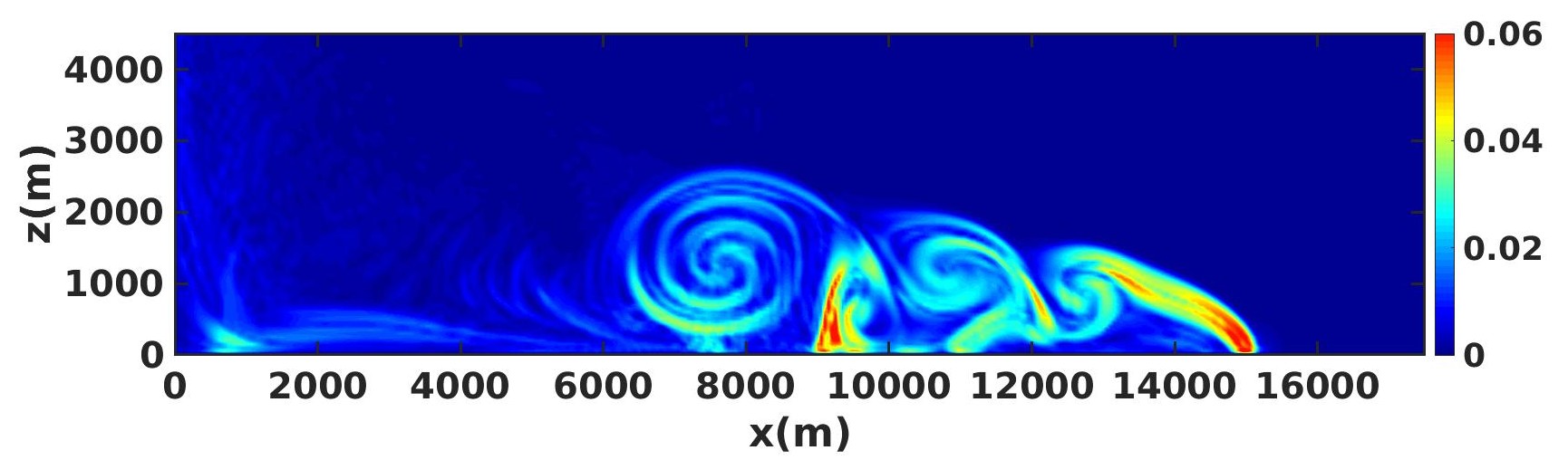}
    \put(115,50){\textcolor{white}{\footnotesize{$a_D$, $\alpha = 12$, $t = 900$ s}}}
    \end{overpic} \\ \vspace{0.2cm}
    \caption{
    Density current, mesh $h = 50$ m:
time evolution of $a_S$ with $\alpha=11$ m (left) and $a_D$ with $\alpha=12$ m (right).
    }
    \label{fig::dc_12.5_aU}
\end{figure}

Next, in Fig.~\ref{fig::dc_12.5_aU_1} we report with a comparison of indicator functions $a_S$ and $a_D$ for the simulations in Fig.~\ref{fig::SmagoRayModel125_25} (left) and \ref{fig:deconRayModel125_25} (left).
We see see that, while both indicator functions have larger values (red to yellow shades) at the bottom of the largest recirculation, at a given time $a_S$ has larger regions of intermediate values (light blue shade) than $a_D$. 
This is due to the fact that $a_S$ is a less selective indicator function, as mentioned earlier. This is more evident on mesh $h = 50$ m: see Fig.~\ref{fig::dc_12.5_aU} for the plots of 
the indicator function for $a_S$ and $a_D$ for the simulations in Fig.~\ref{fig::SmagoRayModel50} (left) and
\ref{fig:deconRayModel50_12}. The higher selectivity of $a_D$ results in much smaller regions of high and intermediate values (red to green shades). Finally, we note that the maximum magnitude of the indicator function is higher for the finer mesh ($h = 25$ m), which gives rise to more localized and higher peaks than the coarser mesh ($h = 50$ m). This is in line with what already observed in \cite{BQV}.

We conclude with a comment on the computational cost. Table \ref{tab:6} reports the computational time taken by the evolve step and filter step per time step and total simulation time for the EFR algorithm with indicator functions $a_L$, $a_S$, and $a_D$ and the specified values of $\alpha$ for meshes $h=50,25$ m. All the simulations were run on a common laptop (AMD Ryzen 7 5700U, 16GB RAM). As expected, the total computational cost increases when switching from the linear filter to nonlinear filters, with the deconvolution-based indicator function being the most expensive. In fact, while $a_S$ requires a simple post-processing of the velocity field, $a_D$ in \eqref{eq:a_D0_a_D1} requires one application
of the linear Helmholtz filter.
Despite this increased cost, the simulation with mesh $h = 25$ m and $a_D$ takes a little less than 1 hour and 10 minutes, which means that the solver is rather efficient. We note that when using
$a_D$, the filter step takes about half of the time needed for the evolve step. 
This might seem counter-intuitive given the relative complexity of the problems solved at the two steps. However, it can be explained with the solver choices.
During the evolve step, the majority of the computational cost is spent to solve the equation for $p'$ (recall we adopt a splitting scheme detailed in \cite{GQR_OF_clima}) with the Diagonal incomplete Cholesky preconditioned conjugate gradient method, which is the same method used for the filter step. Since the the equation of mass conservation is treated fully explicitly, it is very inexpensive to solve. The solver for the conservation of energy equation uses the bi-conjugate gradient stabilized method with a diagonal-based incomplete LU preconditioner. 
The accuracy for the resolution of all the linear system is set to $1e-8$.
The computational cost of the evolve step is also contained by not performing a momentum predictor step.

\begin{table}[htb!]
\begin{center}
\begin{tabular}{ | c | c | c | c | c | c | }
\hline
Model & $h$ (m) & $\alpha$ (m) & Evolve (s) & Filter (s) & Total (s)\\

 \hline
 EFR, $a_L$ & 25    & 2.7  & 0.3 &  0.070  & 3492 \\
 \hline
EFR, $a_S$ & 25    & 8  & 0.3 &  0.108  & 3880 \\
\hline
EFR, $a_D$   & 25  & 10  & 0.3 &  0.159  & 4163 \\
\hline
EFR, $a_L$ & 50  & 2.7  & 0.06 &  0.015 & 707\\
\hline
EFR, $a_S$ & 50  & 11  & 0.06 &  0.023  & 772\\
\hline
EFR, $a_D$ & 50  & 12  & 0.06 &  0.029  & 790\\
\hline 
\end{tabular}
\caption{Density current:
computational time taken by the evolve step and filter step per time step and total simulation time for the EFR algorithm with indicator functions $a_L$, $a_S$, and $a_D$ and the specified values of $\alpha$ for meshes $h=50,25$ m.}\label{tab:6}
\end{center}
\end{table}

\section{Concluding remarks}\label{sec:concl}
In this paper, we presented a filter stabilization technique for the mildly compressible Euler equations that is realized through a three step algorithm called Evolve-Filter-Relax (EFR). 
While filter stabilization and the EFR algorithm have been widely investigated for the incompressible Navier-Stokes equations, this work is one of the few papers that applies them to the 
Euler equations. We showed that the EFR algorithm is equivalent to an eddy viscosity model in LES and we considered three indicator functions to tune the amount and location of eddy viscosity: a constant indicator function, an indicator function proportional to the velocity gradient norm that recovers a Smagorinsky-like model, and indicator function based on approximate deconvolution operators. The first indicator function corresponds to a linear filter, which is known to be overdiffusive, while the other two lead to nonlinear filter. 

We tested our EFR approach with two well-known benchmarks for atmospheric flow: the rising thermal bubble and the density current. For both benchmarks, we showed that the linear filter provides results in excellent agreement with data in the literature obtained by setting an ah-hoc eddy viscosity. We showed that with the nonlinear filters we can capture a larger amount of vortical structures in the flows, which is expected and in line with other publshed data obtained with LES models. Of the three indicator functions under consideration, the deconvolution-based indicator function was shown to be more selective in identifying the regions of the domain where artificial diffusion is needed. Finally, we
commented about the computational efficiency of our approach, highlighting that the filter step is computationally cheap with respect to the evolve step.  


More work is needed to improve the EFR algorithm proposed in this paper. 
A parametric study of the relaxation parameters $\chi$ and $\xi$ (set equal to 1 in this work) would inform us on their ``optimal'' value to improve the accuracy and would help to mitigate the sensitivity to the filtering radius $\alpha$. In addition, the role of the order of the deconvolution $N$ needs to be investigated. 

\section*{Acknowledgements}
We acknowledge the support provided by the European Research Council Executive Agency by the Consolidator Grant project AROMA-CFD ``Advanced Reduced Order Methods with Applications in Computational Fluid Dynamics" - GA 681447, H2020-ERC CoG 2015 AROMA-CFD, PI G. Rozza, and INdAM-GNCS 2019-2020 projects.
This work was also partially supported by US National Science Foundation through grant DMS-1953535 (PI A.~Quaini). A.~Quaini acknowledges 
support from the Radcliffe Institute for Advanced Study at Harvard University where she has been the 2021-2022 William and Flora Hewlett Foundation Fellow.

\bibliographystyle{plain}
\bibliography{bibliography.bib}

\begin{thebibliography}{10}

\bibitem{ABGRALL2001277}
R.~Abgrall.
\newblock Toward the ultimate conservative scheme: Following the quest.
\newblock {\em Journal of Computational Physics}, 167(2):277--315, 2001.

\bibitem{ahmadLindeman2007}
N.~{Ahmad} and J.~{Lindeman}.
\newblock Euler solutions using flux-based wave decomposition.
\newblock {\em Int. J. Numer. Meth. Fluids}, 54:47--72, 2007.

\bibitem{ahmad2018}
N.~N. {Ahmad}.
\newblock High-resolution wave propagation method for stratified flows.
\newblock In {\em {AIAA} Aviation Forum, Atlanta, GA. AIAA}, 2018.

\bibitem{bazilevsCaloCottrellHughesRealiScovazzi2007}
Y.~Bazilevs, V.~Calo, J.~A. Cottrell, T.~J.~R. Hughes, A.~Reali, and
  G.~Scovazzi.
\newblock Variational multiscale residual-based turbulence modeling for large
  eddy simulation of incompressible flows.
\newblock {\em Comput. Methods Appl. Mech. Engrg.}, 197:173--201, 2007.

\bibitem{BQV}
L.~Bertagna, A.~Quaini, and A.~Veneziani.
\newblock {Deconvolution-based nonlinear filtering for incompressible flows at
  moderately large Reynolds numbers}.
\newblock {\em International Journal for Numerical Methods in Fluids},
  81(8):463--488, 2016.

\bibitem{Borggaard2009}
J.~Borggaard, T.~Iliescu, and J.P. Roop.
\newblock A bounded artificial viscosity large eddy simulation model.
\newblock {\em SIAM Journal on Numerical Analysis}, 47:622--645, 2009.

\bibitem{Bowers2012}
A.~L. Bowers, L.~G. Rebholz, A.~Takhirov, and C.~Trenchea.
\newblock Improved accuracy in regularization models of incompressible flow via
  adaptive nonlinear filtering.
\newblock {\em International Journal for Numerical Methods in Fluids},
  70(7):805--828, 2012.

\bibitem{abigail_CMAME}
A.L. Bowers and L.G. Rebholz.
\newblock Numerical study of a regularization model for incompressible flow
  with deconvolution-based adaptive nonlinear filtering.
\newblock {\em Comput. Methods Appl. Mech. Eng.}, 258:1--12, 2013.

\bibitem{Boyd1998283}
J.P. Boyd.
\newblock Two comments on filtering (artificial viscosity) for {Chebyshev} and
  {Legendre} spectral and spectral element methods: Preserving boundary
  conditions and interpretation of the filter as a diffusion.
\newblock {\em Journal of Computational Physics}, 143(1):283 -- 288, 1998.

\bibitem{carpenterDroegemeier1990}
R.~{Carpenter}, K.~{Droegemeier}, P.~{Woodward}, and C.~{Hane}.
\newblock Application of the piecewise parabolic method ({PPM}) to
  meteorological modeling.
\newblock {\em Mon. Wea. Rev.}, 118:586--612, 1990.

\bibitem{Chehab2021}
J.-P. Chehab.
\newblock Damping, stabilization and numerical filtering for the modeling and
  the simulation of time dependent {PDEs}.
\newblock {\em Discrete \& Continuous Dynamical Systems - Series S},
  14(8):2693--2728, 2021.

\bibitem{codina2002}
R.~{Codina}.
\newblock Stabilized finite element approximation of transient incompressible
  flows using orthogonal subscales.
\newblock {\em Comput. Methods Appl. Mech. Engrg.}, 191:4295--4321, 2002.

\bibitem{Codinaetal2017}
Ramon Codina, Santiago Badia, Joan Baiges, and Javier Principe.
\newblock {\em Variational Multiscale Methods in Computational Fluid Dynamics},
  pages 1--28.
\newblock John Wiley \& Sons, Ltd, 2017.

\bibitem{Dunca2005}
A.~Dunca and Y.~Epshteyn.
\newblock On the {Stolz-Adams} deconvolution model for the large-eddy
  simulation of turbulent flows.
\newblock {\em SIAM Journal on Mathematical Analysis}, 37(6):1890--1902, 2005.

\bibitem{doi:10.1137/120867482}
Alexandre Ern and Jean-Luc Guermond.
\newblock Weighting the edge stabilization.
\newblock {\em SIAM Journal on Numerical Analysis}, 51(3):1655--1677, 2013.

\bibitem{Ervin2012}
V.J. Ervin, W.J. Layton, and M.~Neda.
\newblock Numerical analysis of filter-based stabilization for evolution
  equations.
\newblock {\em SIAM Journal on Numerical Analysis}, 50(5):2307--2335, 2012.

\bibitem{Feng2021}
Y.~Feng, J.~Miranda-Fuentes, J.~Jacob, and P.~Sagaut.
\newblock Hybrid lattice boltzmann model for atmospheric flows under anelastic
  approximation.
\newblock {\em Physics of Fluids}, 33(3):036607, 2021.

\bibitem{Fischer2001265}
P.~Fischer and J.~Mullen.
\newblock Filter-based stabilization of spectral element methods.
\newblock {\em Comptes Rendus de l'Academie des Sciences - Series I -
  Mathematics}, 332(3):265 -- 270, 2001.

\bibitem{B-garnier}
E.~Garnier, N.~Adams, and P.~Sagaut.
\newblock {\em Large Eddy Simulation for Compressible Flows}.
\newblock Springer, Berlin, 2009.

\bibitem{GEA}
{GEA - Geophysical and Environmental Applications}.
\newblock \url{https://github.com/GEA-Geophysical-and-Environmental-Apps/GEA}.

\bibitem{giraldo_2008}
F.~X. Giraldo and M.~Restelli.
\newblock {A study of spectral element and discontinuous Galerkin methods for
  the Navier-Stokes equations in nonhydrostatic mesoscale atmospheric modeling:
  Equation sets and test cases}.
\newblock {\em J. Comput. Phys.}, 227:3849--3877, 2008.

\bibitem{giraldoRestelli2008b}
F.~X. {Giraldo} and M.~{Restelli}.
\newblock A conservative discontinuous galerkin semi-implicit formulation for
  the navier-stokes equations in nonhydrostatic mesoscale modeling.
\newblock {\em SIAM J. Sci. Comp.}, 31:2231--2257, 2009.

\bibitem{GIRFOGLIO201927}
Michele Girfoglio, Annalisa Quaini, and Gianluigi Rozza.
\newblock A finite volume approximation of the navier-stokes equations with
  nonlinear filtering stabilization.
\newblock {\em Computers \& Fluids}, 187:27--45, 2019.

\bibitem{GQR_GEA}
Michele Girfoglio, Annalisa Quaini, and Gianluigi Rozza.
\newblock {GEA: a new finite volume-based open source code for the numerical
  simulation of atmospheric and ocean flows}.
\newblock {\em https://arxiv.org/abs/2303.10499}, 2023.

\bibitem{GQR_OF_clima}
Michele Girfoglio, Annalisa Quaini, and Gianluigi Rozza.
\newblock {Validation of an OpenFOAM\textsuperscript{\textregistered}-based
  solver for the Euler equations with benchmarks for mesoscale atmospheric
  modeling}.
\newblock {\em https://arxiv.org/abs/2302.04836}, 2023.

\bibitem{Guermond_pasquetti_popov_JCP_2011}
J~L. Guermond, R.~Pasqueti, and B.~Popov.
\newblock Entropy viscosity method for nonlinear conservation laws.
\newblock {\em J. Comput. Phys.}, 230(11):4248--4267, 2011.

\bibitem{guermondPasquetti2008}
J~L. Guermond and R.~Pasquetti.
\newblock Entropy-based nonlinear viscosity for {Fourier} approximations of
  conservation laws.
\newblock {\em C. R. Acad. Sci., Ser. I}, 346:801--806, 2008.

\bibitem{Guermond_Popov_2014}
J~L. Guermond and B.~Popov.
\newblock Viscous regularization of the {E}uler equations and entropy
  principles.
\newblock {\em SIAM J. Appl. Math.}, 74(2):284--305, 2014.

\bibitem{Hesthaven_Warburton}
Jan~S. Hesthaven and Tim Warburton.
\newblock {\em Nodal Discontinuous Galerkin Methods: Algorithms, Analysis, and
  Applications}.
\newblock Springer Publishing Company, Incorporated, 2007.

\bibitem{Euler_comp_orig}
D.~D. Holm.
\newblock {Averaged Lagrangians and the mean effects of fluctuations in ideal
  fluid dynamics}.
\newblock {\em Physica D: Nonlinear Phenomena}, 170:253--286, 2002.

\bibitem{hughes1995}
T.~{Hughes}.
\newblock Multiscale phenomena: Green's functions, the {D}irichlet-to-{N}eumann
  formulation, subgrid scale models, bubbles and the origins of stabilized
  methods.
\newblock {\em Comput. Methods Appl. Mech. and Engrg.}, 127:387--401, 1995.

\bibitem{hughesFeijoo1998}
T~J~R Hughes, G.~{Feij\'oo}, L.~{Mazzei}, and J.~{Quincy}.
\newblock The variational multiscale method -- {A} paradigm for computational
  mechanics.
\newblock {\em Comput. Methods Appl. Mech. Engrg.}, 166:3--24, 1998.

\bibitem{O-hunt1988}
J.C. Hunt, A.A. Wray, and P.~Moin.
\newblock Eddies stream and convergence zones in turbulent flows.
\newblock Technical Report CTR-S88, CTR report, 1988.

\bibitem{jasakphd}
H.~Jasak.
\newblock {\em Error analysis and estimation for the finite volume method with
  applications to fluid flows}.
\newblock PhD thesis, Imperial College, University of London, 1996.

\bibitem{kellyGiraldo2012}
J.~F. {Kelly} and F.~X. {Giraldo}.
\newblock Continuous and discontinuous {G}alerkin methods for a scalable
  three-dimensional nonhydrostatic atmospheric model: limited-area mode.
\newblock {\em J. Comput. Phys.}, 231:7988--8008, 2012.

\bibitem{klockner_warburton_hesthaven_2011}
A.~Kloeckner, T.~Warburton, and J.~S. Hesthaven.
\newblock Viscous shock capturing in a time-explicit discontinuous galerkin
  method.
\newblock {\em Mathematical Modelling of Natural Phenomena}, 6(3):57–83,
  2011.

\bibitem{KURGANOV20128114}
Alexander Kurganov and Yu~Liu.
\newblock New adaptive artificial viscosity method for hyperbolic systems of
  conservation laws.
\newblock {\em Journal of Computational Physics}, 231(24):8114--8132, 2012.

\bibitem{layton_JMFM}
W.~Layton, L.G. Rebholz, and C.~Trenchea.
\newblock Modular nonlinear filter stabilization of methods for higher
  {Reynolds} numbers flow.
\newblock {\em Journal of Mathematical Fluid Mechanics}, 14:325--354, 2012.

\bibitem{LAYTON20113183}
W.~Layton, L.~R{\"o}he, and H.~Tran.
\newblock Explicitly uncoupled vms stabilization of fluid flow.
\newblock {\em Computer Methods in Applied Mechanics and Engineering},
  200(45):3183--3199, 2011.

\bibitem{marrasEtAl2013a}
S.~{Marras}, M.~{Moragues}, M.~{V\'azquez}, O.~{Jorba}, and G.~{Houzeaux}.
\newblock A {V}ariational {M}ultiscale {S}tabilized finite element method for
  the solution of the {E}uler equations of nonhydrostatic stratified flows.
\newblock {\em J. Comput. Phys.}, 236:380--407, 2013.

\bibitem{marrasNazarovGiraldo2015}
S.~{Marras}, M.~{Nazarov}, and F.~X. {Giraldo}.
\newblock {Stabilized high-order Galerkin methods based on a parameter-free
  dynamic SGS model for LES}.
\newblock {\em J. Comput. Phys.}, 301:77--101, 2015.

\bibitem{Mathew2003}
J.~Mathew, R.~Lechner, H.~Foysi, J.~Sesterhenn, and R.~Friedrich.
\newblock An explicit filtering method for large eddy simulation of
  compressible flows.
\newblock {\em Physics of Fluids}, 15(8):2279--2289, 2003.

\bibitem{CNM:CNM219}
J.~Mullen and P.~Fischer.
\newblock Filtering techniques for complex geometry fluid flows.
\newblock {\em Communications in Numerical Methods in Engineering},
  15(1):9--18, 1999.

\bibitem{nazarovHoffman2013}
M.~{Nazarov} and J.~{Hoffman}.
\newblock Residual-based artificial viscosity for simulation of turbulent
  compressible flow using adaptive finite element methods.
\newblock {\em Int. J. Numer. Methods Fluids}, 71:339--357, 2013.

\bibitem{Olshanskii2013}
M.A. Olshanskii and X.~Xiong.
\newblock A connection between filter stabilization and eddy viscosity models.
\newblock {\em Numerical Methods for Partial Differential Equations},
  29(6):2061--2080, 2013.

\bibitem{perssonPeraire2006}
P.-O. {Persson} and J.~{Peraire}.
\newblock Sub-cell shock capturing for discontinuous {Galerkin} methods.
\newblock {\em Proc. of the 44th AIAA Aerospace Sciences Meeting and Exhibit},
  AIAA-2006-112, 2006.

\bibitem{rispoliSaavedra2006}
F.~{Rispoli} and R.~{Saavedra}.
\newblock A stabilized finite element method based on sgs models for
  compressible flows.
\newblock {\em Comp. Meth. Appl. Mech. Engrg.}, 196:652--664, 2006.

\bibitem{Euler_comp}
P.~Secchi.
\newblock {An alpha model for compressible fluids}.
\newblock {\em Discrete and Continuous Dynamical System - S}, 3:351--359, 2008.

\bibitem{smagorinsky1963}
J.~{Smagorinsky}.
\newblock General circulation experiments with the primitive equations: {I}.
  the basic experiement.
\newblock {\em Mon. Wea. Rev.}, 91:99--164, 1963.

\bibitem{Stolz1999}
S.~Stolz and N.A. Adams.
\newblock An approximate deconvolution procedure for large-eddy simulation.
\newblock {\em Physics of Fluids}, 11(7):1699--1701, 1999.

\bibitem{Stolz2001}
S.~Stolz, N.A. Adams, and L.~Kleiser.
\newblock An approximate deconvolution model for large-eddy simulation with
  application to incompressible wall-bounded flows.
\newblock {\em Physics of Fluids}, 13(4):997--1015, 2001.

\bibitem{strakaWilhelmson1993}
J.~{Straka}, R.~{Wilhelmson}, L.~{Wicker}, J.~{Anderson}, and K.~{Droegemeier}.
\newblock Numerical solution of a nonlinear density current: a benchmark
  solution and comparisons.
\newblock {\em Int. J. Num. Meth. in Fluids}, 17:1--22, 1993.

\bibitem{Visbal2002}
M.R. Visbal and D.P. Rizzetta.
\newblock Large eddy simulation on curvilinear grids using compact differencing
  and filtering schemes.
\newblock {\em J. Fluids Eng.}, 124:836--847, 2002.

\bibitem{Weller1998}
H.~G. Weller, G.~Tabor, H.~Jasak, and C.~Fureby.
\newblock A tensorial approach to computational continuum mechanics using
  object-oriented techniques.
\newblock {\em Computers in physics}, 12(6):620--631, 1998.

\end{thebibliography}
\end{document}